\documentclass[twoside,11pt,a4paper]{article}
\usepackage{amsmath,pifont,amsfonts,amssymb,latexsym, wasysym}
\usepackage{a4wide,xspace}
\pdfoutput=1
\usepackage{bbm}
\usepackage{framed}
\usepackage{graphicx}
\usepackage{epic,color}      
\usepackage{theorem}
\usepackage{fancyhdr,enumerate}
\usepackage[colorlinks=true, linkcolor=blue, citecolor=darkgreen]{hyperref}
\renewcommand{\arraystretch}{1.28}
\newtheorem{Theorem}{Theorem}
\newtheorem{Lemma}{Lemma}[section]

\newtheorem{Corollary}[Lemma]{Corollary}
\newtheorem{Definition}[Lemma]{Definition}
\newtheorem{Example}[Lemma]{Example}
\newtheorem{Proposition}[Lemma]{Proposition}

\newtheorem{Remark}[Lemma]{Remark}
\newcommand{\point}{\hspace{-0.2cm}{\bf . }}
\newcommand{\pointpar}{\point\ \par}
\newcommand{\li}{\left}
\newcommand{\re}{\right}
\begin{document}
\definecolor{darkgreen}{rgb}{0,0.65,0}
%%%%%%%%%%%%%%%%%%%%%%%%%%%%%%%%%%%%%%%%%%%%%%%%%%%%%%%%%%%%%%%%%%%%%%%%%%%%%%%%%%%%%%%%%%%%%%%%%%%     Haeder

\newpage
\title{{\huge The historical Moran model}}
\author{{\Large Peter Seidel}\footnote{Friedrich-Alexander-Universit\"at Erlangen-N\"urnberg, Department Mathematik, Cauerstr. $11$, $91058$ Erlangen, Germany, seidel@math.fau.de} \footnote{The author was partly supported by the DFG Priority Programme SPP 1590 ''Probabilistic Structures in Evolution''.}}
\date{{\Large \today}}
\maketitle
\bigskip
%%%%%%%%%%%%%%%%%%%%%%%%%%%%%%%%%%%%%%%%%%%%%%%%%%%%%%%%%%%%%%%%%%%%%%%%%%%%%%%%%%%%%%%%%%%%%%%%%%% Abstract

\begin{abstract}
\noindent  We consider a multi-type Moran model (in continuous time) with selection and type-dependent mutation. This paper is concerned with the evolution of genealogical information forward in time. For this purpose we define and analytically characterize a path-valued Markov process that contains in its state at time $t$ the extended ancestral lines (adding genealogical distances) of the population alive at time $t$.

The main result is a representation for the conditional distribution of the extended ancestral lines of a subpopulation alive at a fixed time $T$ (present time) given the type information of the subpopulation at time $T$ in terms of the distribution of the sample paths (up to time $T$) of a special Markov process (different from the ancestral selection graph) to which we refer as backward process. This representation allows us both to prove that the extended ancestral lines converge in the limit $T \to \infty$ if the type information converges in the limit $T \to \infty$ and to study the resulting limit of the extended ancestral lines by means of the backward process.

The limit theorem has two applications: First, we can represent the stationary type distribution of the common ancestor type process in terms of the equilibrium distribution of a functional of the backward process, where in the two type case we recover the common ancestor process of Fearnhead (\cite{F02}) if we let the population size tend to infinity. Second, we obtain that the conditioned genealogical distance of two individuals given the types of the two individuals is distributed as a certain stopping time of a further functional of the backward process which is a new approach towards a proof that genealogical distances are stochastically smaller under selection. 
\end{abstract}



%%%%%%%%%%%%%%%%%%%%%%%%%%%%%%%%%%%%%%%%%%%%%%%%%%%%%%%%%%%%%%%%%%%%%%%%%%%%%%%%%%%%%%%%%%%%%%%%%%%      Introduction, the model

\newpage
\tableofcontents
\newpage
\section{Introduction}
\noindent The evolution of a population is commonly described by a Markov process on a suitable state space (for example, see 
\cite{D93}, \cite{DK96}, \cite{DK99}) whose complexity increases with the complexity of the information that is considered. A fundamental question in population genetics which recently undergoes rapid development is how to code {\it genealogical information} in order to study its evolution in time. In this article we approach this complex issue. The innovation concerns the development of a {\it new} representation for genealogical information which allows to model its evolution in time {\it by a Markov process} on the one hand and to investigate {\it fixed time genealogies} in a manageable way on the other hand.

In this paper we consider a version of the {\it Moran Model} in continuous time. This is a model of a population of fixed size $N \in \mathbb{N}$ in which each individual has a {\it genetic type} that determines the {\it fitness} of the individual, where we allow to have $d \in \mathbb{N}$ different genetic types and more than two levels of fitness. The evolution of the population is given by {\it type-dependent mutation} and a change of generation mechanism ({\it resampling}) which is affected by the fitness of the involved individuals ({\it selection}). This means:
\begin{enumerate}
\item \label{itemMutation} {\bf (Mutation)} The type of each individual evolves as a continuous time Markov chain according to a general stochastic matrix.
\item \label{itemResampling} {\bf (Resampling)} Each pair of individuals dies after an exponential holding time and gives birth to a new pair of individuals, called descendants, which choose an  ancestor from the dying pair from which both descendants inherit the type. 

{\bf (Selection)} The choice of the ancestor is not purely random but there is a bias depending on the fitness of the individuals in the dying pair.
\end{enumerate}
%%%%%%%%%%%%%%%%%%%%%%%%%%%%%%%%%%%%%%%%%%%%%%%%%%%%%%%%%%%%%%%%%%%%%%%%%%%%%%%%%%%%%%%%%%%%%%%%%%%      Introduction, diffusion limit

In order to describe the evolution of the type information it is convenient to describe the relative type frequencies of the population alive at time $t$ by a probability measure on the type space. The large population limit ($N \to \infty$) of this measure-valued Markov process leads to the {\it measure-valued Fleming-Viot process} (to the {\it Wright-Fisher diffusion} if in the two type case one considers the relative type frequency of one of the two types). \bigskip
%%%%%%%%%%%%%%%%%%%%%%%%%%%%%%%%%%%%%%%%%%%%%%%%%%%%%%%%%%%%%%%%%%%%%%%%%%%%%%%%%%%%%%%%%%%%%%%%%%%      Introduction, our main goal

The present article, which grew out of \cite{Seidel15},
is devoted to the {\it evolution} of genealogical information {\it forward in time}, in particular, the evolution of both {\it ancestral lines} (coding types of ancestors through time) and {\it genealogical distances}. For this purpose we define a Markov process, which we call the {\it historical Moran model} (short {\bf HMM}), that {\it contains in its state} at time $t$ the {\it extended ancestral lines} of the population alive at time $t$ and can analytically be characterized by means of a well-posed martingale problem. 
%%%%%%%%%%%%%%%%%%%%%%%%%%%%%%%%%%%%%%%%%%%%%%%%%%%%%%%%%%%%%%%%%%%%%%%%%%%%%%%%%%%%%%%%%%%%%%%%%%%      Introduction, background

Our motivation goes back to \cite{DP91} ({\it historical processes}) and \cite{GLW05} ({\it historical interacting Wright-Fisher diffusions}) on the one hand and to \cite{GPW13} and \cite{DGP12} ({\it tree-valued Fleming-Viot process}) on the other hand, but note that the {\it new concept} of extended ancestral lines is different from the concept of a {\it marked metric measure space}. In particular, since we do not pass to equivalence classes (exchangeability property), the extended ancestral lines contain {\it more} information even if the metric measure space is  {\it path-marked} (see \cite{GSW} and \cite{gswtop}). What is important here is that the HMM includes the {\it tree-valued Moran model} as a functional, and hence (see \cite{DGP12}, Theorem 3) the genealogical distance of two individuals sampled from the tree-valued Fleming-Viot process can by analyzed if we let $N \to \infty$ in our HMM. However, this paper does not treat a diffusion limit theorem for the HMM. 

We emphasize that our approach is different from the {\it ancestral selection graph} (\cite{KN97} and \cite{NK97}) and from the {\it look-down construction} (\cite{DK96} and \cite{DK99}) in the following respect. Namely, we give a precise definition of the ancestral lines of the population alive at time $t$ and how these ancestral lines evolve {\it forwards} in time, and hence it is  possible to define the so-called (see \cite{F02}, \cite{Tay07}, \cite{KHB13} and \cite{LKBW15}) {\it common ancestor type process} (short {\bf CAT}) for the Moran model {\it rigorously} as a functional of a {\it forward} evolving population.
%%%%%%%%%%%%%%%%%%%%%%%%%%%%%%%%%%%%%%%%%%%%%%%%%%%%%%%%%%%%%%%%%%%%%%%%%%%%%%%%%%%%%%%%%%%%%%%%%%%      Introduction, main result

The main result in the present paper is a new approach to investigate the distribution of genealogical information of a subpopulation alive at a fixed time $T$, an important issue in population genetics. Namely, the analytic characterization of the HMM allows us to state and prove a {\it strong stochastic representation} for the {\it conditioned} extended ancestral lines. This means, we express the conditional distribution of the extended ancestral lines of a subpopulation alive at time $T$ {\it given the type information of this subpopulation at time $T$} in terms of the  distribution of the sample paths (up to time $T$) of a suitable Markov process which starts with the type information of the subpopulation in the HMM at time $T$ and generates the genealogical information {\it backwards in time}. Note that we use the notation "strong" to emphasize that we represent the distribution of a random variable (in our case conditioned) and not only the expectation of certain functionals like moments, and the notation "stochastic representation" to emphasize that we represent the expectation of certain functionals, but not in the form of a ({\it Feynman-Kac}) {\it duality} (see \cite{JK14}, \cite{DG13} and \cite{EK86}). Furthermore observe that the conditional distribution of genealogical information at time $T$ (given the type information at time $T$) has not been investigated in this general form before (except in \cite{DSt03}, see also \cite{EG09}, but only for the  case where the {\it type information of the Moran model is in equilibrium}, and once again without an explicit definition of genealogical information).  
%%%%%%%%%%%%%%%%%%%%%%%%%%%%%%%%%%%%%%%%%%%%%%%%%%%%%%%%%%%%%%%%%%%%%%%%%%%%%%%%%%%%%%%%%%%%%%%%%%%      Introduction, furhter results and applications

The strong stochastic representation at time $T$ can be used to obtain a limit distribution for the  extended ancestral lines in the limit $T \to \infty$ which is the foundation to investigate the {\it stationary type distribution} of the CAT (in \cite{LKBW15} the notation {\it common ancestor type distribution} is used) and genealogical distances for the multi-type Moran model with population size $N$, selection and type-dependent mutation in a mathematically rigorous way. In order to show that we can already do explicit calculations with our machinery, we consider the case where we have two types and we also let $N \to \infty$. 
\begin{itemize}
\item  On the one hand we recover the common ancestor process of Fearnhead (in \cite{KHB13} this object is called {\it pruned ancestral selection graph}) and therefore we get the stationary type distribution of the CAT for the Wright-Fisher diffusion.
\item On the other hand we make a first step (compare with \cite{DGP12}, Section 3.6$.$) in showing that for {\it any selection strength} genealogical distances (in the tree-valued Fleming-Viot process in equilibrium) are {\it stochastically smaller} under selection. 
\end{itemize}
Note that the present article only contains a brief discussion of the case with more than two types which will be the subject of future research.\bigskip
%%%%%%%%%%%%%%%%%%%%%%%%%%%%%%%%%%%%%%%%%%%%%%%%%%%%%%%%%%%%%%%%%%%%%%%%%%%%%%%%%%%%%%%%%%%%%%%%%%%      Introduction, life-sites and extended ancestral lines

Now we explain the {\it important new concept} of extended ancestral lines that allows us to consider {\it both} ancestral lines and genealogical distances. 

Let
\begin{equation}
I = \{1, \dots, N\}
\end{equation}
and
\begin{equation}
K = \{0,1, \dots,d-1 \}
\end{equation}
where the set $I$, which we call the set of {\it life-sites}, describes the population of size $N$ and the set $K$ models $d$ different genetic types. When we speak of a life-site $i \in I$ at time $t$ (or simply $i \in I$ at time $t$) in what follows then we mean the individual that occupies the life-site $i$ at time $t$. 

The extended ancestral line of $i \in I$ at time $t$ is an element  in
\begin{equation}
\mathcal{D}([0,t], K \times I) \;\;\;\;\mbox{ \bf (cadlag paths from $[0,t]$ to $K \times I$)}
\end{equation}
that describes the genealogy of $i \in I$ at time $t$, namely, evaluated at time $s \in [0,t]$ this extended ancestral line gives us for $i$ besides the type also the life-site of the ancestor alive at time $s$. Hence the whole collection of {\it extended ancestral lines} at time $t$ obviously includes the information of the whole collection of ancestral lines at time $t$. In addition, we can consider genealogical distances. Namely, half the genealogical distance of the life-sites $i$ and $j$ at time $t$ is equal to $t$ minus the supremum over all $s \in [0,t]$ for which the extended ancestral lines of $i$ and $j$ at time $s$ coincide. \bigskip
%%%%%%%%%%%%%%%%%%%%%%%%%%%%%%%%%%%%%%%%%%%%%%%%%%%%%%%%%%%%%%%%%%%%%%%%%%%%%%%%%%%%%%%%%%%%%%%%%%%      Introduction, goals

The first task in the present paper is to describe how the extended ancestral lines {\it evolve forwards} in time. For this purpose, each extended ancestral line is defined on the whole real line, i.e$.$  at time $t$ it is continued on $\mathbb{R} \setminus [0,t]$ as the constant path. In addition, we use a {\it {\sc Time}-space} process that takes values in
\begin{equation}
\mathbb{R} \times (\mathcal{D}(\mathbb{R}, K \times I))^{I}
\end{equation}
to describe the HMM as a piecewise deterministic Markov jump process and to characterize it by means of a {\it well-posed martingale problem} (Theorem \ref{TheoremAnalyticalHMM}). Note that we distinguish between {\it time} which is the parameter of the {\sc Time}-space process and {\sc Time} which is the part of the state space that indicates where the action takes place in the extended ancestral lines. 

Then we focus our attention on the extended ancestral lines of a {\it tagged subset} $J \subset I$ at a fixed time horizon $T$, where we choose $-T$ as the initial {\sc Time} of the HMM. This means we consider extended ancestral lines from {\sc Time} $0$  back to {\sc Time} $-T$. In other words, we describe the situation where we can trace the extended ancestral lines of a present population back into the past. 
%%%%%%%%%%%%%%%%%%%%%%%%%%%%%%%%%%%%%%%%%%%%%%%%%%%%%%%%%%%%%%%%%%%%%%%%%%%%%%%%%%%%%%%%%%%%%%%%%%%      Introduction, BP

In order to investigate the distribution of the extended ancestral lines of $J$ at time $T$, we will introduce a pure Markov jump process on
\begin{equation}
(K \times I)^{J} \times \li(2^{K}\re)^{I}
\end{equation}
which we call {\it backward process} (short {\bf BP}), where $2^{K}$ denotes the power set of $K$. There are four important items concerning the BP we have to mention at this point.
\begin{enumerate}
\item The BP reverses migration and resampling and generates the extended ancestral lines from {\sc Time} $0$  back to {\sc Time} $-T$. In contrast to the  ancestral selection graph of Krone and Neuhauser which produces the genealogy in a three-step procedure the BP is a true Markov process.
%-----------------------------------------------------------------------
\item The first component of the BP is driven by a special kind of coalescing mechanism, but is {\it different from Kingman coalescent}, even in the case without selection and mutation. The reasons for this are  the extended ancestral lines, the goal to {\it condition} on the type information and the fact that we can incorporate mutation and selection in this BP.

The second component of the BP is solely due to selection and interacts with the first one.
%-------------------------------------------------------------------
\item There is a relation between the distribution of the extended ancestral lines of $J$ at time $T$ which is an element in
\begin{equation}
\mathcal{M}_{1}\li((\mathcal{D}([-T,0], K \times I))^{J}\re) %\;\;\;\;\mbox{ \bf (probability measures)}
\end{equation} 
and the distribution of the {\it sample paths} of the BP up to time $T$ which is an element in
\begin{equation}\label{equationM1D0T}
\mathcal{M}_{1}\li(\mathcal{D}([0,T], (K \times I)^{J} \times \li(2^{K}\re)^{I})\re)\;,
\end{equation}
where $\mathcal{M}_{1}$ denotes the set of probability measures.
%----------------------------------------------------------------------
\item There is a Feynman-Kac duality (Proposition \ref{PropositionDualityTypInfoHMM}) between the HMM and the BP with which we can express probabilities for the type information of $J$ in terms of the BP. So, in contrast to the ancestral selection graph we define {\it forward as well as backward Markovian dynamics} and relate both in terms of a {\it duality function.} 
\end{enumerate}
%%%%%%%%%%%%%%%%%%%%%%%%%%%%%%%%%%%%%%%%%%%%%%%%%%%%%%%%%%%%%%%%%%%%%%%%%%%%%%%%%%%%%%%%%%%%%%%%%%%      Introduction, summary results (more precise)

The first result concerning the relation between the extended ancestral lines and the BP is a {\it stochastic representation} (Theorem \ref{TheoremStochasticRepresentation}) that allows us to express  the expectation of certain functionals of the extended ancestral lines of $J$ at time $T$ in terms of the expectation of suitable functionals of the {\it sample paths} of the BP up to time $T$. 

The main result (Theorem \ref{TheoremStrongStochasticRepresentation}) is the aforementioned {\it strong} stochastic representation for the conditioned extended ancestral lines. It says that the conditional distribution of the {\it extended ancestral lines of $J$} at time $T$ {\it given the type information of $J$} at time $T$ is equal to the distribution of a special functional of the sample paths of a {\it transformation} of the BP up to time $T$, where the initial state of the transformed BP depends on the type information of $J$ at time $T$. This transformed BP is in general a {\it time-inhomogeneous} Markov process and arises by a special {\it change of measure} (Theorem \ref{TheoremTransformedBP}), i.e$.$ the distribution of its sample paths up to time $T$ is again an element in (\ref{equationM1D0T}). However, if the type information of the Moran model is in equilibrium, then  it turns out that the transformed BP is {\it time-homogeneous} and can be regarded as a compensated $h$-transform  which has been introduced  in \cite{FS04}.
%%%%%%%%%%%%%%%%%%%%%%%%%%%%%%%%%%%%%%%%%%%%%%%%%%%%%%%%%%%%%%%%%%%%%%%%%%%%%%%%%%%%%%%%%%%%%%%%%%%      Introduction,  applications

The main application of the strong stochastic representation is a limit theorem (Theorem \ref{TheoremLongtime}) for the extended ancestral lines. More precisely, if in the limit as $T \to \infty$ the type information at time $T$ converges to a unique stationary distribution, then the distribution of the extended ancestral lines at time $T$ (considered from {\sc Time} $0$ back to {\sc Time} $-T$) converges (in the limit $T\to \infty$) to a unique element
\begin{equation}
\mathbb{P} \in \mathcal{M}_{1}\li((\mathcal{D}((-\infty,0], K \times I))^{I}\re)
\end{equation} 
that can be represented in terms of the distribution of the sample paths (up to time $\infty$) of the {\it time-homogeneous} transformed BP. In other words, we obtain a strong stochastic representation for the extended ancestral lines in the limit $T \to \infty$ that has to two important applications: 
%---------------------------------------- CAT --------------------------- 
\begin{enumerate}
\item The stationary type distribution of the CAT can be represented in terms of the equilibrium distribution of a functional of the time-homogeneous transformed BP, where we use that the stationary type distribution of the CAT is given by
\begin{equation}
\lim_{t \to -\infty}\mathbb{P}(\mbox{the ancestral line of an arbitrary $i$ at {\sc Time} $t$} \in \, \cdot  \,) \in \mathcal{M}_{1}(K) \;.
\end{equation}
In the two type case $K = \{0,1\}$ (Proposition \ref{PropositionCAT}) this functional can be identified as a pure Markov jump process on
\begin{equation}
K \times \li\{0, \dots, N-1\re\}
\end{equation}
which converges weakly (in the limit $N \to \infty$) to the common ancestor process of Fearnhead whose equilibrium distribution has been investigated in \cite{F02}, Section 3.
%---------------------------------------- genealogical distance----------
\item  The conditioned (given the types) genealogical distance of two individuals in equilibrium can be represented in terms of a stopping time of a further functional of the time-homogeneous transformed BP.

In the two type case $K = \{0,1\}$ (Proposition \ref{PropositionGenalogicalDistance}) this functional can be identified as a pure Markov jump process on 
\begin{equation}
\{\bigtriangleup \} \cup K \times K \times \{0, \dots , N-2\}\;,
\end{equation}
and therefore the conditioned genealogical distance of two individuals in equilibrium is distributed as the first time at which this jump process reaches the absorbing state $\bigtriangleup$. 

Observe that in the limit $N \to \infty$ this jump process converges to a jump process on
\begin{equation}
\{\bigtriangleup \} \cup K \times K \times \mathbb{N}_{0}
\end{equation}
that will be used to study the tail distribution function (Proposition \ref{Proposition}) of both the genealogical distance and the conditioned genealogical distance of two individuals sampled from the tree-valued Fleming-Viot process in equilibrium. In the absence of selection we can explicitly determine the  tail distribution function of the  conditioned genealogical distance from which we recover that the genealogical distance is exponential distributed (Kingman subtree length distribution). Note that it is not clear which formal argument (maybe the notion of intertwining) can be used to get from the conditioned genealogical distance to the Kingman subtree and vice verca.

In the presence of selection we do a Taylor expansion for the tail distribution function at $0$ up to degree $3$ and we obtain that near zero (depending on the selection strength) the tail distribution function is smaller than the tail distribution function of the exponential distribution. This means we cannot yet show to the full extent that genealogical distances are stochastically smaller under selection, but we lay the foundation for proving this in the future.
\end{enumerate}

The key tool to deduce the results in this article is a Feynman-Kac duality (Theorem \ref{TheoremFeynmanKacDualityHMMandHBP}) between the HMM and the {\it historical backward process} (short {\bf HBP}) which is the {\it path process} associated to the BP (see \cite{D93} and \cite{DP91}). Observe that the HBP considered at time $t$ gives us the sample paths of the BP up to time $t$ which is the reason that we can relate the extended ancestral lines in the HMM with the sample paths of the BP by means of a duality function. This Feynman-Kac duality is used to prove the uniqueness of the martingale problem in Theorem \ref{TheoremAnalyticalHMM} and can be verified itself by a generator relation. Furthermore, it includes the stochastic representation for the extended ancestral lines (Theorem \ref{TheoremStochasticRepresentation}) and the Feynman-Kac duality between the type information of the HMM and the BP (Proposition \ref{PropositionDualityTypInfoHMM}). We emphasize that the {\it concept of life-sites} (which only plays a minor role in the applications) is the key ingredient to relate the extended ancestral lines in the HMM with the sample paths of the BP by means of a duality function. \bigskip

This paper is organized as follows.  In Section \ref{secFormulationmodels} we define the HMM and the BP. Then we state our results in Section \ref{secFormulationresults}. In Section \ref{secHBP} we introduce the HBP and provide our key tool, the Feynman-Kac duality between the HMM and the HBP. Finally, we prove our results in Section \ref{secProofs}.
%%%%%%%%%%%%%%%%%%%%%%%%%%%%%%%%%%%%%%%%%%%%%%%%%%%%%%%%%%%%%%%%%%%%%%%%%%%%%%%%%%%%%%%%%%%%%%%%%%%      Formulation of the main results

\setcounter{equation}{0}
\section{Formulation of the models}\label{secFormulationmodels}
\noindent This section is concerned with the introduction of the HMM and the BP which includes an analytical characterization in each case.

Let $B$ be a non-negative real number that represents the {\it rate of mutation}. In order to describe type-dependent mutation we consider  a general stochastic matrix
\begin{equation}
b(u,v), \;\; u,v \in K\;.
\end{equation}
Without loss of generality we can assume that resampling occurs at rate $1$.

Furthermore, recall that $K = \{0,1, \dots,d-1 \}$, let $S \in [0,N]$ be the {\it selection coefficient} and 
\begin{equation}\label{equationDefChi}
\chi: K \to [0,1] \;\;\mbox{ with }\;\; 0 = \chi(0) < \dots <\chi(d-1) = 1
\end{equation}
be the function that assigns to each type its level of fitness,  i.e$.$ there are $d$ different levels of fitness, type $d-1$ has the highest level of fitness and type $0$ has the lowest level of fitness. 
%%%%%%%%%%%%%%%%%%%%%%%%%%%%%%%%%%%%%%%%%%%%%%%%%%%%%%%%%%%%%%%%%%%%%%%%%%%%%%%%%%%%%%%%%%%%%%%%%%%      Formulation of the models, The HMM 

\subsection{The historical Moran Model (HMM)}
\noindent In order to describe the evolution of the extended ancestral lines forward in time we define the HMM by means of a {\it piecewise deterministic Markov jump process} on
\begin{equation}
\mathcal{E} \subset \mathbb{R} \times (\mathcal{D}(\mathbb{R}, K \times I))^{I}\;\;\; \mbox{\bf (state space)}  
\end{equation}
where
\begin{equation}
\mathcal{E} := \li\{\eta = (\eta^{\mbox{\tiny \sc Time}}, \eta^{\mbox{\tiny $\mathcal{D}$}}) = \li(\eta^{\mbox{\tiny \sc Time}}, \li(\eta^{\mbox{\tiny $\mathcal{D}$}}_{i}\re)_{i \in I}\re): \eta^{\mbox{\tiny $\mathcal{D}$}}_{i,s} = \li(\eta^{\mbox{\tiny $\mathcal{D}$}}_{i,\eta^{\mbox{\tiny \sc Time}},K},i\re) \mbox{ for all } i \in I,  s \ge \eta^{\mbox{\tiny \sc Time}}\re\}.
\end{equation}
An element $\eta \in \mathcal{E}$ describes the collection of extended ancestral lines at {\sc Time} $\eta^{\mbox{\tiny \sc Time}}$. For $s \le \eta^{\mbox{\tiny \sc Time}}$,
\begin{equation}
\eta^{\mbox{\tiny $\mathcal{D}$}}_{i,s} = \li(\eta^{\mbox{\tiny $\mathcal{D}$}}_{i,s,K},\eta^{\mbox{\tiny $\mathcal{D}$}}_{i,s,I}\re)
\end{equation}
is the pair of type and life-site of the ancestor of $i$ alive at {\sc Time} $s$, where
\begin{equation}
\eta^{\mbox{\tiny $\mathcal{D}$}}_{i,\eta^{\mbox{\tiny \sc Time}}} = \li(\mbox{type of $i$ at {\sc Time} $\eta^{\mbox{\tiny \sc Time}}$}, i\re)
\end{equation}
since evidently $i$ is the life-site of the ancestor of $i$ at {\sc Time} $\eta^{\mbox{\tiny \sc Time}}$. We need  the function
\begin{equation}\label{equationProjectionOnTypes}
(\cdot)^{*}: \mathcal{E} \to K^{I}, \; \eta \mapsto \li(\eta^{*}_{i}\re)_{i \in I} =  \li(\eta^{\mbox{\tiny $\mathcal{D}$}}_{i,\eta^{\mbox{\tiny \sc Time}},K}\re)_{i \in I}\;\;\; \mbox{ {\bf (projection on types)}} 
\end{equation}
to obtain the types of the population at {\sc Time} $\eta^{\mbox{\tiny \sc Time}}$. Observe that the function $(\cdot)^{*}$ is continuous since, by definition of $\mathcal{E}$, each extended ancestral is constant on $[\eta^{\mbox{\tiny \sc Time}},\infty)$. 
%%%%%%%%%%%%%%%%%%%%%%%%%%%%%%%%%%%%%%%%%%%%%%%%%%%%%%%%%%%%%%%%%%%%%%%%%%%%%%%%%%%%%%%%%%%%%%%%%%%      Formulation of the models, The HMM , Description of states

First we provide a description of changes of states in $\mathcal{E}$. \begin{enumerate}
\item {\bf (Mutation)} For $\eta \in \mathcal{E}$, $i \in I$ and $u \in K$ we define the element
\begin{equation}
 \eta^{i;u} \in \mathcal{E} \;\mbox{ by }\; (\eta^{i;u})^{\mbox{\tiny \sc Time}} := \eta^{\mbox{\tiny \sc Time}}  \mbox{ and } (\eta^{i;u})^{\mbox{\tiny $\mathcal{D}$}}_{l,s}  := \li\{ \begin{array}{ccl} \eta^{\mbox{\tiny $\mathcal{D}$}}_{l,s} &,& s \in \mathbb{R}, l \not= i\\ \eta^{\mbox{\tiny $\mathcal{D}$}}_{i,s} &,& s < \eta^{\mbox{\tiny \sc Time}}, l = i\\ (u,i) &,&  s \ge \eta^{\mbox{\tiny {\sc Time}}}, l = i\end{array} \re. \;.
\end{equation}
\item {\bf (Resampling)} For $\eta \in \mathcal{E}$ and $i,j \in I$ we define the element
\begin{equation}
 \eta^{i\to j} \in \mathcal{E} \;\mbox{ by }\; (\eta^{i \to j})^{\mbox{\tiny \sc Time}} := \eta^{\mbox{\tiny \sc Time}}  \mbox{ and } (\eta^{i\to j})^{\mbox{\tiny $\mathcal{D}$}}_{l,s}  := \li\{ \begin{array}{ccl} \eta^{\mbox{\tiny $\mathcal{D}$}}_{l,s} &,& s \in \mathbb{R}, l \not= j\\ \eta^{\mbox{\tiny $\mathcal{D}$}}_{i,s} &,& s < \eta^{\mbox{\tiny \sc Time}}, l = j\\ (\eta^{*}_{i},j) &,&  s \ge \eta^{\mbox{\tiny {\sc Time}}}, l = j\end{array} \re.\;.
\end{equation}
\end{enumerate}
%%%%%%%%%%%%%%%%%%%%%%%%%%%%%%%%%%%%%%%%%%%%%%%%%%%%%%%%%%%%%%%%%%%%%%%%%%%%%%%%%%%%%%%%%%%%%%%%%%%      Formulation of the models, The HMM, Definition

Now we can define the HMM according to the description in item \ref{itemMutation} and item \ref{itemResampling} on page \pageref{itemMutation}. Observe that  we initial assign to each life-site $i \in I$ a constant path through $(u_{i},i)$, where the types $(u_{i})_{i \in I}$ a chosen according to a distribution on $K^{I}$. 
\begin{Definition}[HMM]\label{DefinitionHMM}\pointpar
\noindent The HMM is a piecewise deterministic Markov jump process on $\mathcal{E}$ whose {\bf initial state} can be characterized by a probability measure $\mu$ on $\mathcal{E}$ of the form
\begin{equation}\label{equationMu}
\mu(\{\eta\}) = \sum_{(u_{i})_{i \in I} \in K^{I}}\delta_{c}(\{\eta^{\mbox{\tiny \sc Time}}\}) \prod_{i \in I}\delta_{[s \mapsto (u_{i},i)]}(\{\eta^{\mbox{\tiny $\mathcal{D}$}}\}) \mu^{*}(\{(u_{i})_{i \in I}\})\;,
\end{equation}
where 
\begin{itemize}
\item $c$ is a real number describing the initial {\sc Time},
\item $[s \mapsto (u_{i},i)]$ denotes the constant path through $(u_{i},i)$ and
\item $\mu^{*}$ is a probability measure on $K^{I}$ describing the initial type distribution.
\end{itemize}
The evolution of the HMM is given as follows:
%%%%%%%%%%%%%%%%%%%%%%%%%%%%%%%%%%%%%%%%%%%%%%%%%%%%%%%%%%%%%%%%%%%%%%%%%%%%%%%%%%%%%%%%%%%%%%%%%%%      Formulation of the models, The HMM, Definition

The {\sc Time} coordinate grows with unit speed. If $\eta \in \mathcal{E}$ is the current state of the HMM, then the following transitions, depending on the {\sc Time}, occur independently for all life-sites, independently for all types and independently of each other:
\begin{enumerate}
\item {\bf (Mutation)} For each $i \in I$ and each $u \in K$ the transition
\begin{equation}
\eta\to \eta^{i ; u}\;\;\mbox{ occurs at rate }\;\;B b(\eta^{*}_{i},u)\;.
\end{equation}
\item {\bf (Resampling)} For each $i,j  \in I$ the transition
\begin{equation}
\eta \to \eta^{i \to j} \;\;\mbox{ occurs at rate }\;\;  \frac{1}{2} + \frac{S}{2N}\li[\chi(\eta^{*}_{i}) - \chi(\eta^{*}_{j})\re]\;,
\end{equation}
where the rate depends on the fitness of $i$ and $j$ {\bf (Selection)}.
\end{enumerate}
\end{Definition}
\begin{Remark}\point\label{RemarkDefinitionHMM} Definition \ref{DefinitionHMM} implies a natural property for the extended ancestral lines, namely, whenever the $I$-coordinate of two extend ancestral lines coincide at {\sc Time} $t$ then also the $K$-coordinate does and the two extend ancestral lines coincide at all {\sc Times} $s < t$.
\end{Remark} 
%%%%%%%%%%%%%%%%%%%%%%%%%%%%%%%%%%%%%%%%%%%%%%%%%%%%%%%%%%%%%%%%%%%%%%%%%%%%%%%%%%%%%%%%%%%%%%%%%%%      Formulation of the models, The BP

\subsection{The backward process (BP)}
\noindent Here we introduce the BP, a pure Markov jump process, which shall be used to determine for a general initial type distribution $\mu^{*}$ the expectation of certain functionals of the extended ancestral lines of a tagged $J \subset I$ at a fixed time $T$ backwards in time. The state space for the BP consists of two components, namely
\begin{equation}\label{equationOverlineMathcalE}
\overline{\mathcal{E}} := \li\{\overline{\eta} = (\overline{\eta}^{\mbox{\tiny $J$}}, \overline{\eta}^{\mbox{\tiny $I$}}) \in (K \times I)^{J} \times \li( 2^{K} \setminus \{\emptyset\}\re)^{I}: \overline{\eta}^{\mbox{\tiny $J$}}_{i} = \overline{\eta}^{\mbox{\tiny $J$}}_{j} \mbox{ whenever }  \overline{\eta}^{\mbox{\tiny $J$}}_{i,I} = \overline{\eta}^{\mbox{\tiny $J$}}_{j,I}\re\}\,,
\end{equation}
where we allow (this is a slight change compared to  Definition 1.9 in \cite{Seidel15}) to have only non-empty subsets of $K$ in the second component of the BP  which is a technical assumption useful later on (see Subsubsection \ref{subsubTransformeddualprocess}).  Before we give a rigorous definition we illustrate and explain the BP in simple terms. \bigskip

The first component of the BP (more precisely, the sample paths of first component up to time $T$) shall describe {\it the reversed extended ancestral lines of $J$ at time $T$} (i.e$.$ the initial types in the first component of the BP are the types of $J$ in the HMM at time $T$), where the restriction in (\ref{equationOverlineMathcalE}) reflects the natural property of extended ancestral lines from Remark \ref{RemarkDefinitionHMM}. The second component of the BP is solely due to {\it selection} and stays constant for all times if $S = 0$. An important property of the BP is that each first component
\begin{equation}
\overline{\eta}^{\mbox{\tiny $J$}} \in (K \times I)^{J}
\end{equation}
defines a {\it partition} of $J$ in which each partition element is {\it marked} by a type in $K$ and {\it located} at a life-site in $I$.  
%%%%%%%%%%%%%%%%%%%%%%%%%%%%%%%%%%%%%%%%%%%%%%%%%%%%%%%%%%%%%%%%%%%%%%%%%%%%%%%%%%%%%%%%%%%%%%%%%%%      Formulation of the models, The BP,  absence of selction

In the absence of selection the evolution of these $K$-marked partition elements on $I$ can be described by a system of instantaneously coalescing random walks on $I$ with a {\it special coalescing rule}. Namely, each $K$-marked partition element undergoes a random walk on $I$, where the mark itself evolves as a Markov chain on $K$. At a jump time the partition element chooses its new location uniformly from $I$ (in the forward view this represents the choice of an ancestor in a resampling event as descirbed in item \ref{itemResampling} on page \pageref{itemResampling}) and immediately coalesce if the chosen life-site is occupied by an other partition element. The special coalescing rule says that a partition element can only jump to an occupied life-site if the partition element (which is located there) has the same mark in $K$.  This means that our BP is {\it different from Kingman coalescent}, even in the absence of selection.  

In the presence of selection each coalescing event additionally leads to a change in the second component of the BP which has two important consequences. On the one hand, the evolution of the first component depends on the state of the second component, e.g$.$ the jump of a partition element to a life-site $i \in I$ depends on the state of $i$ in the second component, i.e$.$ on
\begin{equation}
\overline{\eta}^{\mbox{\tiny $I$}}_{i} \in  2^{K}\setminus \{\emptyset\}\;.
\end{equation} 
On the other hand, the evolution of the second component depends on the first component of the BP, in particular, we call  
\begin{equation}
\mbox{ the life-site $i \in I$ in the second component {\it active} }
\end{equation}
if and only if $i$ is not the location of a partition element. This means that the jump of a partition element from {\it $j$ to $i$ activates the life-site $j$ and deactivates the life-site $i$} in the second component.

To put it in a nutshell, the jump of a partition element located on $j$ to the life-site $i$ shall describe the reversal of a resampling event between the life-sites $i$ and $j$ in the HMM in which $i$ is chosen as ancestor (i.e$.$ the individual located on $j$ dies), where the special coalescing rule is due to the fact that both life-sites have the same type after a resampling event. The {\it active} (in the backward view) life-sites in the second component of the BP shall code the types of {\it dead individuals} in the HMM (forward view). Namely, in the case with selection the rate of a resampling event does depend on the type of the dying individual, the information which gets lost in a resampling event, and hence the BP has to reproduce this information in some form, for a detailed explanation see item \ref{itemAnnotation3} subsequent to Definition \ref{DefinitionBP}.   \bigskip
%%%%%%%%%%%%%%%%%%%%%%%%%%%%%%%%%%%%%%%%%%%%%%%%%%%%%%%%%%%%%%%%%%%%%%%%%%%%%%%%%%%%%%%%%%%%%%%%%%%      Formulation of the models, The BP, marked partition

In order to introduce the BP we first have to rigorously define the concept of a $K$-marked $J$-partition on $I$ included in the first component of the BP and the concept of active life-sites.
\begin{itemize}
\item {\bf $K$-marked $J$-partition on $I$:} To each state $\overline{\eta} \in \overline{\mathcal{E}}$ we can assign a partition of $J$ denoted by 
\begin{equation}
\Gamma(\overline{\eta})\;,
\end{equation}
formally
\begin{equation}
\Gamma : \overline{\mathcal{E}} \to \li\{\mathcal{C}\in 2^{J} \setminus \{\emptyset\} : \gamma \cap \gamma' = \emptyset \; \forall \gamma \not= \gamma' \in \mathcal{C}\re\}\;.
\end{equation}
A partition element $\gamma \in \Gamma\li(\overline{\eta}\re)$ is defined by the rule
\begin{equation}
 j, i \in \gamma \;\iff\; \overline{\eta}^{\mbox{\tiny $J$}}_{j} = \overline{\eta}^{\mbox{\tiny $J$}}_{i} \;.
\end{equation}
Due to  the restriction in (\ref{equationOverlineMathcalE}) for each element $\overline{\eta} \in \overline{\mathcal{E}}$ we can assign to each partition element $\gamma \in \Gamma(\overline{\eta})$ a unique mark in $K$ and a unique location in $I$ by
\begin{equation}
\overline{\eta}^{\mbox{\tiny $J$}}_{\gamma} = \li(\overline{\eta}^{\mbox{\tiny $J$}}_{\gamma,K},\overline{\eta}^{\mbox{\tiny $J$}}_{\gamma,I}\re) := \li(\overline{\eta}^{\mbox{\tiny $J$}}_{i,K},\overline{\eta}^{\mbox{\tiny $J$}}_{i,I}\re)  \in K \times I\;,
\end{equation}
where $i \in \gamma$ is arbitrary.  This means that $\gamma$ represents the descendants of an ancestor that has type $\overline{\eta}^{\mbox{\tiny $J$}}_{\gamma,K}$ and occupies the life-site $\overline{\eta}^{\mbox{\tiny $J$}}_{\gamma,I}$ in the HMM.

In addition, we write
\begin{equation}
\overline{\eta}^{\mbox{\tiny $J$}}_{\gamma} = (u,i) \; \iff \; \li(\overline{\eta}^{\mbox{\tiny $J$}}_{i,K},\overline{\eta}^{\mbox{\tiny $J$}}_{i,I}\re) =  (u,i) \; \mbox{ for all } i \in \gamma\;.
\end{equation}  
\item {\bf The active life-sites} code the type information of dead individuals in the HMM. Formally, they are given by the set
\begin{equation}
\tilde{\Gamma}(\overline{\eta}) := I \setminus \li\{\overline{\eta}^{\mbox{\tiny $J$}}_{\gamma,I} : \gamma \in  \Gamma(\overline{\eta}) \re\} \; , 
\end{equation}
that is, the life-site $i \in I$ in the second component of the BP is an active life-site if and only if $i$ is not the location of a partition element.
\end{itemize}
%%%%%%%%%%%%%%%%%%%%%%%%%%%%%%%%%%%%%%%%%%%%%%%%%%%%%%%%%%%%%%%%%%%%%%%%%%%%%%%%%%%%%%%%%%%%%%%%%%%      Formulation of the models, The BP, Definition

Now we come to the definition of the BP. We want to reverse the dynamics of the HMM in order to follow extended ancestral lines back into the past. In the following Definition \ref{DefinitionBP} we first specify the transitions that reverse mutation and then the transitions that reverse resampling, where those transitions which occur only in the presence of selection are indicated by ($\star$). An overview of all the different kinds of transitions is given in Figure \ref{FigureTransitionsBP} at the end of this subsection. For additional information about these transitions see also the explanations subsequent to Definition \ref{DefinitionBP}. Recall also (\ref{equationDefChi}) and note that $\mathbbm{1}\{\mbox{statement}\}$ is $1$ if the statement is true and $0$ if the statement is false. 

\begin{Definition}[The BP]\label{DefinitionBP}\pointpar
\noindent The BP is a pure Markov jump process on $\overline{\mathcal{E}}$. 

Its {\bf initial state} is given by the element
\begin{equation}\label{equationoverlinexi}
\overline{\xi^{*}} = (\overline{\xi^{*}}^{\mbox{\tiny $J$}},\overline{\xi^{*}}^{\mbox{\tiny $I$}}) \in \overline{\mathcal{E}} \;\;\mbox{ defined by }\;\; \overline{\xi^{*}}^{\mbox{\tiny $J$}}_{j} = (\xi^{*}_{j},j) \; \forall j \in J \;\mbox{ and }\; \overline{\xi^{*}}^{\mbox{\tiny $I$}}_{i} = K  \; \forall i\in I\;, 
\end{equation} 
where the element $\xi^{*} \in K^{J}$ will represent the types on which we condition $J$ in the HMM at time $T$ later on. This means
\begin{equation}
\Gamma(\overline{\xi^{*}}) = \{\{j\} : j \in J\} \;\mbox{ and }\; \tilde{\Gamma}(\overline{\xi^{*}}) = I \setminus J\;,
\end{equation} 
where the partition element $\{j\}$ is located at $j$ and marked by $\xi^{*}_{j}$.

The transition rates are given by a transition matrix
\begin{equation}
\overline{\mathcal{K}}:\overline{\mathcal{E}} \times \overline{\mathcal{E}} \to [0,\infty)
\end{equation}
that is positive on the following pairs of states, where we first state the rate and then specify the new element to which the transition leads: 
\begin{enumerate}
%------------------------------------------ reversal mutation -----
\item {\bf (Reversal of mutation)}
\begin{enumerate}
\item \label{itemdualmutationJ}{\bf [Transitions for partition elements]}\\    
For each $\gamma \in \Gamma(\overline{\eta})$ and each $u \in K$,
\begin{equation}
\overline{\mathcal{K}}(\overline{\eta}, \overline{\eta}^{\gamma;u})  = B b(u, \overline{\eta}^{\mbox{\tiny $J$}}_{\gamma,K})\;,
\end{equation}
where the element $\overline{\eta}^{\gamma;u} \in \overline{\mathcal{E}}$ is defined by 
\begin{equation} 
\li(\overline{\eta}^{\gamma;u}\re)^{\mbox{\tiny $J$}}_{\gamma'} := \li\{ \begin{array}{ccl} (u,\overline{\eta}^{\mbox{\tiny $J$}}_{\gamma,I}) &,& \gamma' = \gamma  \\ \overline{\eta}^{\mbox{\tiny $J$}}_{\gamma'} &,& \gamma' \not= \gamma \end{array} \re. \;\mbox{ and }\; (\overline{\eta}^{\gamma;u})^{\mbox{\tiny $I$}} := \overline{\eta}^{\mbox{\tiny $I$}}\;.
\end{equation}
%-----------------------------------------------------------------------
\item \label{itemdualmutationI}{\bf [Transitions for active life-sites]}   
\begin{enumerate}
\item ($\star$) \label{itemdualmutationIcup}  For each $i \in \tilde{\Gamma}(\overline{\eta})$, each $v \in \overline{\eta}^{\mbox{\tiny $I$}}_{i}$ and each $u \in K \setminus \overline{\eta}^{\mbox{\tiny $I$}}_{i}$, 
\begin{equation}
\overline{\mathcal{K}}(\overline{\eta}, \overline{\eta}^{i;\cup\{u\}})  = B b(u,v)\;,
\end{equation}
where the element $\overline{\eta}^{i;\cup\{u\}} \in \overline{\mathcal{E}}$ is defined by 
\begin{equation} 
\li(\overline{\eta}^{i;\cup\{u\}}\re)^{\mbox{\tiny $J$}} :=  \overline{\eta}^{\mbox{\tiny $J$}} \;\mbox{ and }\; \li(\overline{\eta}^{i;\cup\{u\}}\re)^{\mbox{\tiny $I$}}_{l} := \li\{ \begin{array}{ccl} \overline{\eta}^{\mbox{\tiny $I$}}_{i} \cup \{u\} &,& l = i \\ \overline{\eta}^{\mbox{\tiny $I$}}_{l} &,& l \not= i \end{array} \re.\;.
\end{equation}
%-----------------------------------------------------------------------
\item  ($\star$) \label{itemdualmutationIsetminus} For each $i \in \tilde{\Gamma}(\overline{\eta})$, each $v \in K\setminus\overline{\eta}^{\mbox{\tiny $I$}}_{i}$ and each $u \in \overline{\eta}^{\mbox{\tiny $I$}}_{i}$, 
\begin{equation}
\overline{\mathcal{K}}(\overline{\eta}, \overline{\eta}^{i;\setminus\{u\}})  = \mathbbm{1}\{|\overline{\eta}^{\mbox{\tiny $I$}}_{i}| > 1\} B b(u,v)\;,
\end{equation}
where the element $\overline{\eta}^{i;\setminus\{u\}} \in \overline{\mathcal{E}}$ is defined by 
\begin{equation} 
\li(\overline{\eta}^{i;\setminus \{u\}}\re)^{\mbox{\tiny $J$}} :=  \overline{\eta}^{\mbox{\tiny $J$}} \;\mbox{ and }\; \li(\overline{\eta}^{i;\setminus \{u\}}\re)^{\mbox{\tiny $I$}}_{l} := \li\{ \begin{array}{ccl} \overline{\eta}^{\mbox{\tiny $I$}}_{i} \setminus \{u\} &,& l = i \\ \overline{\eta}^{\mbox{\tiny $I$}}_{l} &,& l \not= i \end{array} \re.\;.
\end{equation}
\end{enumerate}
\end{enumerate}
%-------------------------------------reversal resampling J \times J
\item {\bf (Reversal of resampling)}
\begin{enumerate}
\item \label{itemdualresamplingJJ} {\bf [Interactions between partition elements]}
\begin{enumerate}
%-----------------------------------------------------------------------
\item \label{itemdualresamplingJJK}  For each $\gamma, \gamma' \in  \Gamma(\overline{\eta})$ with $\gamma \not=\gamma'$,
\begin{equation}
\overline{\mathcal{K}}(\overline{\eta}, \overline{\eta}^{\gamma \to \gamma'})  = \mathbbm{1}\{\overline{\eta}^{\mbox{\tiny $J$}}_{\gamma,K} =\overline{\eta}^{\mbox{\tiny $J$}}_{\gamma',K}\} \li(\frac{1}{2} + \frac{S}{2N}\li[\chi(\overline{\eta}^{\mbox{\tiny $J$}}_{\gamma,K}) - 1\re]\re)\;,
\end{equation}
where the element $\overline{\eta}^{\gamma \to \gamma'} \in \overline{\mathcal{E}}$ is defined by 
\begin{equation} 
\li(\overline{\eta}^{\gamma \to \gamma'}\re)^{\mbox{\tiny $J$}}_{\gamma''} := \li\{ \begin{array}{ccl} \overline{\eta}^{\mbox{\tiny $J$}}_{\gamma'} &,& \gamma'' = \gamma \\ \overline{\eta}^{\mbox{\tiny $J$}}_{\gamma''} &,& \gamma'' \not= \gamma \end{array} \re. 
\end{equation}
and
\begin{equation}
\li(\overline{\eta}^{\gamma \to \gamma'}\re)^{\mbox{\tiny $I$}}_{l} := \li\{ \begin{array}{ccl} K &,& l = \overline{\eta}^{\mbox{\tiny $J$}}_{\gamma, I} \\ \overline{\eta}^{\mbox{\tiny $I$}}_{l} &,& l \not= \overline{\eta}^{\mbox{\tiny $J$}}_{\gamma,I} \end{array} \re.\;.
\end{equation}
This means
\begin{equation}
\Gamma(\overline{\eta}^{\gamma \to \gamma'}) = \Gamma(\overline{\eta}) \cup \{\gamma \cup \gamma'\} \setminus \{\gamma,\gamma'\} \;\mbox{ and }\; \tilde{\Gamma}(\overline{\eta}^{\gamma \to \gamma'}) = \tilde{\Gamma}(\overline{\eta}) \cup\{\overline{\eta}^{\mbox{\tiny $J$}}_{\gamma,I}\}.
\end{equation}
%-----------------------------------------------------------------------
\item  ($\star$)  \label{itemdualresamplingJJw}  For each $\gamma,\gamma' \in  \Gamma(\overline{\eta})$  with $\gamma \not=\gamma' \in  \Gamma(\overline{\eta})$ and each $w \in \{0,\dots, d-2\}$,
\begin{equation}
\overline{\mathcal{K}}(\overline{\eta}, \overline{\eta}^{\gamma \stackrel{w}{\to} \gamma'})  =  \mathbbm{1}\{\overline{\eta}^{\mbox{\tiny $J$}}_{\gamma,K} =\overline{\eta}^{\mbox{\tiny $J$}}_{\gamma',K}\} \frac{ S}{2N}\li[\chi(w+1) - \chi(w)\re]\;,
\end{equation}
where the element $\overline{\eta}^{\gamma \stackrel{w}{\to} \gamma'} \in \overline{\mathcal{E}}$ is defined by 
\begin{equation} 
\li(\overline{\eta}^{\gamma \stackrel{w}{\to} \gamma'}\re)^{\mbox{\tiny $J$}}:= \li(\overline{\eta}^{\gamma \to \gamma'}\re)^{\mbox{\tiny $J$}} 
\end{equation}
and
\begin{equation}
\li(\overline{\eta}^{\gamma \stackrel{w}{\to} \gamma'}\re)^{\mbox{\tiny $I$}}_{l} := \li\{ \begin{array}{ccl} \{0,\dots, w\} &,& l = \overline{\eta}^{\mbox{\tiny $J$}}_{\gamma, I} \\ \overline{\eta}^{\mbox{\tiny $I$}}_{l} &,& l \not= \overline{\eta}^{\mbox{\tiny $J$}}_{\gamma,I} \end{array} \re.\;.
\end{equation}
\end{enumerate}
%------------------------------------reversal resampling J vs I
\item \label{itemdualresamplingJI} {\bf [Interactions: partition elements $\to$ active life-sites]}
\begin{enumerate}
%-----------------------------------------------------------------------
\item \label{itemdualresamplingJIK}  For each $\gamma \in  \Gamma(\overline{\eta})$ and each  $i \in \tilde{\Gamma}(\overline{\eta})$,
\begin{equation}
\overline{\mathcal{K}}(\overline{\eta}, \overline{\eta}^{\gamma \to i})  = \mathbbm{1}\{\overline{\eta}^{\mbox{\tiny $J$}}_{\gamma,K}  \in \overline{\eta}^{\mbox{\tiny $I$}}_{i}\} \li(\frac{1}{2} + \frac{S}{2N}\li[\chi(\overline{\eta}^{\mbox{\tiny $J$}}_{\gamma,K}) - 1\re]\re)\;,
\end{equation}
where the element $\overline{\eta}^{\gamma \to i} \in \overline{\mathcal{E}}$ is defined (recall \ref{itemdualresamplingJJK}) by 
\begin{equation} 
\li(\overline{\eta}^{\gamma \to i}\re)^{\mbox{\tiny $J$}}_{\gamma'} := \li\{ \begin{array}{ccl} (\overline{\eta}^{\mbox{\tiny $J$}}_{\gamma,K},i) &,& \gamma' = \gamma \\ \overline{\eta}^{\mbox{\tiny $J$}}_{\gamma'} &,& \gamma' \not= \gamma \end{array} \re. \;\mbox{ and }\; (\overline{\eta}^{\gamma \to i})^{\mbox{\tiny $I$}} := (\overline{\eta}^{\gamma \to \gamma'})^{\mbox{\tiny $I$}}\; .
\end{equation}
This means
\begin{equation}
\Gamma(\overline{\eta}^{\gamma \to i}) = \Gamma(\overline{\eta})  \;\mbox{ and }\; \tilde{\Gamma}(\overline{\eta}^{\gamma \to i}) = \tilde{\Gamma}(\overline{\eta}) \cup \{\overline{\eta}^{\mbox{\tiny $J$}}_{\gamma,I}\} \setminus \{i\}\;.
\end{equation}
%-----------------------------------------------------------------------
\item  ($\star$)  \label{itemdualresamplingJIw}  For each $\gamma \in  \Gamma(\overline{\eta})$, each  $i \in \tilde{\Gamma}(\overline{\eta})$ and each $w \in \{0,\dots, d-2\}$,
\begin{equation}
\overline{\mathcal{K}}(\overline{\eta}, \overline{\eta}^{\gamma \stackrel{w}{\to} i})  =  \mathbbm{1}\{\overline{\eta}^{\mbox{\tiny $J$}}_{\gamma,K} \in \overline{\eta}^{\mbox{\tiny $I$}}_{i}\} \frac{S}{2N}\li[\chi(w+1) - \chi(w)\re]\;,
\end{equation}
where the element $\overline{\eta}^{\gamma \stackrel{w}{\to} i} \in \overline{\mathcal{E}}$ is defined (recall \ref{itemdualresamplingJJw}) by 
\begin{equation} 
\li(\overline{\eta}^{\gamma \stackrel{w}{\to} i}\re)^{\mbox{\tiny $J$}}_{\gamma'} := \li(\overline{\eta}^{\gamma \to i}\re)^{\mbox{\tiny $J$}}_{\gamma'}  \;\mbox{ and }\; \li(\overline{\eta}^{\gamma \stackrel{w}{\to} i}\re)^{\mbox{\tiny $I$}} := \li(\overline{\eta}^{\gamma \stackrel{w}{\to} \gamma'}\re)^{\mbox{\tiny $I$}}\;.
\end{equation}
\end{enumerate}
%------------------------------------reversal resampling I vs J
\item \label{itemdualresamplingIJ} {\bf [Interactions: active life-sites $\to$ partition elements]}
\begin{enumerate}
%-----------------------------------------------------------------------
\item ($\star$)\label{itemdualresamplingIJK}  For each $i \in \tilde{\Gamma}(\overline{\eta})$ and each  $\gamma \in  \Gamma(\overline{\eta})$,
\begin{equation}
\overline{\mathcal{K}}(\overline{\eta}, \overline{\eta}^{i;K})  = \mathbbm{1}\{\overline{\eta}^{\mbox{\tiny $J$}}_{\gamma,K}  \in \overline{\eta}^{\mbox{\tiny $I$}}_{i}\} \li(\frac{1}{2} + \frac{S}{2N}\li[\chi(\overline{\eta}^{\mbox{\tiny $J$}}_{\gamma,K}) -1\re]\re)\;,
\end{equation}
where the element $\overline{\eta}^{i;K} \in \overline{\mathcal{E}}$ is defined by 
\begin{equation} 
\li(\overline{\eta}^{i;K}\re)^{\mbox{\tiny $J$}} :=  \overline{\eta}^{\mbox{\tiny $J$}} \;\mbox{ and }\; \li(\overline{\eta}^{i;K}\re)^{\mbox{\tiny $I$}}_{l} := \li\{ \begin{array}{ccl} K &,& l = i\\ \overline{\eta}^{\mbox{\tiny $I$}}_{l} &,& l \not= i \end{array} \re.\;.
\end{equation}
%-----------------------------------------------------------------------
\item ($\star$)\label{itemdualresamplingIJw}  For each $i \in \tilde{\Gamma}(\overline{\eta})$, each $\gamma \in  \Gamma(\overline{\eta})$ and each $w \in \{0,\dots, d-2\}$,
\begin{equation}
\overline{\mathcal{K}}(\overline{\eta}, \overline{\eta}^{i;w})  = \mathbbm{1}\{\overline{\eta}^{\mbox{\tiny $J$}}_{\gamma,K} \in \overline{\eta}^{\mbox{\tiny $I$}}_{i}\} \frac{S}{2N}\li[\chi(w+1) - \chi(w)\re]\;,
\end{equation}
where the element $\overline{\eta}^{i;w} \in \overline{\mathcal{E}}$ is defined by 
\begin{equation} 
\li(\overline{\eta}^{i;w}\re)^{\mbox{\tiny $J$}} :=  \overline{\eta}^{\mbox{\tiny $J$}} \;\mbox{ and }\; \li(\overline{\eta}^{i;w}\re)^{\mbox{\tiny $I$}}_{l} := \li\{ \begin{array}{ccl} \{0,\dots, w\} &,& l = i \\ \overline{\eta}^{\mbox{\tiny $I$}}_{l} &,& l \not= i \end{array} \re.\;.
\end{equation}
\end{enumerate}
%------------------------------------reversal resampling I 
\item\label{itemdualresamplingII} {\bf [Interactions between active life-sites]}
\begin{enumerate}
%-----------------------------------------------------------------------
\item  ($\star$) \label{itemdualresamplingIIK}  For each $i,j \in \tilde{\Gamma}(\overline{\eta})$ with $i \not=j$, 
\begin{equation}
 \overline{\mathcal{K}}(\overline{\eta}, \overline{\eta}^{i \cap j})  = \mathbbm{1}\{\overline{\eta}^{\mbox{\tiny $I$}}_{i} \cap \overline{\eta}^{\mbox{\tiny $I$}}_{j}   \not\in \{\emptyset, K\} \}  \li(\frac{1}{2} + \frac{S}{2N}\li[\chi(\min \overline{\eta}^{\mbox{\tiny $I$}}_{i} \cap \overline{\eta}^{\mbox{\tiny $I$}}_{j}) - 1\re]\re),
\end{equation}
where the element $\overline{\eta}^{i\cap j} \in \overline{\mathcal{E}}$ is defined by 
\begin{equation} 
\li(\overline{\eta}^{i\cap j}\re)^{\mbox{\tiny $J$}} :=  \overline{\eta}^{\mbox{\tiny $J$}} \;\mbox{ and }\; \li(\overline{\eta}^{i \cap j}\re)^{\mbox{\tiny $I$}}_{l} := \li\{ \begin{array}{ccl} \overline{\eta}^{\mbox{\tiny $I$}}_{i} \cap \overline{\eta}^{\mbox{\tiny $I$}}_{j} &,& l = i \\  K &,& l = j \\ \overline{\eta}^{\mbox{\tiny $I$}}_{l} &,& l \not=  i,j \end{array} \re.\;.
\end{equation}
%-----------------------------------------------------------------------
\item  ($\star$)\label{itemdualresamplingIIw} For each $i,j \in \tilde{\Gamma}(\overline{\eta})$ with $i \not=j$ and each $w \in \{0,\dots, d-2\}$,
\begin{equation}
\overline{\mathcal{K}}(\overline{\eta}, \overline{\eta}^{i \stackrel{w}{\cap} j})  = \mathbbm{1}\{\overline{\eta}^{\mbox{\tiny $I$}}_{i} \cap \overline{\eta}^{\mbox{\tiny $I$}}_{j}   \not\in \{\emptyset, K\}\} \frac{S}{2N}\li[\chi(w+1) - \chi(w)\re]\;,
\end{equation}
where the element $\overline{\eta}^{i\stackrel{w}{\cap} j} \in \overline{\mathcal{E}}$ is defined by 
\begin{equation} 
\li(\overline{\eta}^{i\stackrel{w}{\cap} j}\re)^{\mbox{\tiny $J$}} :=  \overline{\eta}^{\mbox{\tiny $J$}} \;\mbox{ and }\; \li(\overline{\eta}^{i \stackrel{w}{\cap} j}\re)^{\mbox{\tiny $I$}}_{l} := \li\{ \begin{array}{ccl} \overline{\eta}^{\mbox{\tiny $I$}}_{i} \cap \overline{\eta}^{\mbox{\tiny $I$}}_{j} &,& l = i \\  \{0,\dots,w\} &,& l = j \\ \overline{\eta}^{\mbox{\tiny $I$}}_{l} &,& l \not=  i,j \end{array} \re.\;.
\end{equation}
%-----------------------------------------------------------------------
\item \label{itemdualresamplingIIv}  ($\star$) For each $i,j \in \tilde{\Gamma}(\overline{\eta})$ with $i \not=j$ and each $v \in \overline{\eta}^{\mbox{\tiny $I$}}_{i} \cap \overline{\eta}^{\mbox{\tiny $I$}}_{j}$  with $v \not= \min \overline{\eta}^{\mbox{\tiny $I$}}_{i} \cap \overline{\eta}^{\mbox{\tiny $I$}}_{j}$,
\begin{equation}
\overline{\mathcal{K}}(\overline{\eta}, \overline{\eta}^{i \stackrel{v}{\cap} j})  = \mathbbm{1}\{\overline{\eta}^{\mbox{\tiny $I$}}_{i} \cap \overline{\eta}^{\mbox{\tiny $I$}}_{j}   \not\in \{\emptyset, K\}\} \frac{S}{2N}\li[\chi(v) - \chi(v^{<})\re]\;,
\end{equation}
where 
\begin{equation}
v^{<} := \max\{v' \in \overline{\eta}^{\mbox{\tiny $I$}}_{i} \cap \overline{\eta}^{\mbox{\tiny $I$}}_{j}:v' < v\}
\end{equation}
is the greatest element in the intersection of $\overline{\eta}^{\mbox{\tiny $I$}}_{i}$ and $\overline{\eta}^{\mbox{\tiny $I$}}_{j}$ that is less than $v$, and the element $\overline{\eta}^{i\stackrel{v}{\cap} j} \in \overline{\mathcal{E}}$ is defined by 
\begin{equation} 
\li(\overline{\eta}^{i\stackrel{v}{\cap} j}\re)^{\mbox{\tiny $J$}} :=  \overline{\eta}^{\mbox{\tiny $J$}} \mbox{ and } \li(\overline{\eta}^{i \stackrel{v}{\cap} j}\re)^{\mbox{\tiny $I$}}_{l} := \li\{ \begin{array}{ccl}\{v,\dots\}&,& l = i \\ K&,& l = j \\ \overline{\eta}^{\mbox{\tiny $I$}}_{l} &,& l \not=  i,j \end{array} \re. \;,
\end{equation} 
where
\begin{equation}
\{v,\dots\} := \{v' \in \overline{\eta}^{\mbox{\tiny $I$}}_{i} \cap \overline{\eta}^{\mbox{\tiny $I$}}_{j}:v' \ge v\}
\end{equation}
is the set of elements which are in the intersection of $\overline{\eta}^{\mbox{\tiny $I$}}_{i}$ and $\overline{\eta}^{\mbox{\tiny $I$}}_{j}$ and which are equal or greater than $v$. 
\end{enumerate}
\end{enumerate}
\end{enumerate}
\end{Definition}
%%%%%%%%%%%%%%%%%%%%%%%%%%%%%%%%%%%%%%%%%%%%%%%%%%%%%%%%%%%%%%%%%%%%%%%%%%%%%%%%%%%%%%%%%%%%%%%%%%%      Formulation of the models, The BP, Annotations

As already announced we now give additional explanations (divided into five items) that are essential for the understanding of the BP.
\begin{enumerate}
%-----------------------------------------------------------------------
\item {\bf (Reversal of mutation)}\label{itemAnnotation1} The description in \ref{itemdualmutationJ} says that  the mark of a partition element in the BP  mutates  from $v$ to $u$  at the same rate as the type of a live-site in the HMM mutates from $u$ to $v$. Since the matrix $b(\cdot, \cdot)$ is in general not doubly stochastic, this way of reversing the mutation mechanism gives rise to a Feynman-Kac term in the duality relation depending on $b(\cdot, \cdot)$.
%-----------------------------------------------------------------------
\item  {\bf (The case $S = 0$, without selection)}\label{itemAnnotation2} The transitions rates specified in  \ref{itemdualresamplingJJw}, \ref{itemdualresamplingJIw}, \ref{itemdualresamplingIJw}, \ref{itemdualresamplingIIw} and \ref{itemdualresamplingIIv} are equal to $0$ and hence the second component of the BP is equal to $\overline{\xi^{*}}^{\mbox{\tiny $I$}}$ from (\ref{equationoverlinexi}) for all times. This in turn implies that the transitions  in \ref{itemdualmutationIcup} and \ref{itemdualmutationIsetminus} do not occur, the transition rate  \ref{itemdualresamplingIIK} is equal to $0$ and the transition in \ref{itemdualresamplingIJK} does not change the state. This means that the BP is a pure jump process on
\begin{equation}
\li\{\overline{\eta}  \in  \overline{\mathcal{E}}: \overline{\eta}^{\mbox{\tiny $I$}}_{i} = K \; \mbox{ for all } \; i \in I\re\}
\end{equation}
whose transitions are given by \ref{itemdualmutationJ}, \ref{itemdualresamplingJJK} and \ref{itemdualresamplingJIK}. 

In other words, each $K$-marked partition element undergoes  a random walk on $I$ according to  \ref{itemdualresamplingJIK}. It instantaneously coalesce (see \ref{itemdualresamplingJJK}) if it jumps to a life that is occupied by a partition element with the same mark, where in this case the set $K$ is assigned to the emerging active life-site which, however, leaves the second component of the BP unchanged due to the initial state given in (\ref{equationoverlinexi}). In addition, the mark of each partition element itself evolves as a Markov chain on $K$ according to \ref{itemdualmutationJ}.

Furthermore, if one takes out the information concerning the life-sites (this can be done when one only considers ancestral lines and genealogical distances in the case where the type information is exchangeable), then one obtains a $K$-marked $J$-partition  whose partition elements coalesce at rate $1$ if their marks coincide as required in \ref{itemdualresamplingJJK}, where again each mark evolves as a Markov chain on $K$ according to \ref{itemdualmutationJ}.
%-----------------------------------------------------------------------
\item{\bf (The case $S > 0$, with selection)}\label{itemAnnotation3} The additional coalescing event in \ref{itemdualresamplingJJw} (also the event in \ref{itemdualresamplingJIw} and in \ref{itemdualresamplingIJw}) generates a {\it strict and non-empty subset of $K$} that is assigned to the emerging active life-site. This has the effect that on the one hand the jump of a partition element to a life-site $i$ depends on $\overline{\eta}^{\mbox{\tiny $I$}}_{i}$ (consider the rate in \ref{itemdualresamplingJIw}) and on the other hand the active live-sites evolve according to  \ref{itemdualmutationIcup} and  \ref{itemdualmutationIsetminus} and interact with each other according to \ref{itemdualresamplingIIK}, \ref{itemdualresamplingIIw} and \ref{itemdualresamplingIIv}, where the transition \ref{itemdualresamplingIIv} can only occur if $|K| > 2$. Moreover note that the BP in fact takes values in (\ref{equationOverlineMathcalE}) since by definition the occurrence of empty sets is prevented in \ref{itemdualmutationIsetminus} \ref{itemdualresamplingIIK}, \ref{itemdualresamplingIIw} and \ref{itemdualresamplingIIv}.
 
The key idea for these additional events including the occurrence of strict and non-empty subsets in the second component of the BP comes from the fact that we reverse the resampling events. Consider (in the HMM) the conditional rate  of the event  
\begin{equation}
\li\{ \begin{array}{c} \mbox{(type of $i$ at time $t$)} = u \\ \mbox{(type of $j$ at time $t$)} =  u \end{array} \re\}
\end{equation} 
given both the information that $i$ replaces the type of $j$ at time $t$ and the types of the life-sites $i$ and $j$  just before this resampling event has occurred. A heuristic argument yields that this conditional rate is equal to
\begin{equation}
\sum_{v \in K}  \li(\frac{1}{2} + \frac{S}{2N}\li[\chi(u) - \chi(v)\re]\re) \mathbbm{1}\li\{ \begin{array}{c}\mbox{(type of $i$ at time $t-$)} = u \\ \mbox{(type of $j$ at time $t-$)} =  v\end{array} \re\}
\end{equation} 
which can be reformulated (recall (\ref{equationDefChi})) as
\begin{eqnarray}
&&\li(\frac{1}{2} + \frac{S}{2N}\li[\chi(u) - 1\re]\re) \mathbbm{1}\li\{ \begin{array}{c}\mbox{(type of $i$ at time $t-$)} = u \\ \mbox{(type of $j$ at time $t-$)} \in K\end{array} \re\}\\
&+&\sum_{w = 0}^{d-2}\frac{S}{2N}\li[\chi(w+1) - \chi(w)\re] \mathbbm{1}\li\{ \begin{array}{l}\mbox{(type of $i$ at time $t-$)} = u\\ \mbox{(type of $j$ at time $t-$)} \in \{0,\dots, w\}\end{array} \re\}
\end{eqnarray}
which leads to the rates and the transitions in \ref{itemdualresamplingJJK} and  \ref{itemdualresamplingJJw}.
%-----------------------------------------------------------------------
\item{\bf (The case $K = \{0,1\}$, two types)}\label{itemAnnotation4} If $K = \{0,1\}$, then the transitions in \ref{itemdualresamplingJJw}, in \ref{itemdualresamplingJIw} and also in \ref{itemdualresamplingIIw} generate the subset $\{0\}$, the transition in \ref{itemdualresamplingIIv} does not occur since $|K| = 2$. In addition,  \ref{itemdualmutationIcup} leads to a transition from $\{0\}$ to $K$ and \ref{itemdualmutationIsetminus} does not occur. Hence, in the two type case the BP with initial state  (\ref{equationoverlinexi}) takes values in
\begin{equation}
\li\{\overline{\eta} \in \overline{\mathcal{E}}: \overline{\eta}^{\mbox{\tiny $I$}}_{i} \in \{ K, \{0\}\} \mbox{ for all } i \in I \re\} 
\end{equation}
which will be important in Subsection \ref{subsecApplications}.
%-----------------------------------------------------------------------
\item{\bf (Ancestral lines and genealogical distances)}\label{itemAnnotation5} When we are interested in ancestral lines and genealogical distances in the case where the type information is exchangeable (this will be the case in Subsection \ref{subsecApplications}), then the BP {\it reduces} as follows (see also the brief discussion at the end of item \ref{itemAnnotation2}). In the first component we only need the partition elements together with their marks in $K$. In the second component we only have to consider  for each strict and non-empty subset of $K$ the number of active life sites to which this subset is assigned.  Observe that this reduction of the BP, formally a functional of the BP, is again a Markov process.

So, in the two type case (see item \ref{itemAnnotation4} above), for the second component we only have to consider
\begin{equation}\label{equationAnzahlNuller}
\overline{\eta}^{\mbox{\tiny $I$}}_{|\{0\}|} := \li|\li\{i \in  \Gamma(\overline{\eta}) : \overline{\eta}^{\mbox{\tiny $I$}}_{i}= \{0\}\re\}\re|\;,
\end{equation}
the number of active life sites to which the subset $\{0\}$ is assigned.
\end{enumerate}
%%%%%%%%%%%%%%%%%%%%%%%%%%%%%%%%%%%%%%%%%%%%%%%%%%%%%%%%%%%%%%%%%%%%%%%%%%%%%%%%%%%%%%%%%%%%%%%%%%%      Formulation of the models, The BP,overview transitions

\begin{figure}[h]\caption{Transitions of the BP}
\label{FigureTransitionsBP}
\begin{center}
\begin{tabular}{|l|c|c|}\hline
& a transition occur for each &  and changes the state of\\ \hline 
\ref{itemdualmutationJ} &  $\gamma \in \Gamma(\overline{\eta})$ & $\overline{\eta}^{\mbox{\tiny $J$}}_{\gamma,K}$\\ \hline
\ref{itemdualmutationIcup} &  $i \in \tilde{\Gamma}(\overline{\eta})$, $v \in \overline{\eta}^{\mbox{\tiny $I$}}_{i}$ and  $u \in K \setminus \overline{\eta}^{\mbox{\tiny $I$}}_{i}$ & $\overline{\eta}^{\mbox{\tiny $I$}}_{i}$\\ \hline
\ref{itemdualmutationIsetminus} & $i \in \tilde{\Gamma}(\overline{\eta})$,  $v \in K\setminus\overline{\eta}^{\mbox{\tiny $I$}}_{i}$ and  $u \in \overline{\eta}^{\mbox{\tiny $I$}}_{i}$&  $\overline{\eta}^{\mbox{\tiny $I$}}_{i}$\\ \hline
\ref{itemdualresamplingJJK} & $\gamma\not= \gamma' \in  \Gamma(\overline{\eta})$ & $\overline{\eta}^{\mbox{\tiny $J$}}_{\gamma}$ and $\overline{\eta}^{\mbox{\tiny $I$}}_{l}$ with $l =  \overline{\eta}^{\mbox{\tiny $J$}}_{\gamma,I}$\\ \hline
\ref{itemdualresamplingJJw} & $\gamma \not= \gamma' \in  \Gamma(\overline{\eta})$ and $w \in \{0,\dots, d-2\}$ & $\overline{\eta}^{\mbox{\tiny $J$}}_{\gamma}$ and $\overline{\eta}^{\mbox{\tiny $I$}}_{l}$ with $l =  \overline{\eta}^{\mbox{\tiny $J$}}_{\gamma,I}$\\ \hline
\ref{itemdualresamplingJIK} & $\gamma \in  \Gamma(\overline{\eta})$ and $i \in \tilde{\Gamma}(\overline{\eta})$ & $\overline{\eta}^{\mbox{\tiny $J$}}_{\gamma}$ and $\overline{\eta}^{\mbox{\tiny $I$}}_{l}$ with $l =  \overline{\eta}^{\mbox{\tiny $J$}}_{\gamma,I}$\\ \hline
\ref{itemdualresamplingJIw} & $\gamma \in  \Gamma(\overline{\eta})$, $i \in \tilde{\Gamma}(\overline{\eta})$ and $w \in \{0,\dots, d-2\}$ & $\overline{\eta}^{\mbox{\tiny $J$}}_{\gamma}$ and $\overline{\eta}^{\mbox{\tiny $I$}}_{l}$ with $l =  \overline{\eta}^{\mbox{\tiny $J$}}_{\gamma,I}$\\ \hline
\ref{itemdualresamplingIJK} & $i \in \tilde{\Gamma}(\overline{\eta})$ and $\gamma \in  \Gamma(\overline{\eta})$&  $\overline{\eta}^{\mbox{\tiny $I$}}_{i}$\\ \hline
\ref{itemdualresamplingIJw} & $i \in \tilde{\Gamma}(\overline{\eta})$, $\gamma \in  \Gamma(\overline{\eta})$  and $w \in \{0,\dots, d-2\}$ &  $\overline{\eta}^{\mbox{\tiny $I$}}_{i}$\\ \hline
\ref{itemdualresamplingIIK} &  $i \not=j \in \tilde{\Gamma}(\overline{\eta})$ & $\overline{\eta}^{\mbox{\tiny $I$}}_{i}$ and $\overline{\eta}^{\mbox{\tiny $I$}}_{j}$ \\ \hline
\ref{itemdualresamplingIIw} &  $i \not=j \in \tilde{\Gamma}(\overline{\eta})$  and $w \in \{0,\dots, d-2\}$ & $\overline{\eta}^{\mbox{\tiny $I$}}_{i}$ and $\overline{\eta}^{\mbox{\tiny $I$}}_{j}$ \\ \hline
\ref{itemdualresamplingIIv} &  $i \not=j \in \tilde{\Gamma}(\overline{\eta})$ and  $v \in \overline{\eta}^{\mbox{\tiny $I$}}_{i} \cap \overline{\eta}^{\mbox{\tiny $I$}}_{j}: v\not= \min \overline{\eta}^{\mbox{\tiny $I$}}_{i} \cap \overline{\eta}^{\mbox{\tiny $I$}}_{j}$ & $\overline{\eta}^{\mbox{\tiny $I$}}_{i}$ and $\overline{\eta}^{\mbox{\tiny $I$}}_{j}$ \\ \hline
\end{tabular} 
\end{center}
\end{figure}
%%%%%%%%%%%%%%%%%%%%%%%%%%%%%%%%%%%%%%%%%%%%%%%%%%%%%%%%%%%%%%%%%%%%%%%%%%%%%%%%%%%%%%%%%%%%%%%%%%%      Formulation of the models, Analytical characterization

\subsection{Analytical characterization}\label{subsecAnalyticalCharacterization}
\noindent In this subsection we first describe the HMM by a Markov process that arises as the unique solution of a well-posed martingale problem. This gives an analytic characterization of the HMM and is the main tool to prove the duality with help of the generator. 

Then we characterize for $J \subset I$ the BP by a Markov process that uniquely corresponds to a bounded linear operator.\bigskip

{\it The HMM:} Let $\Omega = \mathcal{D}([0,\infty),\mathcal{E})$ be the space of cadlag paths from $[0,\infty)$ to $\mathcal{E}$, $\mathcal{F} =  \mathcal{B}(\Omega)$ be the Borel $\sigma$-field on $\Omega$ and $(X_{t})_{t \ge 0}$ be the canonical $\mathcal{E}$-valued process on $(\Omega, \mathcal{F})$ with canonical right continuous filtration $(\mathcal{F}_{t})_{t \ge 0}$. In order to characterize the HMM analytically we specify a measure determining set
\begin{equation}
\mathcal{A} \subset C_{b}(\mathcal{E}), 
\end{equation}
define a linear operator 
\begin{equation}
L:\mathcal{A}  \to  C_{b}(\mathcal{E})
\end{equation}
and use the following concept, the martingale problem.
\begin{Definition}[Solution and well-posedness]\pointpar \label{DefinitionMartingaleproblem} \noindent Let $\mu$ be a distribution on $\mathcal{E}$. A probability measure $\mathbb{Q}$ on $\Omega$ is a solution of the $\Omega$-martingale problem for $(L,\mu)$ w.r.t$.$ $\mathcal{A}$ if and only if
\begin{equation}\label{equationMartingalproblem1}
\li(f(X_{t}) - \int_{0}^{t}Lf(X_{s})ds\re)_{t \ge 0} \; \mbox{ is a martingale under $\mathbb{Q}$ for all $f \in \mathcal{A}$}
\end{equation}
and
\begin{equation}\label{equationMartingalproblem2}
\mathbb{Q} \circ (X_{0})^{-1} = \mu\;.
\end{equation}
The $\Omega$-martingale problem for $(L,\mu)$ w.r.t$.$ $\mathcal{A}$ is called well-posed if $\mathbb{Q}$ is uniquely determined by (\ref{equationMartingalproblem1}) and (\ref{equationMartingalproblem2}). 

Moreover, the $\Omega$-martingale problem for $L$ w.r.t$.$ $\mathcal{A}$ is called well-posed if the $\Omega$-martingale problem for $(L,\mu)$ w.r.t$.$ $\mathcal{A}$ is well-posed for all distributions $\mu$ on $\mathcal{E}$.
\end{Definition}

As measure determining set $\mathcal{A}$, see also \cite{D93}, we take the algebra that is generated by functions of the form 
\begin{equation}\label{equationFunctionsInA}
f(\eta) = g(\eta^{\mbox{\tiny \sc Time}})\prod_{j \in J}\mathbbm{1}\{\eta^{*}_{j} = u_{j}\}\prod_{n = 1}^{m_{j}}\int_{r_{j}^{n}}^{t_{j}^{n}}F_{j}^{n}(\eta^{\mbox{\tiny $\mathcal{D}$}}_{j,s})ds
\end{equation}
where 
\begin{equation}
g \in C^{1}_{b}(\mathbb{R}),\; J \subset I, \; u_{j} \in K, \; r_{j}^{1} < t_{j}^{1} < \dots < r_{j}^{m_{j}} < t_{j}^{m_{j}} \mbox{ and  } F_{j}^{n} \in C_{b}(K \times I)\;.
\end{equation}

\begin{Remark}\point There is an important reason why we consider functions of the form (\ref{equationFunctionsInA}) instead of functions of the form
\begin{equation}\label{equationFunctionsNotInA}
\eta \mapsto g(\eta^{\mbox{\tiny \sc Time}})\prod_{j \in J}\mathbbm{1}\{\eta^{*}_{j} = u_{j}\}\prod_{n = 1}^{m_{j}}F_{j}^{n}(\eta^{\mbox{\tiny $\mathcal{D}$}}_{j,t_{j}^{n}})\;.
\end{equation}
Namely, in contrast to a function of the form (\ref{equationFunctionsNotInA}) which evaluates extended ancestral lines at different {\sc Time} points, a function of the form (\ref{equationFunctionsInA}) is an element in $C_{b}(\mathcal{E})$ since the extended ancestral lines are evaluated in terms of integrals over different {\sc Time} intervals, where one should have in mind that $(\cdot)^{*}$, the projection on the current types, is continuous due to the special form of $\mathcal{E}$ as mentioned before. 

Furthermore, the special form (\ref{equationFunctionsInA}) of the functions in $\mathcal{A}$ will be used to relate extended ancestral lines with the sample-paths of the BP by means of a duality function $H$ in Subsection \ref{subsecFeynmanKacHMMandHBP}.
\end{Remark} 

The linear operator $L= L^{\mbox{\tiny \sc Time}} + L^{\mbox{\tiny $\mathcal{D}$}}$ on $\mathcal{A}$ is defined by
\begin{equation}
L^{\mbox{\tiny \sc Time}}f(\eta) = \lim_{\epsilon \downarrow 0} \frac{f(\eta^{\mbox{\tiny \sc Time}} + \epsilon, \eta^{\mbox{\tiny $\mathcal{D}$}}) - f(\eta^{\mbox{\tiny \sc Time}}, \eta^{\mbox{\tiny $\mathcal{D}$}})}{\epsilon}
\end{equation}
and
\begin{eqnarray}
L^{\mbox{\tiny $\mathcal{D}$}}f(\eta) &=&  B\sum_{i \in I}\sum_{u \in K} b(\eta^{*}_{i},u)\li[f(\eta^{i;u})-f(\eta)\re]\\
 && +  \sum_{i,j \in I}\li(\frac{1}{2} + \frac{S}{2N}\li[\chi(\eta^{*}_{i}) - \chi(\eta^{*}_{j})\re]\re)\li[f(\eta^{i\to j})-f(\eta)\re]\;.
\end{eqnarray}
Observe that for $f \in \mathcal{A}$, $Lf$ is indeed {\it continuous} although the map $\eta \mapsto \eta^{i \to j}$ is not. The reason for this is once more the evaluation of the extended ancestral lines in terms of integrals over different {\sc Time} intervals.
%%%%%%%%%%%%%%%%%%%%%%%%%%%%%%%%%%%%%%%%%%%%%%%%%%%%%%%%%%%%%%%%%%%%%%%%%%%%%%%%%%%%%%%%%%%%%%%%%%%      Caution:  canonical filtration vs   canonical right continuous filtration  , see Patric !!!!!!!!!!!!!!!!!!!!!
%%%%%%%%%%%%%%%%%%%%%%%%%%%%%%%%%%%%%%%%%%%%%%%%%%%%%%%%%%%%%%%%%%%%%%%

\begin{Theorem}[Analytical characterization of the HMM]\point\label{TheoremAnalyticalHMM}
\begin{enumerate}[a)]
\item\label{TheoremAnalyticalHMMa} The $\Omega$-martingale problem for $L$ w.r.t$.$ $\mathcal{A}$ is well-posed. 
\item\label{TheoremAnalyticalHMMb}  If  the unique solution for $(L,\delta_{\eta})$ is denoted by $\mathbb{P}_{\eta}$, then 
\begin{equation}
 X =  \li( \Omega , \mathcal{F}, (\mathcal{F}_{t})_{t \ge 0},(X_{t})_{t \ge 0}, \{ \mathbb{P}_{\eta} : \eta \in \mathcal{E}\} \re) 
\end{equation}
is a Borel strong Markov process.
\item\label{TheoremAnalyticalHMMc}  For any $\mu \in \mathcal{M}_{1}(\mathcal{E})$ the probability measure 
\begin{equation}\label{equationPmu}
\mathbb{P}_{\mu}(\cdot) = \int_{\mathcal{E}} \mathbb{P}_{\eta}(\cdot) \mu(d\eta)
\end{equation}
is the unique solution for $(L,\mu)$.
\item\label{TheoremAnalyticalHMMd}   If $\mu$ is the initial state defined in (\ref{equationMu}), then $\mathbb{P}_{\mu}$ is equal to the law of the piecewise deterministic Markov jump process in Definition \ref{DefinitionHMM}.
\end{enumerate}
\end{Theorem}
Observe that throughout this paper we use 
\begin{equation}
\mbox{$\mathbb{E}_{\mu}$ to denote the expectation w.r.t$.$  $\mathbb{P}_{\mu}$}
\end{equation}
respectively $\mathbb{E}_{\eta}$ to denote the expectation w.r.t$.$ $\mathbb{P}_{\eta}$.

Finally, we consider two important functionals of the HMM:
\begin{enumerate}
\item For $t \ge 0$ and $i,j \in I$,
\begin{equation}\label{equationDt}
\mathbb{D}_{t}(i,j):= 2 \li|(X_{t})^{\mbox{\tiny {\sc Time}}} - \sup\li\{ s \in \li[(X_{0})^{\mbox{\tiny {\sc Time}}},(X_{t})^{\mbox{\tiny {\sc Time}}}\re]:  (X_{t})^{\mbox{\tiny $\mathcal{D}$}}_{i,s} = (X_{t})^{\mbox{\tiny $\mathcal{D}$}}_{j,s} \, \re\}\re|
\end{equation}
 is the genealogical distance of the life-sites $i$ and $j$ at time $t$, where we use the convention $\sup \emptyset := X^{\mbox{\tiny {\sc Time}}}_{0}$. Hence
\begin{equation}
\li(I, (\mathbb{D}_{t}(i,j))_{i,j \in I}, \frac{1}{N}\sum_{i = 1}^{N}\delta_{(i,(X_{t})^{*}_{i})}  \re)_{t \ge 0}
\end{equation} 
is a version of the tree-valued Moran model introduced in \cite{DGP12}.
\item The CAT of the Moran Model is a stochastic process on $K$ we denote by $(K_{t})_{t \ge 0}$.  In \cite{KHB13} the CAT at time $t$ is defined {\it as the type of the unique individual that is, at some time $s > t$, ancestral to the whole population}, where it is assumed that the type information of the population is in equilibrium. In contrast we can always explicitly define the CAT at time $t$ (as a functional of our forward evolving HMM) by
\begin{equation}\label{equationCAT}
K_{t} := \li\{u \in K: \exists s \ge t, j\in I \mbox{ such that } \li(X_{s}\re)^{\mathcal{D}}_{i,\li(X_{t}\re)^{\mbox{\tiny {\sc Time}}}} = \li(u,j\re) \mbox{ for all } i \in I \re\}\;.
\end{equation} 
Of course it is difficult to work with this representation. In order to obtain the stationary type distribution of the CAT we will consider the extended ancestral lines in equilibrium and trace a single ancestral line back to {\sc Time} $-\infty$ in Subsection \ref{subsecApplications}. Observe that this is in the spirit of \cite{F02}.\bigskip
\end{enumerate}
%%%%%%%%%%%%%%%%%%%%%%%%%%%%%%%%%%%%%%%%%%%%%%%%%%%%%%%%%%%%%%%%%%%%%%%%%%%%%%%%%%%%%%%%%%%%%%%%%%%      Formulation of the models, Analytical characterization

{\it The BP:} Let $\overline{\Omega} = \mathcal{D}([0,\infty),\overline{\mathcal{E}})$ be the space of cadlag paths from $[0,\infty)$ to $\overline{\mathcal{E}}$, $\overline{\mathcal{F}} =  \mathcal{B}(\overline{\Omega})$ be the Borel $\sigma$-field on $\overline{\Omega}$ and $(\overline{X}_{t})_{t \ge 0}$ be the canonical $\overline{\mathcal{E}}$-valued process on $(\overline{\Omega}, \overline{\mathcal{F}})$ with canonical right continuous filtration $(\overline{\mathcal{F}}_{t})_{t \ge 0}$ . In order to characterize the BP analytically we define, according to Definition \ref{DefinitionBP}, the map $\overline{L}: C_{b}(\overline{\mathcal{E}})  \to  C_{b}(\overline{\mathcal{E}})$ by 
\begin{equation}
\overline{L}f(\overline{\eta}) = \sum_{\overline{\zeta} \in \overline{\mathcal{E}}} \overline{\mathcal{K}}(\overline{\eta}, \overline{\zeta})\li[f(\overline{\zeta}) - f(\overline{\eta})\re]\;.
\end{equation}
%%%%%%%%%%%%%%%%%%%%%%%%%%%%%%%%%%%%%%%%%%%%%%%%%%%%%%%%%%%%%%%%%%%%%%%%%%%%%%%%%%%%%%%%%%%%%%%%%%%      Formulation of the models, Analytical characterization, dual

Since $\overline{L}$ is a bounded linear operator, the family $\{\exp(t\overline{L}\,\cdot\,): t \ge 0\}$ of linear operators on $C_{b}(\overline{\mathcal{E}})$ is a strongly continuous, positive, contraction semigroup, and the $\overline{\Omega}$-martingale problem for $\overline{L}$ w.r.t$.$ $C_{b}(\overline{\mathcal{E}})$ is well-posed. Hence
\begin{equation}\label{equationOverlineX}
 \overline{X} =  \li( \overline{\Omega} , \overline{\mathcal{F}}, (\overline{\mathcal{F}}_{t})_{t \ge 0},(\overline{X}_{t})_{t \ge 0}, \{\overline{\mathbb{P}}_{\overline{\eta}} : \overline{\eta} \in \overline{\mathcal{E}}\}\re) 
\end{equation}
is a Borel strong Markov process if $\overline{\mathbb{P}}_{\overline{\eta}}$ denotes the unique solution for $(\overline{L},\delta_{\overline{\eta}})$. In analogy to the HMM we write $\overline{\mathbb{E}}_{\overline{\eta}}$ if we take the expectation w.r.t$.$ $\overline{\mathbb{P}}_{\overline{\eta}}$. Finally, the probability measure $\overline{\mathbb{P}}_{\overline{\xi^{*}}}$ is equal to the law of the Markov jump process in Definition \ref{DefinitionBP} if  $\overline{\xi^{*}}$ is the initial state defined in (\ref{equationoverlinexi}).
%%%%%%%%%%%%%%%%%%%%%%%%%%%%%%%%%%%%%%%%%%%%%%%%%%%%%%%%%%%%%%%%%%%%%%%%%%%%%%%%%%%%%%%%%%%%%%%%%%%      Formulation of the main results

\setcounter{equation}{0}
\section{Formulation of the results}\label{secFormulationresults}
\noindent  This section is organized as follows. In Subsection \ref{subsecDualityTypInfoHMM} we state the Feynman-Kac duality between the HMM and the BP with which we can represent the type information of a tagged $J \subset I$ in the HMM in terms of the BP. This  Feynman-Kac duality serves as a warm up for the results concerning the relation between the extended ancestral lines of $J$ and the sample paths of the BP presented in Subsection \ref{subsecHMMsamplepathsBP}. Subsection \ref{subsecErgodicity} contains the limit theorem for the extended ancestral lines which allows us to study both the stationary type distribution of the {\it CAT} and the {\it conditioned genealogical distance} of two individuals in Subsection \ref{subsecApplications}.
%%%%%%%%%%%%%%%%%%%%%%%%%%%%%%%%%%%%%%%%%%%%%%%%%%%%%%%%%%%%%%%%%%%%%%%%%%%%%%%%%%%%%%%%%%%%%%%%%%%      results,Feynman-Kac HMM and BP

\subsection{Feynman-Kac duality between the HMM and the BP}\label{subsecDualityTypInfoHMM}
\noindent Here we show how we can use the BP to determine the probability that the life-sites in $J$ have given {\it types} at a given time. The tool is a Feynman-Kac duality, a duality with respect to a Feynman-Kac function
\begin{equation}
V: \overline{\mathcal{E}} \to \mathbb{R}
\end{equation}
and a duality function 
\begin{equation}
H^{*}:\mathcal{E} \times \overline{\mathcal{E}} \to [0,1]\; , 
\end{equation}
where for $\overline{\eta} \in \overline{\mathcal{E}}$ the value of $H^{*}(\eta, \overline{\eta})$ only depends on $\eta^{*}$ (see (\ref{equationProjectionOnTypes})).

Namely, the duality function is given by
\begin{equation}\label{equationHstar}
H^{*}(\eta, \overline{\eta})  =  \li(\prod_{\gamma \in \Gamma(\overline{\eta})} \mathbbm{1}\{ \eta^{*}_{\overline{\eta}^{\mbox{\tiny $J$}}_{\gamma,I}} = \overline{\eta}^{\mbox{\tiny $J$}}_{\gamma,K}\}\re) \li(\prod_{i \in \tilde{\Gamma}(\overline{\eta})} \mathbbm{1}\{\eta^{*}_{i} \in \overline{\eta}^{\mbox{\tiny $I$}}_{i}\}\re)\;.
\end{equation}
and can be used to consider the types of $J$ in the HMM at a given time $t$ by
\begin{equation}
H^{*}(X_{t}, \overline{\xi^{*}}) = \prod_{j \in J} \mathbbm{1}\{ (X_{t})^{*}_{j} = \xi^{*}_{j}\} = \mathbbm{1}\{ \li((X_{t})^{*}_{j}\re)_{j \in J} = \xi^{*}\}
\end{equation}
if $\overline{\xi^{*}} \in \overline{\mathcal{E}}$ is defined as in (\ref{equationoverlinexi}). 
%%%%%%%%%%%%%%%%%%%%%%%%%%%%%%%%%%%%%%%%%%%%%%%%%%%%%%%%%%%%%%%%%%%%%%%%%%%%%%%%%%%%%%%%%%%%%%%%%%%      Results, Feynman-Kac HMM and BP
 
The Feynman-Kac function has in the neutral case (recall the items \ref{itemAnnotation1} and \ref{itemAnnotation2} subsequent to Definition \ref{DefinitionBP}) the form
\begin{equation}
V(\overline{\eta})
= B\sum_{\gamma \in \Gamma(\overline{\eta}) }\li(\sum_{u \in K} b(u,\overline{\eta}^{\mbox{\tiny $J$}}_{\gamma,K})-1\re) -  \frac{1}{2}\sum_{\gamma \not= \gamma' \in \Gamma(\overline{\eta}) } \mathbbm{1}\{\overline{\eta}^{\mbox{\tiny $J$}}_{\gamma,K} \not= \overline{\eta}^{\mbox{\tiny $J$}}_{\gamma',K}\}\;,
\end{equation}
where the first term arises due to the fact that $b(\cdot, \cdot)$ is in general not doubly stochastic and the second term,  which is due to the "having same type" requirement in \ref{itemdualresamplingJJK}, is the difference between the total interaction rate in the BP and total resampling rate in the HMM.

In the case with selection we again have the first term and the difference, but also additional effects appear due to the restrictions $\mathbbm{1}\{|\overline{\eta}^{\mbox{\tiny $I$}}_{i}| > 1\}$ in \ref{itemdualmutationIsetminus} and $\mathbbm{1}\{\overline{\eta}^{\mbox{\tiny $I$}}_{i} \cap \overline{\eta}^{\mbox{\tiny $I$}}_{j} \not\in \{\emptyset, K\}\}$ in \ref{itemdualresamplingIIv}. So, we set 
\begin{eqnarray}\label{equationV}
V(\overline{\eta})
&=& B\sum_{\gamma \in \Gamma(\overline{\eta}) }\li(\sum_{u \in K} b(u,\overline{\eta}^{\mbox{\tiny $J$}}_{\gamma,K})-1\re) - B \sum_{i \in \tilde{\Gamma}(\overline{\eta})}\mathbbm{1}\{|\overline{\eta}^{\mbox{\tiny $I$}}_{i}| = 1\}\sum_{u \in \overline{\eta}^{\mbox{\tiny $I$}}_{i}}\sum_{v \in K \setminus \{u\}}b(u,v)\\
&& +\sum_{\gamma \not= \gamma' \in \Gamma(\overline{\eta}) } \li(\mathbbm{1}\{\overline{\eta}^{\mbox{\tiny $J$}}_{\gamma,K} = \overline{\eta}^{\mbox{\tiny $J$}}_{\gamma',K}\}\li[\frac{1}{2} + \frac{S}{2N}\chi(\overline{\eta}^{\mbox{\tiny $J$}}_{\gamma,K})\re] - \frac{1}{2}\re)\\
&&+2\sum_{\gamma \in \Gamma(\overline{\eta}) ,  i \in \tilde{\Gamma}(\overline{\eta})} \li(\mathbbm{1}\{\overline{\eta}^{\mbox{\tiny $J$}}_{\gamma,K} \in \overline{\eta}^{\mbox{\tiny $I$}}_{i}\}\li[\frac{1}{2} + \frac{S}{2N}\chi(\overline{\eta}^{\mbox{\tiny $J$}}_{\gamma,K})\re]- \frac{1}{2}\re)\\
&&+\sum_{i \not= j \in \tilde{\Gamma}(\overline{\eta})} \li(\mathbbm{1}\{\overline{\eta}^{\mbox{\tiny $I$}}_{i} \cap \overline{\eta}^{\mbox{\tiny $I$}}_{j} \not= K\}\frac{S}{2N}\chi(\max \overline{\eta}^{\mbox{\tiny $I$}}_{i} \cap \overline{\eta}^{\mbox{\tiny $I$}}_{j}) -  \frac{1}{2}\mathbbm{1}\{\overline{\eta}^{\mbox{\tiny $I$}}_{i} \cap \overline{\eta}^{\mbox{\tiny $I$}}_{j} = \emptyset\}\re)
 \end{eqnarray}
and have completed our preparations. 
%%%%%%%%%%%%%%%%%%%%%%%%%%%%%%%%%%%%%%%%%%%%%%%%%%%%%%%%%%%%%%%%%%%%%%%%%%%%%%%%%%%%%%%%%%%%%%%%%%%      Results, Feynman-Kac HMM and BP, Proposition

The Feynman-Kac duality now reads as follows.
\begin{Proposition}[Feynman-Kac duality for the type information of the HMM]\label{PropositionDualityTypInfoHMM} \pointpar
\noindent Let $\mu$ be as in (\ref{equationMu}) with $c = 0$ (initial {\sc Time}) and general initial type distribution $\mu^{*}$, $V$ as in (\ref{equationV}) and $H^{*}$ as in (\ref{equationHstar}). 

Then
\begin{equation}
\mathbb{E}_{\mu}\li[H^{*}(X_{t}, \overline{\eta})\re] = \overline{\mathbb{E}}_{\overline{\eta}}\li[\mathbb{E}_{\mu}\li[H^{*}(X_{0}, \overline{X}_{t})\re]e^{\int_{0}^{t}V(\overline{X}_{s})ds} \re]
\end{equation}
for all $\overline{\eta} \in \overline{\mathcal{E}}$ and all $t \ge 0$, where $\mathbb{E}_{\mu}$ denotes the expectation w.r.t$.$ $\mathbb{P}_{\mu}$ defined in (\ref{equationPmu}).
\end{Proposition}
%%%%%%%%%%%%%%%%%%%%%%%%%%%%%%%%%%%%%%%%%%%%%%%%%%%%%%%%%%%%%%%%%%%%%%%%%%%%%%%%%%%%%%%%%%%%%%%%%%%      Results, Feynman-Kac HMM and BP, Proposition, Remarks

\begin{Remark}\point 
We have that
\begin{equation}
H^{*}(X_{0}, \overline{X}_{t}) = \li(\prod_{\gamma \in \Gamma(\overline{X}_{t})} \mathbbm{1}\{ (X_{0})^{*}_{(\overline{X}_{t})^{\mbox{\tiny $J$}}_{\gamma,I}} = (\overline{X}_{t})^{\mbox{\tiny $J$}}_{\gamma,K}\} \re) \li(\prod_{i \in \tilde{\Gamma}(\overline{X}_{t})} \mathbbm{1}\{ (X_{0})^{*}_{i} \in (\overline{X}_{t})^{\mbox{\tiny $I$}}_{i}\}\re)\;.
\end{equation}
This first part on the r.h.s$.$ of this equation tells us that the types of the partition elements in the first component of the BP at time $t$ represent the types of the ancestors of $J$ in the HMM alive at time $0$. This means each life-site $j \in J$ in the first component of the BP has to find his way back to the type and the life-site of its ancestor at time $0$.  In the second part it is checked whether the states of the active life-sites at time $t$ (coding the type information of dead individuals in the HMM) fit with the types in the HMM a time $0$.
\end{Remark}
\begin{Remark}\label{RemarkExchangeableTypes}\point The formula in Proposition \ref{PropositionDualityTypInfoHMM} and the behaviour of the BP imply that the type information of $J$ is exchangeable at each time $t  \ge 0$ if $\mu^{*}$ is exchangeable. 

Formally, if $\mathbb{P}_{\mu}((X_{0})^{*} = \zeta^{*}) = \mathbb{P}_{\mu}(((X_{0})^{*}_{\sigma(i)})_{i \in I} = \zeta^{*})$ for all $\zeta^{*} \in K^{I}$ and all bijections $\sigma:I \to I$, then $\mathbb{P}_{\mu}(((X_{t})^{*}_{j})_{j \in J} = \xi^{*}) = \mathbb{P}_{\mu}(((X_{t})^{*}_{\sigma(j)})_{j \in J} = \xi^{*})$ for all $\xi^{*} \in K^{J}$ and all bijections $\sigma:I \to I$. 
\end{Remark}
%%%%%%%%%%%%%%%%%%%%%%%%%%%%%%%%%%%%%%%%%%%%%%%%%%%%%%%%%%%%%%%%%%%%%%%%%%%%%%%%%%%%%%%%%%%%%%%%%%%      Formulation of the main results, 
\subsection{Relation between the HMM and the sample paths of the BP} \label{subsecHMMsamplepathsBP}
\noindent This subsection covers the main result, the strong stochastic representation for the conditioned extended ancestral lines of $J$ in terms of the  sample paths of a transformation of the BP,  which is carried out in the Subsubsections \ref{subsubStochasticRepT} - \ref{subsubStrongduality}. Namely, we first give a stochastic representation with which we can express the expectation of certain functionals of the extended ancestral lines of $J$ at a fixed time $T$ in terms of the expectation of suitable functionals of the sample paths of the BP up to this time $T$. Based on this stochastic representation we then transform the BP and use the new object to obtain the {\it strong} stochastic representation for the conditioned extended ancestral lines of $J$ alive at time $T$.

Throughout this subsection we assume that $\mu$, the initial distribution of the HMM, is defined as in (\ref{equationMu}) with $c = -T$ (initial {\sc Time}) and general initial type distribution $\mu^{*}$. So, under this $\mu$ the extended ancestral lines at time $T$ describe the situation where the extended ancestral lines are considered from {\sc Time} $0$  back to {\sc Time} $-T$. In addition, we need an element $\xi^{*} \in K^{J}$ to describe the type information of $J$ at time $T$, that is, the type information of $J$ at {\sc Time} $0$.
%%%%%%%%%%%%%%%%%%%%%%%%%%%%%%%%%%%%%%%%%%%%%%%%%%%%%%%%%%%%%%%%%%%%%%%%%%%%%%%%%%%%%%%%%%%%%%%%%%%      Formulation of the main results, Duality

\subsubsection{Stochastic representation for extended ancestral lines} \label{subsubStochasticRepT}
\noindent This subsubsection includes the first step for the strong stochastic representation. We want to express certain information of the extended ancestral lines of $J$ between {\sc Time} $-T$  and {\sc Time} $0$ in terms of the sample paths of the BP up to time $T$.

We already know (see Proposition \ref{PropositionDualityTypInfoHMM}) that
\begin{equation}
\mathbb{P}_{\mu}\li( \li((X_{T})^{*}_{j}\re)_{j \in J} = \xi^{*}\re) = \overline{\mathbb{E}}_{\overline{\xi^{*}}}\li[ \mathbb{E}_{\mu}\li[H^{*}(X_{0}, \overline{X}_{T})\re]e^{\int_{0}^{T}V(\overline{X}_{s})ds} \re]\;,
\end{equation}
where $\li(X_{T}\re)^{*}_{j}$ is the type of $j$ at {\sc Time} $0$. We shall state that this equation for the types of $J$ can be refined to one for the extended ancestral lines of $J$. 

In order to determine the distribution of the extended ancestral lines of $J$ between {\sc Time} $-T$  and {\sc Time} $0$ restricted that the type information of $J$ at {\sc Time} $0$ is equal to $\xi^{*}$, we will determine the expectation of functionals of the form (\ref{equationFunctionsInA}) in which the extended ancestral lines of $J$ are evaluated in terms of integrals over different {\sc Time} intervals in $[-T,0]$. This means that for $-T \le r_{j}^{1} < t_{j}^{1} < \dots < r_{j}^{m_{j}} < t_{j}^{m_{j}} \le 0$ and  $F_{j}^{n} \in C_{b}(K \times I)$ we represent the expectation
\begin{equation}\label{equationFunctionalHMM}
\mathbb{E}_{\mu}\li[\li(\prod_{j \in J}\prod_{n = 1}^{m_{j}}\int_{r_{j}^{n}}^{t_{j}^{n}}F_{j}^{n}((X_{T})^{\mbox{\tiny $\mathcal{D}$}}_{j,s})ds\re)\mathbbm{1}\{\li((X_{T})^{*}_{j}\re)_{j \in J} = \xi^{*}\}\re]
\end{equation}
by the expectation of a suitable functional, a functional that also depends on the parameters in (\ref{equationFunctionalHMM}), of the sample paths of the  BP up to time $T$. 
%%%%%%%%%%%%%%%%%%%%%%%%%%%%%%%%%%%%%%%%%%%%%%%%%%%%%%%%%%%%%%%%%%%%%%%%%%%%%%%%%%%%%%%%%%%%%%%%%%%      results, Stochastic representation

\newpage
\begin{Theorem}[Stochastic representation at time $T$]\label{TheoremStochasticRepresentation} \pointpar
\noindent Let $J \subset I$, $T > 0$, $\xi^{*} \in K^{J}$, $\overline{\xi^{*}} \in \overline{\mathcal{E}}$ be defined as in (\ref{equationoverlinexi}) and $\mu$ as in (\ref{equationMu}) with $c = -T$ and general $\mu^{*}$. In addition, let $V$ be defined as in (\ref{equationV}) and $H^{*}$ as in (\ref{equationHstar}).

Then
\begin{equation}
(\ref{equationFunctionalHMM}) = \overline{\mathbb{E}}_{{\overline{\xi^{*}}}}\li[\li(\prod_{j \in J}\prod_{n = 1}^{m_{j}}\int^{t_{j}^{n}}_{r_{j}^{n}} F_{j}^{n}((\overline{X}_{-s})^{\mbox{\tiny $J$}}_{j})ds\re) \mathbb{E}_{\mu}\li[H^{*}(X_{0}, \overline{X}_{T})\re]e^{\int_{0}^{T}V(\overline{X}_{s})ds} \re]
\end{equation} 
for  all $-T \le r_{j}^{1} < t_{j}^{1} < \dots < r_{j}^{m_{j}} < t_{j}^{m_{j}} \le 0$ and all $F_{j}^{n} \in C_{b}(K \times I)$.
\end{Theorem}
%%%%%%%%%%%%%%%%%%%%%%%%%%%%%%%%%%%%%%%%%%%%%%%%%%%%%%%%%%%%%%%%%%%%%%%%%%%%%%%%%%%%%%%%%%%%%%%%%%%      Results, Duality, Remarks

\begin{Remark}\label{RemarkExchangeableLines}\point In analogy to Remark \ref{RemarkExchangeableTypes} the representation in Theorem \ref{TheoremStochasticRepresentation} and the behaviour of the BP imply that the ancestral lines and the genealogical distances of $J$ at time $T$ are exchangeable if $\mu^{*}$ is exchangeable. 
\end{Remark}
%%%%%%%%%%%%%%%%%%%%%%%%%%%%%%%%%%%%%%%%%%%%%%%%%%%%%%%%%%%%%%%%%%%%%%%%%%%%%%%%%%%%%%%%%%%%%%%%%%%      Formulation of the main results, The transformed BP 

\subsubsection{The transformed BP}\label{subsubTransformeddualprocess}
\noindent Here we define a transformation of the law of the path of the BP $\overline{X}$ from (\ref{equationOverlineX}) by changing the collection of probability measures
\begin{equation}
\li\{\overline{\mathbb{P}}_{\overline{\eta}} : \overline{\eta} \in \overline{\mathcal{E}}\re\}\;.
\end{equation}
This transformation shall be used to describe the conditional distribution of the extended ancestral lines of $J$ considered from {\sc Time} $0$ back to {\sc Time} $-T$ {\it given} the type information of $J$ at {\sc Time} $0$. 
%%%%%%%%%%%%%%%%%%%%%%%%%%%%%%%%%%%%%%%%%%%%%%%%%%%%%%%%%%%%%%%%%%%%%%%%%%%%%%%%%%%%%%%%%%%%%%%%%%%      Formulation of the main results, The transformed dual 

Observe that for fixed $\xi^{*} \in K^{J}$,
\begin{eqnarray}
&&\frac{\mathbb{E}_{\mu}\li[\li(\prod_{j \in J}\prod_{n = 1}^{m_{j}}\int_{r_{j}^{n}}^{t_{j}^{n}}F_{j}^{n}((X_{T})^{\mbox{\tiny $\mathcal{D}$}}_{j,s})ds\re) \mathbbm{1}\{((X_{T})^{*}_{j})_{j \in J}  = \xi^{*}\}\re]}{\mathbb{P}_{\mu}(((X_{T})^{*}_{j})_{j \in J} = \xi^{*})} \\
&=& \frac{\overline{\mathbb{E}}_{\overline{\xi^{*}}}\li[\li(\prod_{j \in J}\prod_{n = 1}^{m_{j}}\int^{t_{j}^{n}}_{r_{j}^{n}} F_{j}^{n}((\overline{X}_{-s})^{\mbox{\tiny $J$}}_{j})ds\re) \mathbb{E}_{\mu}[H^{*}(X_{0}, \overline{X}_{T})]e^{\int_{0}^{T}V(\overline{X}_{s})ds} \re]}{\overline{\mathbb{E}}_{\overline{\xi^{*}}}\li[\mathbb{E}_{\mu}[H^{*} (X_{0}, \overline{X}_{T})]e^{\int_{0}^{T}V(\overline{X}_{s})ds}\re]}
\end{eqnarray}
for all $-T \le r_{j}^{1} < t_{j}^{1} < \dots < r_{j}^{m_{j}} < t_{j}^{m_{j}} \le 0,$ and all $F_{j}^{n} \in C_{b}(K \times I)$ due to Theorem \ref{TheoremStochasticRepresentation}. This means we will {\it re-weight} the sample paths of the BP  by the functional
\begin{equation}
\li(\mathbb{E}_{\mu}\li[H^{*}(X_{0}, \overline{X}_{T-t})\re]e^{\int_{0}^{T-t}V(\overline{X}_{s})ds}\re)_{t \in [0,T]}
\end{equation}
to obtain a new {\it Markov process} that describes the reversed extended ancestral lines. 
%%%%%%%%%%%%%%%%%%%%%%%%%%%%%%%%%%%%%%%%%%%%%%%%%%%%%%%%%%%%%%%%%%%%%%%%%%%%%%%%%%%%%%%%%%%%%%%%%%      Formulation of the main results, The transformed BP 

The crucial point is that the new Markov process has the {\it same transitions} as the BP, but now appearing at {\it different rates}, namely changed by the time-space potential $h^{T}$ that is defined as follows:
\begin{Definition}[Time-space potential]\pointpar 
\noindent For $T> 0$ let
\begin{equation}\label{equationHT}
h^{T}(t,\overline{\eta}) := \overline{\mathbb{E}}_{\overline{\eta}}\li[\mathbb{E}_{\mu}\li[H^{*}(X_{0}, \overline{X}_{T-t})\re]e^{\int_{0}^{T-t}V(\overline{X}_{s})ds} \re]
\end{equation}
for all $t \in [0,T)$ and all $\overline{\eta} \in \overline{\mathcal{E}}$.
\end{Definition}
Note that
\begin{equation}\label{equationHTFormula}
h^{T}(t,\overline{\eta}) = \mathbb{E}_{\mu}\li[H^{*}(X_{T-t},\overline{\eta})\re] 
\end{equation}
due to Proposition \ref{PropositionDualityTypInfoHMM}, that is, at time $t$ the time-space potential can  be represented in terms of the type distribution of the HMM at time $T-t$. So, for a general initial type distribution $\mu^{*}$ the potential $h^{T}$ in fact depends on $t$ and therefore the transformed BP is {\it time-inhomogeneous}. However
\begin{equation}
h^{T}(t,\cdot) =h(\cdot) \; \mbox{ for all }\; 0 \le  t \le T \;\; \mbox{ if  $\mu^{*}$ is stationary }
\end{equation}
which yields a {\it time-homogeneous} transformed BP.

In order to change the probability measures formally, for $t\in [0,T]$ let $\pi_{t}: \overline{\Omega} \to \overline{\Omega}$, $(\overline{\omega}_{s})_{s \ge 0} \mapsto (\overline{\omega}_{s+t})_{s \ge 0}$ be the shift operator which cuts off the path $\overline{\omega}$ before time $t$ and shifts the remaining part in time. In addition, for each $t \in [0,T)$ and each $\overline{\eta} \in \overline{\mathcal{E}}$ consider the function
\begin{equation}
\overline{\mathbb{P}}_{t,\overline{\eta}}^{h^{T}}:\overline{\mathcal{F}}^{[t,T]}\to [0,1], \; C \mapsto \frac{\overline{\mathbb{E}}_{\overline{\eta}}\li[\mathbbm{1}_{\pi_{t}(C)} \mathbb{E}_{\mu}\li[H^{*}(X_{0}, \overline{X}_{T-t})\re]e^{\int_{0}^{T-t}V(\overline{X}_{s})ds}\re]}{h^{T}(t,\overline{\eta})}
\end{equation}
which is a probability measure on $(\overline{\Omega}, \overline{\mathcal{F}}^{[t,T]})$, where $\overline{\mathcal{F}}^{[t,T]} = \sigma(\overline{X}_{s}: s \in [t,T])$.  For technical reasons (recall (\ref{equationHTFormula}) and the fact that by definition in (\ref{equationOverlineMathcalE}) there are only non-empty subsets of $K$ in the second component of the BP) we also assume that 
\begin{equation}\label{equationEmuXtge0}
h^{T}:[0,T) \times \overline{\mathcal{E}} \to  (0,\infty)
\end{equation}
which typically holds, for example, if $b(u,v)$ is irreducible or if the initial type distribution $\mu^{*}$ puts positive mass on each element in $K^{J}$. 
%%%%%%%%%%%%%%%%%%%%%%%%%%%%%%%%%%%%%%%%%%%%%%%%%%%%%%%%%%%%%%%%%%%%%%%%%%%%%%%%%%%%%%%%%%%%%%%%%%%      Formulation of the main results, The transformed BP

\begin{Theorem}[Analytical characterization of the transformed BP]\point\label{TheoremTransformedBP}
\begin{enumerate}[a)]
\item\label{TheoremTransformedBPitemA} {\bf(Time-inhomogeneous)} We obtain that
\begin{equation}
\overline{X}^{h^{T}} =  \li\{\li(\overline{\Omega},\overline{\mathcal{F}}^{[t,T]}, (\overline{\mathcal{F}}^{[t,T]}_{s})_{s \in [t,T)},  (\overline{X}_{s})_{s \in [t,T]}, \{\overline{\mathbb{P}}_{t,\overline{\eta}}^{h^{T}}: \overline{\eta} \in \overline{\mathcal{E}}\} \re): t \in [0,T)\re\}
\end{equation} 
is a time-inhomogeneous Borel strong Markov process that corresponds to a family 
\begin{equation}\label{equationFamilyLhT}
\li\{\overline{L}^{h^{T}}_{t}: t \in [0,T)\re\}
\end{equation}
of bounded linear operators, where
\begin{equation}\label{equationLhT}
\overline{L}^{h^{T}}_{t}f(\overline{\eta}) = \sum_{\overline{\zeta} \in \overline{\mathcal{E}}} \overline{\mathcal{K}}(\overline{\eta}, \overline{\zeta})\li( \frac{h^{T}(t,\overline{\zeta})}{h^{T}(t,\overline{\eta})}\re)\li[f(\overline{\zeta}) - f(\overline{\eta})\re]\;.
\end{equation}
We call $\overline{X}^{h^{T}}$ the time-inhomogeneous transformed BP.
\item\label{TheoremTransformedBPitemB} {\bf(Time-homogeneous)} If $h:\overline{\mathcal{E}} \to  (0,\infty)$ and $h^{T}(t,\cdot) =h(\cdot)$ for all $0 \le  t < T < \infty$, then
\begin{equation}
\{\overline{\mathbb{P}}_{0,\overline{\eta}}^{h^{T}}: T \in [0,\infty)\} \; \mbox{ is projective for each }\; \overline{\eta} \in \overline{\mathcal{E}}\;, 
\end{equation}
which uniquely defines a family of probability measures
\begin{equation}
\{\overline{\mathbb{P}}_{\overline{\eta}}^{h}: \overline{\eta} \in \overline{\mathcal{E}}\}  \;\mbox{ on }\; (\overline{\Omega},\overline{\mathcal{F}})\;,
\end{equation} 
and
\begin{equation}
\overline{X}^{h} =  \li(\overline{\Omega},\overline{\mathcal{F}}, (\overline{\mathcal{F}}_{t})_{t \ge 0}, (\overline{X}_{t})_{t \ge 0}, \{\overline{\mathbb{P}}_{\overline{\eta}}^{h}: \overline{\eta} \in \overline{\mathcal{E}}\} \re)\;.
\end{equation} 
is a  Borel strong Markov process with bounded generator
\begin{equation}\label{equationLh}
\overline{L}^{h}f(\overline{\eta}) =  \sum_{\overline{\zeta} \in \overline{\mathcal{E}}} \overline{\mathcal{K}}(\overline{\eta}, \overline{\zeta})\li( \frac{h(\overline{\zeta})}{h(\overline{\eta})}\re)\li[f(\overline{\zeta}) - f(\overline{\eta})\re]\;, 
\end{equation}
namely, the compensated $h$-transform as introduced in \cite{FS04}. We refer to $\overline{X}^{h}$ as the time-homogeneous transformed BP.
\end{enumerate}
\end{Theorem}
%%%%%%%%%%%%%%%%%%%%%%%%%%%%%%%%%%%%%%%%%%%%%%%%%%%%%%%%%%%%%%%%%%%%%%%%%%%%%%%%%%%%%%%%%%%%%%%%%%%      Formulation of the main results, strong stochastic representation

\subsubsection{Strong stochastic representation for conditioned extended ancestral lines}\label{subsubStrongduality}
\noindent In this subsubsection we shall see that given the types of $J$ at {\sc Time} $0$ the conditional distribution of the extended ancestral lines of $J$ considered from {\sc Time} $0$ back to {\sc Time} $-T$ is equal to the distribution of a {\it special functional} of the sample paths of the transformed BP up to time $T$. 

Let
\begin{equation}
 \li((X_{T})^{\mbox{\tiny $\mathcal{D}$}}_{j,t}\re)_{t \in [-T,0],\, j \in J} \in (\mathcal{D}([-T,0], K \times I))^{J}
\end{equation}
be the extended ancestral lines of $J$ considered from {\sc Time} $0$ back to {\sc Time} $-T$ and
\begin{equation}
\li((X_{T})^{*}_{j}\re)_{j \in J}  \in  K^{J}
\end{equation}
be the types of $J$ at {\sc Time} $0$. More precisely, we have that in distribution
\begin{equation}
 (X_{T})^{\mbox{\tiny $\mathcal{D}$}}_{j,-t} =  (\overline{X}_{t})^{\mbox{\tiny $J$}}_{j} \;\; \mbox{ for all } t \in [0,T] \mbox{ and all } j \in J\;, 
\end{equation}
where the BP is considered under the new measure $ \overline{\mathbb{P}}_{0,\overline{\xi^{*}}}^{h^{T}}$. 

In order to state this relation between the extended ancestral lines of $J$ and the sample paths of the first component of the BP formally, for each $j \in J$ we first reverse the (right-continuous) path $((\overline{X}_{t})^{\mbox{\tiny $J$}}_{j} )_{t \in [0,T]}$, then  shift it to $[-T,0]$ and finally transform the resulting left-continuous path again into a right-continuous path. For this we use the special  map  
\begin{eqnarray}
\mathbb{F}_{T}: \mathcal{D}([0,T], (K \times I)^{J}) &\to& \li(\mathcal{D}([-T,0], K \times I)\re)^{J}\\ ((\overline{\omega}_{t})^{\mbox{\tiny $J$}})_{t \in [0,T]} &\mapsto& (\eta^{\mbox{\tiny $\mathcal{D}$}}_{j,-t})_{t \in [-T,0], j \in J}
\end{eqnarray}
that is defined by
\begin{equation}\label{equationFT}
\eta^{\mbox{\tiny $\mathcal{D}$}}_{j,-t} :=  \li\{ \begin{array}{ccl} (\overline{\omega}_{t})^{\mbox{\tiny $J$}}_{j} &,& t \in\{s \in [0,T]: (\overline{\omega}_{s})^{\mbox{\tiny $J$}}_{j} = (\overline{\omega}_{s-})^{\mbox{\tiny $J$}}_{j} \} \\ (\overline{\omega}_{t-})^{\mbox{\tiny $J$}}_{j} &,& t \in \{s \in [0,T]: (\overline{\omega}_{s})^{\mbox{\tiny $J$}}_{j} \not= (\overline{\omega}_{s-})^{\mbox{\tiny $J$}}_{j} \} \end{array} \re. \;,
\end{equation}
where $(\overline{\omega}_{t-})^{\mbox{\tiny $J$}}_{j}  := \lim_{\epsilon \downarrow 0}(\overline{\omega}_{t-\epsilon})^{\mbox{\tiny $J$}}_{j} $.
%%%%%%%%%%%%%%%%%%%%%%%%%%%%%%%%%%%%%%%%%%%%%%%%%%%%%%%%%%%%%%%%%%%%%%%%%%%%%%%%%%%%%%%%%%%%%%%%%%%      Results, Strong stochastic representation

Now we can state the strong stochastic representation for the conditioned extended ancestral lines in terms of the sample paths of the transformed BP, where it is useful to recall Theorem \ref{TheoremTransformedBP}.
\begin{Theorem}[Strong stochastic representation at time $T$]\label{TheoremStrongStochasticRepresentation}\pointpar
\noindent Let $J \subset I$, $T > 0$ and $\mu$ as in (\ref{equationMu}) with $c = -T$ (initial {\sc Time}) and  initial type distribution $\mu^{*}$. Let $\xi^{*} \in K^{J}$ and $\overline{\xi^{*}}$ be the initial state defined in (\ref{equationoverlinexi}) and assume that (\ref{equationEmuXtge0}) holds.

We have that
\begin{equation}
\mathbb{P}_{\mu}\li.\li(\li((X_{T})^{\mbox{\tiny $\mathcal{D}$}}_{j,t}\re)_{t \in [-T,0], j \in J} \in \, \cdot \,\re| \li((X_{T})^{*}_{j}\re)_{j \in J}= \xi^{*}\re) 
= \overline{\mathbb{P}}_{0,\overline{\xi^{*}}}^{h^{T}}\li( \mathbb{F}_{T}(\li((\overline{X}_{t})^{\mbox{\tiny $J$}}\re)_{t \in [0,T]}) \in \, \cdot \,\re)
\end{equation} 
if $\mu^{*}$ is general and
\begin{equation}
\mathbb{P}_{\mu}\li.\li(\li((X_{T})^{\mbox{\tiny $\mathcal{D}$}}_{j,t}\re)_{t \in [-T,0], j \in J} \in \, \cdot \,\re| \li((X_{T})^{*}_{j}\re)_{j \in J}= \xi^{*}\re)  = \overline{\mathbb{P}}_{\overline{\xi^{*}}}^{h}\li( \mathbb{F}_{T}(\li((\overline{X}_{t})^{\mbox{\tiny $J$}}\re)_{t \in [0,T]}) \in \, \cdot \,\re) 
\end{equation}
if $\mu^{*}$ is stationary.
\end{Theorem}
%%%%%%%%%%%%%%%%%%%%%%%%%%%%%%%%%%%%%%%%%%%%%%%%%%%%%%%%%%%%%%%%%%%%%%%%%%%%%%%%%%%%%%%%%%%%%%%%%%%      Formulation of the main results, strong stochastic representation, first application

\begin{Remark}\label{RemarkExchangeableConditionedLines}\point In analogy to  Remark \ref{RemarkExchangeableLines} we get that the conditioned ancestral lines and the conditioned genealogical distances of $J$ at time $T$ are exchangeable if $\mu^{*}$ is exchangeable. 
\end{Remark}

\begin{Remark}\label{RemarkExchangeableConditionedLines2}\point If we want to study  conditioned ancestral lines and  conditioned genealogical distances in the exchangeable situation of Remark \ref{RemarkExchangeableConditionedLines},  then we can reduce the transformed BP as described in item \ref{itemAnnotation5} on page \pageref{itemAnnotation5}.  This means that we only need the partition elements with their marks in $K$ (first component) and for each strict and non-empty subset of $K$ the number of active life sites to which this subset is assigned (second component).
\end{Remark}
%%%%%%%%%%%%%%%%%%%%%%%%%%%%%%%%%%%%%%%%%%%%%%%%%%%%%%%%%%%%%%%%%%%%%%%%%%%%%%%%%%%%%%%%%%%%%%%%%%%      Formulation of the main results, strong stochastic representation

In order to illustrate the strong stochastic representation we give a simple example in which we contrast the genealogical distance of two individuals at time $T$ with the conditioned genealogical distance.

\begin{Example}[Genealogical distances]\label{Example}\pointpar
\noindent Recall (\ref{equationDt}) and consider the  case where $B = 0$ (no mutation), $S =0$ (no selection), $\mu^{*} = \nu^{\otimes N}$ with $\nu \in \mathcal{M}_{1}(K)$ (exchangeable initial type distribution), $J = \{i,j\}$  and $u,v \in K$.

In this simple situation the  probability that the genealogical distance of $i$ and $j$ at time $T$ is greater than $2t$ is well known (e.g$.$ see Corollary 3.4 in \cite{GPW13} for a rigorous argument). Namely for each $t \ge 0$,
\begin{equation}\label{equationExpMintT}
\mathbb{P}_{\mu}\li(\mathbb{D}_{T}(i,j) > 2t\re) = e^{-t}\mathbbm{1}\{t < T\}\;,
\end{equation}
where $\mathbb{P}_{\mu}(\mathbb{D}_{T}(i,j) = 2T) = e^{-T}$ describes the probability that $i$ and $j$ have no common ancestor until time $T$.
%%%%%%%%%%%%%%%%%%%%%%%%%%%%%%%%%%%%%%%%%%%%%%%%%%%%%%%%%%%%%%%%%%%%%%%%%%%%%%%%%%%%%%%%%%%%%%%%%%%      Formulation of the main results, strong stochastic representation

Now we consider the conditional probability of this event, where it is helpful to recall item \ref{itemAnnotation2} on page \pageref{itemAnnotation2}. The strong stochastic representation yields (set $J = \{i,j\}$ and $\xi^{*} = (u,v)$) that for each $t \ge 0$,
\begin{equation}
\mathbb{P}_{\mu}(\mathbb{D}_{T}(i,j) > 2t \,| (X_{T})^{*}_{i} = u, (X_{T})^{*}_{j} = v) = \overline{\mathbb{P}}_{0, \overline{\xi^{*}}}^{h^{T}}(\inf\{r \in [0,T]: (\overline{X}_{r})^{\mbox{\tiny $J$}}_{i} = (\overline{X}_{r})^{\mbox{\tiny $J$}}_{j} \} > t)
\end{equation}
where we use the convention $\inf\{\emptyset\} = T$. Since $(\overline{X}_{r})^{\mbox{\tiny $J$}}_{i} = (\overline{X}_{r})^{\mbox{\tiny $J$}}_{j}$ if and only if $ \{i,j\} \in \Gamma(\overline{X}_{r})$, we have to determine the first time the partition elements $\{i\}$ and $\{j\}$ (in the first component of the time-inhomogeneous transformed BP) coalesce. So, if $u \not= v$, then the partition elements cannot coalesce. If $u = v$, then at time $r \in [0,T]$ the partition elements $\{i\}$ and $\{j\}$ coalesce at rate
\begin{equation}
\frac{\mathbb{P}_{\mu}((X_{T-r})^{*}_{i} = u)}{\mathbb{P}_{\mu}((X_{T-r})^{*}_{i} = (X_{T-r})^{*}_{j} = u)} = \frac{\nu(\{u\})}{\nu(\{u\})(1-e^{r-T}) + \nu(\{u\})^{2}e^{r-T}} \;.%=  \frac{1}{1 - e^{r-T}\nu(K \setminus \{u\})} \;.
\end{equation}
Hence
\begin{equation}\label{equationExpMintTConditioned}
\mathbb{P}_{\mu}(\mathbb{D}_{T}(i,j) > 2t \,| (X_{T})^{*}_{i} = u, (X_{T})^{*}_{j} = v) = \li\{\begin{array}{ccc} \frac{e^{-t} - e^{-T}\nu(K \setminus \{u\})}{1 - e^{-T}\nu(K \setminus \{u\})} \mathbbm{1}\{t < T\} &,& u = v \\ \mathbbm{1}\{t < T\} &,&  u \not= v \end{array}\re. .
\end{equation}
Note that by calculations the r.h.s$.$ of (\ref{equationExpMintT}) can be recovered by the r.h.s$.$ of (\ref{equationExpMintTConditioned}), but a formal argument is missing.
\end{Example}

Finally, the formula in (\ref{equationExpMintTConditioned}) implies the following
\begin{Corollary}[Conditioned exponential random variable]\pointpar 
\noindent Assume that $S=0$, $B = 0$ and $\mu^{*} = \nu^{\otimes N}$. If type $u$ is rare in the initial population (more precisely, if $\nu(\{u\}) \to 0$), then the  conditional distribution of half the genealogical distance of $i$ and $j$ at time $T$ given that both  $i$ and $j$ have type $u$ at time $T$ is equal to the law of  an exponential random variable which is conditioned to take values in $[0,T]$.
\end{Corollary}
%%%%%%%%%%%%%%%%%%%%%%%%%%%%%%%%%%%%%%%%%%%%%%%%%%%%%%%%%%%%%%%%%%%%%%%%%%%%%%%%%%%%%%%%%%%%%%%%%%%      Formulation of the main results, Ergodicity and applications

\subsection{Longtime behaviour of extended ancestral lines}\label{subsecErgodicity}
\noindent The main application of the strong stochastic representation concerns the longtime behaviour of the HMM. Remember that $\overline{X}^{h}$ is the time-homogeneous transformed BP given in Theorem \ref{TheoremTransformedBP} \ref{TheoremTransformedBPitemB}) and that $\mu$ stands for the initial distribution of the HMM defined in (\ref{equationMu}) with $c = -T$ (initial {\sc Time}) and  initial type distribution $\mu^{*}$. Furthermore, let 
\begin{equation}
\mathbb{F}: \mathcal{D}([0,\infty), (K \times I)^{J})  \to (\mathcal{D}((-\infty,0], K \times I))^{J}
\end{equation}
be defined in analogy to $\mathbb{F}_{T}$ in (\ref{equationFT}).
%%%%%%%%%%%%%%%%%%%%%%%%%%%%%%%%%%%%%%%%%%%%%%%%%%%%%%%%%%%%%%%%%%%%%%%%%%%%%%%%%%%%%%%%%%%%%%%%%%%      Formulation ..., strong stochastic representation, long time

\begin{Theorem}[Extended ancestral lines in the limit $T \to \infty$]\label{TheoremLongtime}\pointpar
\noindent Assume that $B >0$ (with mutation) and $b(\cdot,\cdot)$ is irreducible, i.e$.$ 
\begin{equation}
\lim_{T \to \infty}\mathbb{P}_{\mu}((X_{T})^{*} \in \, \cdot \,) =: \mathbb{P}((X_{\infty})^{*} \in \, \cdot \,) \;\in\;  \mathcal{M}_{1}(K^{I})
\end{equation}
for each initial type distribution $\mu^{*}$, where the limit is the unique stationary type distribution that puts positive mass on each element in $K^{I}$, i.e$.$
\begin{equation}
h(\overline{\eta}) := \mathbb{E}\li[H^{*}(X_{\infty},\overline{\eta})\re] > 0 \;\; \mbox{ for all }\; \overline{\eta}  \in  \overline{\mathcal{E}} \;.
\end{equation}
Then the following holds:
\begin{enumerate}[a)]
\item\label{TheoremLongtimeItemA} For each initial type distribution $\mu^{*}$,
\begin{equation}
\lim_{T \to \infty}\mathbb{P}_{\mu}\li((X_{T})^{\mbox{\tiny $\mathcal{D}$}}\in \, \cdot \,\re) =: \mathbb{P} \li((X_{\infty})^{\mbox{\tiny $\mathcal{D}$}} \in \, \cdot \,\re) \;\in\;  \mathcal{M}_{1}((\mathcal{D}(\mathbb{R}, K \times I))^{I})\;.
\end{equation}
\item\label{TheoremLongtimeItemB} For $J \subset I$ and $\xi^{*} \in K^{J}$,
\begin{equation}
\mathbb{P}\li.\li(\li((X_{\infty})^{\mbox{\tiny $\mathcal{D}$}}_{j,t}\re)_{t \in (-\infty,0], j \in J} \in \, \cdot \,\re| \li((X_{\infty})^{*}_{j}\re)_{j \in J}= \xi^{*}\re) = \overline{\mathbb{P}}_{\overline{\xi^{*}}}^{h}\li( \mathbb{F}(\li((\overline{X}_{t})^{\mbox{\tiny $J$}}\re)_{t \ge 0}) \in \, \cdot \,\re)\;.
\end{equation}
\end{enumerate}
\end{Theorem}
%%%%%%%%%%%%%%%%%%%%%%%%%%%%%%%%%%%%%%%%%%%%%%%%%%%%%%%%%%%%%%%%%%%%%%%%%%%%%%%%%%%%%%%%%%%%%%%%%%%      Formulation ..., strong stochastic representation, long time, Remark

\begin{Remark}[Fixation]\label{RemarkTheoremLongtime}\point If $B=0$ (no mutation), then one can show that
\begin{equation}
\lim_{T \to \infty}\mathbb{P}_{\mu}\li((X_{T})^{\mbox{\tiny $\mathcal{D}$}}\in \, \cdot \,\re) = \sum_{u \in K} \mathbb{P}\li((X_{\infty})^{*}_{i} = u \mbox{ for all }i \in I\re) \mathbb{P}_{u} \li((X_{\infty})^{\mbox{\tiny $\mathcal{D}$}} \in \, \cdot \,\re) \;,
\end{equation}
where $\mathbb{P}\li((X_{\infty})^{*}_{i} = u \mbox{ for all }i \in I\re)$ is the probability (depending on the initial type distribution $\mu^{*}$) that eventually all individual have type $u$ and $\mathbb{P}_{u}$ can be described as follows. Under $\mathbb{P}_{u}$ the ancestral lines are constant paths through $u$,  i.e$.$
\begin{equation}
\mathbb{P}_{u}\li(\li((X_{\infty})^{\mbox{\tiny $\mathcal{D}$}}_{i,t}\re)_{t \in \mathbb{R}, i \in I} \in \, (\mathcal{D}(\mathbb{R}, \{u\} \times I))^{I}\,\re) = 1\;,
\end{equation}
and
\begin{equation}
\mathbb{P}_{u}\li(\li((X_{\infty})^{\mbox{\tiny $\mathcal{D}$}}_{i,t,I}\re)_{t \in (-\infty,0], i \in I} \in \, \cdot \,\re) \;\in\;  \mathcal{M}_{1}((\mathcal{D}((-\infty,0], I))^{I})
\end{equation}
arises as the law of a system of $N$ instantaneously coalescing random walks on $I$.
\end{Remark}

\begin{Remark}\label{RemarkExchangeableConditionedLinesInfty}\point The unique stationary type distribution is exchangeable due to the mutation and resampling dynamics. Hence  under $\mathbb{P}$ the conditioned ancestral lines and the conditioned genealogical distances are exchangeable and for calculations we can reduce the time-homogeneous transformed BP in analogy to Remark \ref{RemarkExchangeableConditionedLines2}
\end{Remark}
%%%%%%%%%%%%%%%%%%%%%%%%%%%%%%%%%%%%%%%%%%%%%%%%%%%%%%%%%%%%%%%%%%%%%%%%%%%%%%%%%%%%%%%%%%%%%%%%%%%      Applications

\subsection{Applications}\label{subsecApplications}
In this subsection we consider the CAT and the genealogical distance of two individuals in equilibrium, two important issues we can study rigorously within the framework of the HMM.  We first work out particular questions and open problems. Then we investigate these questions concretely in the two type case ($|K| = 2$) in the Subsubsections \ref{subsubCAT} and \ref{subsubConditionedGenealogicalDistance}, where we also let the population size $N$ tend to infinity.

Assume that we are in the situation of Theorem \ref{TheoremLongtime}. Hence, see Remark \ref{RemarkExchangeableConditionedLinesInfty}, the ancestral lines and the genealogical distances are exchangeable and we can use the reduced time-homogeneous transformed BP (i.e$.$ in the first component we only need the  partition elements with their marks in $K$  and  in the second component only how often  each strict and non-empty subset of $K$ occurs) in the sequel. \bigskip
%%%%%%%%%%%%%%%%%%%%%%%%%%%%%%%%%%%%%%%%%%%%%%%%%%%%%%%%%%%%%%%%%%%%%%%%%%%%%%%%%%%%%%%%%%%%%%%%%%%      Applications, CAT, intro

{\bf The stationary type distribution of the CAT, a reformulation}: First of all note that the analytical characterization given in Theorem \ref{TheoremTransformedBP} \ref{TheoremTransformedBPitemB})  implies that eventually all partition elements in the first component of the time-homogeneous transformed BP merge together almost surely. Hence Theorem \ref{TheoremLongtime} \ref{TheoremLongtimeItemB}) (set $J = I$) implies  that all extended ancestral lines coincide at some time in the past almost surely. Formally,
\begin{equation}\label{equationExtendAncestralLinesCoincide}
\mathbb{P}\li(\exists t \le 0 \mbox{ such that }  (X_{\infty})^{\mbox{\tiny $\mathcal{D}$}}_{i,s} = (X_{\infty})^{\mbox{\tiny $\mathcal{D}$}}_{j,s} \mbox{ for all } s \le t \mbox{ and all } j,i \in I\re) = 1\;.
\end{equation}
Thus we can trace a single ancestral line back to {\sc Time} $-\infty$ to obtain the stationary type distribution of the CAT.

Concretely we fix $j \in J$ and consider the limit
\begin{equation}
\lim_{t \to \infty} \mathbb{P}\li( (X_{\infty})^{\mbox{\tiny $\mathcal{D}$}}_{j,-t,K}  \in  \cdot \re) = \lim_{t \to \infty} \sum_{u \in K}\mathbb{P}\li( \li. (X_{\infty})^{\mbox{\tiny $\mathcal{D}$}}_{j,-t,K}  \in  \cdot \,\re|(X_{\infty})^{*}_{j}=u\re) \mathbb{P}( (X_{\infty})^{*}_{j}=u)\;.
\end{equation}
Then again Theorem \ref{TheoremLongtime} \ref{TheoremLongtimeItemB})  implies (set $J= \{j\}$ and $\xi^{*} = u$) that this limit exists and is given (in terms of the first component of the reduced time-homogeneous transformed BP) by 
\begin{equation}\label{equationLimitDistribution}
\lim_{t \to \infty}  \sum_{u \in K}\overline{\mathbb{P}}_{\overline{\xi^{*}}}^{h}\li( \li(\overline{X}_{t}\re)^{\mbox{\tiny $J$}}_{K} \in \, \cdot \,\re)\mathbb{P}((X_{\infty})^{*}_{j}=u) \;.
\end{equation}
In other words, the equilibrium distribution of the first component of the reduced time-homogeneous transformed BP represents the stationary type distribution of the CAT for the  multi-type Moran model with selection and type-dependent mutation. However, the first component itself is {\it not} Markovian. So, one has to determine the equilibrium distribution of the whole (first and second component) reduced time-homogeneous transformed BP in order to get (\ref{equationLimitDistribution}). But this is a problem on its own whose complexity increases with the number of types in $K$, in particular, when one is interested in an explicit form of (\ref{equationLimitDistribution}).

In the two type case, see also (\ref{equationAnzahlNuller}), we therefore have to understand how the functional
\begin{equation}\label{equationReducedBPforCAT}
\li((\overline{X}_{t})^{\mbox{\tiny $J$}}_{K},(\overline{X}_{t})^{\mbox{\tiny $I$}}_{|\{0\}|}\re)_{t \ge 0} \;\in\; \mathcal{D}([0,\infty), K \times \{0, \dots , N-1\})
\end{equation}
evolves under $\overline{\mathbb{P}}_{\overline{\xi^{*}}}^{h}$.\bigskip
%%%%%%%%%%%%%%%%%%%%%%%%%%%%%%%%%%%%%%%%%%%%%%%%%%%%%%%%%%%%%%%%%%%%%%%%%%%%%%%%%%%%%%%%%%%%%%%%%%%      Applications, genealogica distances, intro

{\bf The genealogical distance of two individuals, a reformulation}: For $i,j \in I$ consider the functional
\begin{equation}
\mathbb{D}_{\infty}(i,j):= 2 \li|\sup\li\{ r \le 0: (X_{\infty})^{\mbox{\tiny $\mathcal{D}$}}_{i,r} = (X_{\infty})^{\mbox{\tiny $\mathcal{D}$}}_{j,r}\re\}\re|
\end{equation}
that describes  under $\mathbb{P}$ the genealogical distance of $i$ and $j$  in equilibrium. Observe that this random variable converges to the genealogical distance of two individuals sampled from the tree-valued Fleming-Viot process in equilibrium if $N \to \infty$. 

It is known, see for example \cite{DGP12}, that the random variable $\mathbb{D}_{\infty}(i,j)$ is exponential distributed (for each $N$ and hence in the limit $N \to \infty$) if there is no selection. In the case with selection the distribution is not known.  However, there is the conjecture (a proof is not yet available) that for each selection coefficient $S >0$ the random variable $\mathbb{D}_{\infty}(i,j)$ is stochastically smaller than for selection coefficient $S=0$, where one cannot expect that this ordering is monotone in $S$.

With our machinery we can approach this open problem. Namely, we express the conditional probability of the event that the genealogical distance of $i$ and $j$ in equilibrium (given the types of $i$ and $j$) is greater than $2t$ in terms of the time-homogeneous transformed BP and thus we can study the  tail distribution function of the conditioned genealogical distance and therefore of the genealogical distance. Formally, the strong stochastic representation in Theorem \ref{TheoremLongtime} \ref{TheoremLongtimeItemB}) implies (set $J = \{i,j\}$ and $\xi^{*} = (u,v)$) that
\begin{equation}
\mathbb{P}(\mathbb{D}_{\infty}(i,j) > 2t \,| (X_{\infty})^{*}_{i} = u, (X_{\infty})^{*}_{j} = v) = \overline{\mathbb{P}}_{\overline{\xi^{*}}}^{h}(\inf\{r \ge 0: (\overline{X}_{r})^{\mbox{\tiny $J$}}_{i} = (\overline{X}_{r})^{\mbox{\tiny $J$}}_{j} \} > t)
\end{equation}
for all $t \ge 0$. This means, see also the example at the end of Subsubsection \ref{subsubStrongduality}, we have to determine the first time  the partition elements $\{i\}$ and $\{j\}$ coalesce. Thus we have to consider the whole reduced time-homogeneous transformed BP which gives rise to the same problems we have discussed for the CAT.  

In the two type case (remember that $\Gamma(\overline{X}_{t}) \in \{\{\{i\},\{j\}\}, \{\{i,j\}\}\}$ denotes the partition contained in the first component of the BP) we therefore have to understand how under $\overline{\mathbb{P}}_{\overline{\xi^{*}}}^{h}$ the functional
\begin{equation}\label{equationReducedBPforDistances}
\li(\Gamma(\overline{X}_{t}), (\overline{X}_{t})^{\mbox{\tiny $J$}}_{i,K}, (\overline{X}_{t})^{\mbox{\tiny $J$}}_{j,K},(\overline{X}_{t})^{\mbox{\tiny $I$}}_{|\{0\}|}\re)_{t \ge 0} 
\end{equation}
evolves until the time the partition elements $\{i\}$ and $\{j\}$ coalesce.
\bigskip

Now we prepare the setting for the Subsubsections \ref{subsubCAT} and \ref{subsubConditionedGenealogicalDistance}, where we also let $N \to \infty$. Observe that $\mathbb{P}$ and $\overline{\mathbb{P}}_{\overline{\xi^{*}}}^{h}$ depend on the population size $N$ in the sequel, although we do not indicate this.

Assume that $K = \{0,1\}$ (i.e$.$ $\chi(0) = 0 < 1 = \chi(1)$ due to (\ref{equationDefChi})), $B > 0$, $b_{0} := b(0,0) = b(1,0)>0$ and $b_{1} := b(0,1) = b(1,1)>0$ (i.e$.$ the assumption of Theorem \ref{TheoremLongtime} is satisfied) and $S \ge 0$ (this allows us to distinguish between the case with selection and without selection). 

First note that the explicit form of the unique stationary distribution $\mathbb{P}((X_{\infty})^{*} \in \cdot)$ is not important here, we only use that
\begin{equation}
\lim_{N \to \infty}\mathbb{P}\circ\li(\frac{1}{N}\sum_{i = 1}^{N}(X_{\infty})^{*}_{i}\re)^{-1} = \pi\;,
\end{equation}
where $\pi \in \mathcal{M}_{1}([0,1])$ is the unique stationary distribution of a Wright-Fisher diffusion with generator
\begin{equation}
\hat{L}f(z) = \li[Bb_{1}(1-z) - Bb_{0}z +Sz(1-z)\re]f'(z)+\frac{1}{2}z(1-z)f''(z),\;\;z \in [0,1]\;,
\end{equation}
which acts on functions $f \in C^{2}([0,1])$.

In order to specify the transitions of (\ref{equationReducedBPforCAT}) and (\ref{equationReducedBPforDistances}) we need, recall (\ref{equationHTFormula}), for each $n,m\in \mathbb{N}_{0}$ with $n+m \le N$ the probability that in equilibrium $n$ of $N$ live-site have type $1$ and $m$ of $N$ have type $0$. So, we set
\begin{equation}\label{equationDefinitionPN}
P_{N}(1^{n},0^{m}) := \mathbb{P}((X_{\infty})^{*}_{1} = 1,\dots,(X_{\infty})^{*}_{n} = 1,(X_{\infty})^{*}_{n+1} = 0, \dots ,(X_{\infty})^{*}_{n+m} = 0)\;,
\end{equation}
where we abbreviate $P_{N}(1^{n},0^{0})$  by $P_{N}(1^{n})$, $P_{N}(1^{n},0^{1})$ by $P_{N}(1^{n},0)$, $P_{N}(1^{0},0^{m})$ by $P_{N}(0^{m})$  and $P_{N}(1^{1},0^{m})$ by $P_{N}(1,0^{m})$ in the sequel.  In the limit $N \to \infty$  these probabilities  converge to mixed moments of the Wright-Fisher diffusion in equilibrium of the form 
\begin{equation}\label{equationDefinitionMixedMoments}
E[1^{n},0^{m}] := \int_{[0,1]}z^{n}(1-z)^{m}\pi(dz) \;\;, \; n,m \in \mathbb{N}_{0}\;,
\end{equation}
where we use the analogue abbreviations $E[1^{n}]$, $E[1^{n},0]$, $E[0^{m}]$ and $E[1,0^{m}]$.
%%%%%%%%%%%%%%%%%%%%%%%%%%%%%%%%%%%%%%%%%%%%%%%%%%%%%%%%%%%%%%%%%%%%%%%%%%%%%%%%%%%%%%%%%%%%%%%%%%%      Applications, CAT

\subsubsection{Stationary type distribution of the CAT}\label{subsubCAT}
\noindent Assume that we are in the two type case $K = \{0,1\}$. First we define a version of (\ref{equationReducedBPforCAT}) in terms of a pure Markov jump process 
\begin{equation}
\overline{Y}^{N} = (\overline{Y}^{N}_{t})_{t \ge 0} \;\mbox{ on }\; K \times \{0, \dots , N-1\}
\end{equation}
that has a unique equilibrium distribution which describes the stationary type distribution of the CAT. Then we let $N \to \infty$ and obtain the common ancestor process of Fearnhead. Note that due to the HMM (modeling the forward-in-time evolution of the ancestral lines) our result is a mathematically rigorous proof that the common ancestor process of Fearnhead follows a single ancestral line back to time $-\infty$.   

The jump process $\overline{Y}^{N}$ is defined (on a suitable probability space) by the  transitions listed in the second column of Figure \ref{FigureTransitionsY}. It is not hard to see that $\overline{Y}^{N}$ is an irreducible Markov chain on a finite set and therefore has a unique equilibrium state $\overline{Y}^{N}_{\infty}$. 

Note, if $S = 0$ (no selection), then $\overline{Y}^{N}$ is a jump process on $K \times \{0\}$ with $(0,0) \to (1,0)$ at rate $Bb_{1}$ and $(1,0) \to (0,0)$  at rate $Bb_{0}$. Hence $P(\overline{Y}^{N}_{\infty} = (0,0)) = b_{0}$ and $P(\overline{Y}^{N}_{\infty} = (1,0)) = b_{1}$.
%%%%%%%%%%%%%%%%%%%%%%%%%%%%%%%%%%%%%%%%%%%%%%%%%%%%%%%%%%%%%%%%%%%%%%%%%%%%%%%%%%%%%%%%%%%%%%%%%%%      Applications, CAT, Figure

\renewcommand{\arraystretch}{1.5}
\begin{figure}[h]\caption{Transitions of $\overline{Y}^{N}$ and $\overline{Y}$}
\label{FigureTransitionsY}
\begin{center}
\begin{tabular}{|l|l|l|}\hline
transition &  rate for $\overline{Y}^{N}$& rate for $\overline{Y}$\\ \hline 
%-------------------------------n to n+1---------------------------------
$(u,n) \to (u,n+1)$  & $S(n+1)(\frac{N-1-n}{N})\frac{P_{N}(1^{u},0^{n+2-u})}{P_{N}(1^{u},0^{n+1-u})}$ & $S(n+1)\frac{E[1^{u},0^{n+2-u}]}{E[1^{u},0^{n+1-u}]}$\\ \hline
%-------------------------------n to n-1---------------------------------
$(u,n) \to (u,n-1)$ & $[Bb_{0}n + {n +1 -u \choose 2}(\frac{N-S}{N})]\frac{P_{N}(1^{u},0^{n-u})}{P_{N}(1^{u},0^{n+1-u})}$ & $[Bb_{0}n + {n +1 -u \choose 2}]\frac{E[1^{u},0^{n-u}]}{E[1^{u},0^{n+1-u}]}$\\ \hline
%-------------------------------u to 1-u---------------------------------
$(u,n) \to (1-u,n)$ & $B[ub_{1}+ (1-u)b_{0}]\frac{P_{N}(1^{1-u},0^{n+u})}{P_{N}(1^{u},0^{n+1-u})}$ & $B[ub_{1}+ (1-u)b_{0}]\frac{E[1^{1-u},0^{n+u}]}{E[1^{u},0^{n+1-u}]}$\\ \hline
\end{tabular} 
\end{center}
\end{figure}
%%%%%%%%%%%%%%%%%%%%%%%%%%%%%%%%%%%%%%%%%%%%%%%%%%%%%%%%%%%%%%%%%%%%%%%%%%%%%%%%%%%%%%%%%%%%%%%%%%%      Applications, CAT

\begin{Proposition}[Stationary type distribution of the CAT]\label{PropositionCAT}\pointpar
\noindent Let $K = \{0,1\}$, $J = \{j\}$ and $\xi^{*} = u \in K$. We have that
\begin{equation}
P((\overline{Y}^{N}_{t})_{t \ge 0} \in \, \cdot \,| \, \overline{Y}^{N}_{0} = (u,0)) =  \overline{\mathbb{P}}_{\overline{\xi^{*}}}^{h}(\li((\overline{X}_{t})^{\mbox{\tiny $J$}}_{K},(\overline{X}_{t})^{\mbox{\tiny $I$}}_{|\{0\}|}\re)_{t \ge 0} \in \, \cdot \,)\;.
\end{equation}
This means
\begin{equation}
\li(\sum_{n = 0}^{N-1} P(\overline{Y}^{N}_{\infty} = (0,n)) , \sum_{n = 0}^{N-1} P(\overline{Y}^{N}_{\infty} = (1,n)) \re)
\end{equation}
is the probability vector describing the stationary type distribution of the CAT.
\end{Proposition}
%%%%%%%%%%%%%%%%%%%%%%%%%%%%%%%%%%%%%%%%%%%%%%%%%%%%%%%%%%%%%%%%%%%%%%%%%%%%%%%%%%%%%%%%%%%%%%%%%%%      Applications, CAT

Now we let $N \to \infty$. In this case the sequence $\{\overline{Y}^{N}: N \in \mathbb{N}\}$ converges to a pure Markov jump process 
\begin{equation}
\overline{Y} = (\overline{Y}_{t})_{t \ge 0} \;\mbox{ on } \; K \times \mathbb{N}_{0} 
\end{equation}
with transitions as given in the third column of Figure \ref{FigureTransitionsY}. If we consider mutation rate $\tfrac{B}{2}$ and selection coefficient $\tfrac{S}{2}$, then (see Section 2 in \cite{F02}) the process
\begin{center}
{\bf $\overline{Y}$ coincides with the common ancestor process of Fearnhead!}
\end{center}
But note that in the definition of the common ancestor process on page 45 in \cite{F02} the factor $n+1$ in the transition rate of $(u,n) \to (u,n+1)$ is missing, but not in the proof as well as in subsequent statements.
  
\begin{Remark}\point  Since the time-homogeneous transformed BP is defined for general finite $K$, we could now work out a representation for  more than two types (which is of interest). However, for general $K$ each strict and non-empty subset of $K$ can occur in the second component of the time-homogeneous transformed BP. This means the stationary type distribution of the CAT is given by the unique equilibrium state of a pure Markov jump process on $K \times (\mathbb{N}_{0})^{\tilde{d}}$ (which has the same form as $\overline{Y}$), where $\tilde{d}$ is the number of strict and non-empty subsets of $K$ (e.g$.$ if $K = \{0,1,2\}$, then $\tilde{d} = 6$). So we do not carry out this program in the present paper and defer it to future work.
\end{Remark}
%%%%%%%%%%%%%%%%%%%%%%%%%%%%%%%%%%%%%%%%%%%%%%%%%%%%%%%%%%%%%%%%%%%%%%%%%%%%%%%%%%%%%%%%%%%%%%%%%%%      Applications, genealogical distance

\subsubsection{Conditioned genealogical distance of two individuals}\label{subsubConditionedGenealogicalDistance}
Assume that we are in the two type case $K = \{0,1\}$. Here we define a pure Markov jump process
\begin{equation}
\overline{Y}^{2,N} = (\overline{Y}^{2,N}_{t})_{t \ge 0} \; \mbox{  on }\;  \{\bigtriangleup \} \cup K \times K \times \{0, \dots , N-2\}
\end{equation}
that is a version of (\ref{equationReducedBPforDistances}) until it reaches its absorbing state $\bigtriangleup$. The jump process  $\overline{Y}^{2,N}$ allows us to express the conditioned genealogical distance of two individuals in equilibrium by means of the  hitting time of the absorbing state $\bigtriangleup$.  Then we let $N \to \infty$ in order to investigate the tail distribution function near zero of both the genealogical distance  and the conditioned genealogical distance of two individuals sampled from the tree-valued Fleming-Viot process in equilibrium,  where the case $S = 0$ (no selection) shall be compared with the case $S > 0$.
%%%%%%%%%%%%%%%%%%%%%%%%%%%%%%%%%%%%%%%%%%%%%%%%%%%%%%%%%%%%%%%%%%%%%%%%%%%%%%%%%%%%%%%%%%%%%%%%%%%      Applications, genealogical distance, Y^{2,N}

For our purposes we identify the state $(0,0) \in K\times K$ by $ \Circle$ and the state $(1,1) \in K\times K$  by  $\CIRCLE$ and we combine (remember that genealogical distances are  exchangeable) the states $(0,1)$ and $(1,0)$ to form a new element denoted  by $\RIGHTcircle$. With this we define $\overline{Y}^{2,N}$ as the jump process on
\begin{equation}
\{\bigtriangleup \} \cup \{\Circle, \CIRCLE, \RIGHTcircle\} \times \{0, \dots , N-2\}
\end{equation}
with {\it transitions} as given in the {\it second} column of Figure \ref{FigureTransitionsY2}.
\renewcommand{\arraystretch}{1.5}
\begin{figure}[h]\caption{Transitions of $\overline{Y}^{2,N}$ and $\overline{Y}^{2}$}
\label{FigureTransitionsY2}
\begin{center}
\begin{tabular}{|l|l|l|}\hline
transition &  rate for  $\overline{Y}^{2,N}$& rate for $\overline{Y}^{2}$\\ \hline 
%-------------------------------n to n+1---------------------------------
$(\Circle,n) \to (\Circle,n+1)$ & $S(n+2)(\frac{N-2-n}{N}) \frac{P_{N}(0^{n+3})}{P_{N}(0^{n+2})}$  & $S(n+2) \frac{E[0^{n+3}]}{E[0^{n+2}]}$\\ \hline
$(\CIRCLE,n) \to (\CIRCLE,n+1)$ & $S(n+2)(\frac{N-2-n}{N}) \frac{P_{N}(1^{2},0^{n+1})}{P_{N}(1^{2},0^{n})}$  & $S(n+2) \frac{E[1^{2},0^{n+1}]}{E[1^{2},0^{n}]}$\\ \hline
$(\RIGHTcircle,n) \to (\RIGHTcircle,n+1)$ & $S(n+2)(\frac{N-2-n}{N}) \frac{P_{N}(1,0^{n+2})}{P_{N}(1,0^{n+1})}$  & $S(n+2) \frac{E[1,0^{n+2}]}{E[1,0^{n+1}]}$\\ \hline
%-------------------------------n to n-1---------------------------------
$(\Circle,n) \to (\Circle,n-1)$ & $(Bb_{0}n + ({n +2 \choose 2}-1)\frac{N-S}{N})\frac{P_{N}(0^{n+1})}{P_{N}(0^{n+2})} $  & $(Bb_{0}n + {n +2 \choose 2} - 1)\frac{E[0^{n+1}]}{E[0^{n+2}]}$\\ \hline
$(\CIRCLE,n) \to (\CIRCLE,n-1)$ & $(Bb_{0}n + {n \choose 2}\frac{N-S}{N})\frac{P_{N}(1^{2},0^{n-1})}{P_{N}(1^{2},0^{n})} $  & $(Bb_{0}n + {n \choose 2})\frac{E[1^{2},0^{n-1}]}{E[1^{2},0^{n}]}$\\ \hline
$(\RIGHTcircle,n) \to (\RIGHTcircle,n-1)$ & $(Bb_{0}n + {n +1 \choose 2}\frac{N-S}{N})\frac{P_{N}(1,0^{n})}{P_{N}(1,0^{n+1})} $  & $(Bb_{0}n + {n +1 \choose 2})\frac{E[1,0^{n}]}{E[1,0^{n+1}]}$\\ \hline
%----------------------------- \Circle to \RIGHTcircle ------------------
$(\Circle,n) \to (\RIGHTcircle,n)$ & $2Bb_{0}\frac{P_{N}(1,0^{n+1})}{P_{N}(0^{n+2})}$ &  $2Bb_{0}\frac{E[1,0^{n+1}]}{E[0^{n+2}]}$\\ \hline
%-----------------------------\CIRCLE to \RIGHTcircle ------------------
$(\CIRCLE,n)\to (\RIGHTcircle,n)$ & $2Bb_{1}\frac{P_{N}(1,0^{n+1})}{P_{N}(1^{2},0^{n})}$ &  $2Bb_{1}\frac{E[1,0^{n+1}]}{E[1^{2},0^{n}]}$\\ \hline
%----------------------------- \RIGHTcircle to \Circle ------------------
$(\RIGHTcircle,n) \to (\Circle,n)$ & $Bb_{1}\frac{P_{N}(0^{n+2})}{P_{N}(1,0^{n+1})}$ &  $Bb_{1}\frac{E[0^{n+2}]}{E[1,0^{n+1}]}$\\ \hline
%--------------------------- \RIGHTcircle to \CIRCLE ------------------
$(\RIGHTcircle,n) \to (\CIRCLE,n)$ & $Bb_{0}\frac{P_{N}(1^{2},0^{n})}{P_{N}(1,0^{n+1})}$ & $Bb_{0}\frac{E[1^{2},0^{n}]}{E[1,0^{n+1}]}$\\ \hline
%-------------------------- \Circle to \bigtriangleup ------------------
$(\Circle,n) \to \bigtriangleup$ & $(\frac{N-S}{N})\frac{P_{N}(0^{n+1})}{P_{N}(0^{n+2})} + \frac{S}{N}$ &  $\frac{E[0^{n+1}]}{E[0^{n+2}]}$\\ \hline
%--------------------------\CIRCLE to \bigtriangleup ------------------
$(\CIRCLE,n)\to \bigtriangleup$ & $\frac{P_{N}(1,0^{n})}{P_{N}(1^{2},0^{n})} + \frac{S}{N}\frac{P_{N}(1,0^{n+1})}{P_{N}(1^{2},0^{n})}$ & $\frac{E[1,0^{n}]}{E[1^{2},0^{n}]}$\\ \hline
\end{tabular} 
\end{center}
\end{figure}
%%%%%%%%%%%%%%%%%%%%%%%%%%%%%%%%%%%%%%%%%%%%%%%%%%%%%%%%%%%%%%%%%%%%%%%%%%%%%%%%%%%%%%%%%%%%%%%%%%%      Applications, genealogical distance, Corollary

\begin{Proposition}[Conditioned genealogical distance]\label{PropositionGenalogicalDistance}\pointpar
\noindent Let $K = \{0,1\}$, $J = \{i,j\}$ and $\xi^{*} = y \in \{\Circle, \CIRCLE, \RIGHTcircle\}$. We have that
\begin{equation}
P(\inf\{r \ge 0: \overline{Y}^{2,N}_{r} = \bigtriangleup \} > t \,|\, \overline{Y}^{2,N}_{0} = (y,0) )= \overline{\mathbb{P}}_{\overline{\xi^{*}}}^{h}(\inf\{r \ge 0: (\overline{X}_{r})^{\mbox{\tiny $J$}}_{i} = (\overline{X}_{r})^{\mbox{\tiny $J$}}_{j} \} > t)
\end{equation}
for all $t \ge 0$.

This means 
\begin{equation}
P(\overline{Y}^{2,N}_{t} \not= \bigtriangleup  \,|\, \overline{Y}^{2,N}_{0} = (y,0))
\end{equation}
is the conditional probability that the genealogical distance of $i$ and $j$ in equilibrium (given the type information $y$) is greater than $2t$, where $y = \Circle$ means that both ($i$ and $j$) have type $0$, $y = \CIRCLE$  that both have type $1$ and $y = \RIGHTcircle$ that both have different types.
\end{Proposition}
%%%%%%%%%%%%%%%%%%%%%%%%%%%%%%%%%%%%%%%%%%%%%%%%%%%%%%%%%%%%%%%%%%%%%%%%%%%%%%%%%%%%%%%%%%%%%%%%%%%      Applications, genealogical distance

Now we let $N \to \infty$ and obtain that the sequence $\{\overline{Y}^{2,N}:N \in \mathbb{N}\}$ converges to a Markov jump process 
\begin{equation}
\overline{Y}^{2} = (\overline{Y}^{2}_{t})_{t \ge 0} \;\mbox{ on } \; \{\bigtriangleup \} \cup \{\Circle, \CIRCLE, \RIGHTcircle\} \times \mathbb{N}_{0}
\end{equation}
with transitions as listed in the third column of Figure \ref{FigureTransitionsY2}. Hence for $y  \in \{\Circle, \CIRCLE, \RIGHTcircle\}$ and $t \ge 0$,
\begin{equation}
P(\overline{Y}^{2}_{t} \not= \bigtriangleup  \,|\, \overline{Y}^{2}_{0} = (y,0))
\end{equation}
is the conditional probability that the  genealogical distance of two individuals (sampled from the tree-valued Fleming-Viot process in equilibrium) is greater than $2t$. 
%%%%%%%%%%%%%%%%%%%%%%%%%%%%%%%%%%%%%%%%%%%%%%%%%%%%%%%%%%%%%%%%%%%%%%%%%%%%%%%%%%%%%%%%%%%%%%%%%%%      Applications, genealogical distance
\bigskip

Finally, the aim is to investigate these conditional probabilities in order to study the tail distribution function of the genealogical distance of two individuals.

For this purpose define for each $(y,n) \in \{\Circle, \CIRCLE, \RIGHTcircle\} \times \mathbb{N}_{0} $, 
\begin{equation}
f_{t}(y,n) := P(\overline{Y}^{2}_{t} \not= \bigtriangleup  \,|\, \overline{Y}^{2}_{0} = (y,n)) \;\; \mbox{ for all }\;t \ge 0 \;,
\end{equation}
and set $f_{t}(\bigtriangleup) = 0$ for all $t \ge 0$. This means  
\begin{equation}
pf_{t}: = f_{t}(\Circle,0) E[0^{2}] + f_{t}(\CIRCLE,0) E[1^{2}] + f_{t}(\RIGHTcircle,0)2E[1,0]
\end{equation}
is the  probability that the  genealogical distance of two individuals sampled from the tree-valued Fleming-Viot process in equilibrium is greater than $2t$.

In the following proposition we shall analyze $f_{t}(y,n)$ and $pf_{t}$, where we contrast the case $S = 0$ (no selection) with the case $S > 0$. Remember also that $Bb_{0}$ is the mutation rate to type $0$ and $Bb_{1}$ the mutation rate to type $1$.
%%%%%%%%%%%%%%%%%%%%%%%%%%%%%%%%%%%%%%%%%%%%%%%%%%%%%%%%%%%%%%%%%%%%%%%%%%%%%%%%%%%%%%%%%%%%%%%%%%%      Applications, genealogical distance

\begin{Proposition}[Genealogical distance of two individuals]\label{Proposition}\pointpar
\noindent Let $K = \{0,1\}$, $Bb_{0} > 0$ and $Bb_{1} > 0$. Then the following holds:
\begin{enumerate}
\item {\bf The case $S = 0$ (without selection):} We have that
\begin{equation}
f_{t}(\Circle,0) = e^{-t}\li(1 + \frac{b_{1}(e^{-2Bt} - 1)}{1+2Bb_{0}}\re),\;\; 
f_{t}(\CIRCLE,0) = e^{-t}\li(1 + \frac{b_{0}(e^{-2Bt} - 1)}{1+2Bb_{1}}\re)
\end{equation}
and 
\begin{equation}
f_{t}(\RIGHTcircle,0) = e^{-t}\li(1 - \frac{(e^{-2Bt} - 1)}{2B}\re)
\end{equation}
for all $t \ge 0$. These explicit formulas imply: 
\begin{enumerate}
\item The genealogical distance is exponentially distributed.\\ Formally, $pf_{t} = e^{-t}$ for all $t \ge 0$.
\item The conditioned genealogical distance given the type information $\Circle$ is stochastically smaller than the conditioned genealogical distance given the type information $\CIRCLE$ if and only if the mutation rate to type $1$ is greater than to type $0$. \\ Formally,  $f_{t}(\Circle,0) <  f_{t}(\CIRCLE,0)$ for all $t > 0  \iff b_{1} > b_{0}$.
\item The conditioned genealogical distance is stochastically smaller than the genealogical distance if and only if the given types coincide.\\ Formally, $\max\li\{f_{t}(\Circle,0),  f_{t}(\CIRCLE,0)\re\} < e^{-t} < f_{t}(\RIGHTcircle,0)$ for all $t  > 0$.
\item The conditioned genealogical distance is exponentially distributed if the mutation rate $B$ tends to infinity.\\ Formally, $\lim_{B \to \infty}f_{t}(\Circle,0) = \lim_{B \to \infty}f_{t}(\CIRCLE,0) = \lim_{B \to \infty}f_{t}(\RIGHTcircle,0) = e^{-t}\;\; \forall\, t  \ge 0$.
\end{enumerate}
%----------------------------------------- with selection ---------------
\item {\bf The case $S > 0$ (with selection):}
\begin{enumerate}
\item We have that
\begin{equation}
\li.\frac{\partial f_{t}(\Circle,0)}{\partial t}\re|_{t = 0}  =  \frac{-E[0]}{E[0^{2}]}\;, \; \li.\frac{\partial f_{t}(\CIRCLE,0)}{\partial t}\re|_{t = 0}  = \frac{-E[1]}{E[1^{2}]} \mbox{ and } \li.\frac{\partial f_{t}(\RIGHTcircle,0)}{\partial t}\re|_{t = 0}  =  0\,.
\end{equation}
This means for each selection coefficient $S$ there is $\epsilon(S) > 0$ such that 
\begin{equation}
 \max\li\{f_{t}(\Circle,0),  f_{t}(\CIRCLE,0)\re\} <  e^{-t}  < f_{t}(\RIGHTcircle,0)\;\;\mbox{ for all }\;  0 <t < \epsilon(S) \;.
\end{equation} 
\item  We have that
\begin{equation}
pf_{0} = \li.-\frac{\partial pf_{t}}{\partial t}\re|_{t = 0} = \li.\frac{\partial^{2}pf_{t}}{\partial t^{2}}\re|_{t = 0} = 1 \;\mbox{ and }\;  \li.\frac{\partial^{3}pf_{t}}{\partial t^{3}}\re|_{t = 0}  = -\li(1 + 2S^{2}E[1,0]\re).
\end{equation}
This means for each selection coefficient $S$ there is $\epsilon(S) > 0$ such that 
\begin{equation}
pf_{t} <  e^{-t}  \;\; \mbox{ for all }\;  0 <t < \epsilon(S) \;. 
\end{equation}
\end{enumerate}
\end{enumerate}
\end{Proposition}
%%%%%%%%%%%%%%%%%%%%%%%%%%%%%%%%%%%%%%%%%%%%%%%%%%%%%%%%%%%%%%%%%%%%%%%%%%%%%%%%%%%%%%%%%%%%%%%%%%%      Applications, remark and preview how the proofs work

\begin{Remark}[Food for thought]\point If $S > 0$ and $Bb_{0} = Bb_{1} = \tfrac{1}{2}$, then one has that
\begin{equation}
\frac{-E[0]}{E[0^{2}]} = \frac{(e^{2S}-2S-1)S}{1-e^{2S}+2S(S+1)} < \frac{(1-e^{2S}+2Se^{2S})S}{1-e^{2S}+2Se^{2S}(S-1)} = \frac{-E[1]}{E[1^{2}]} \;.
\end{equation}
This means the conditioned genealogical distance given that the two individuals are fit is {\rm not} stochastically smaller than the  conditioned genealogical distance given that the two individuals are unfit.
\end{Remark}

\bigskip
Before we prove (see Section \ref{secProofs}) the theorems and propositions stated in this section, in Section \ref{secHBP} we will provide the main technical tool, the Feynman-Kac duality between the HMM and the HBP. 
%%%%%%%%%%%%%%%%%%%%%%%%%%%%%%%%%%%%%%%%%%%%%%%%%%%%%%%%%%%%%%%%%%%%%%%%%%%%%%%%%%%%%%%%%%%%%%%%%%%      Main Tool: The HBP

\setcounter{equation}{0}
\section{The key tool: The historical backward process (HBP)}\label{secHBP}
\noindent  Augmenting its past history to a stochastic process is a powerful tool (see \cite{D93} and \cite{DP91}) we shall use in the present paper to describe and prove the relation between the  HMM and the sample-paths of the BP.

For this purpose we define in this section the HBP, the {\it path process} associated to the BP, which contains at time $t$ the sample paths of the BP up to time $t$. The first task is to give a description of the HBP in terms of a Borel strong Markov process. Then we formulate the relation between the HMM and the HBP which is based on a generator relation.
%%%%%%%%%%%%%%%%%%%%%%%%%%%%%%%%%%%%%%%%%%%%%%%%%%%%%%%%%%%%%%%%%%%%%%%%%%%%%%%%%%%%%%%%%%%%%%%%%%%      The HBP, Definition and analytical characterization

\subsection{Definition and analytical characterization}
\noindent In this subsection we introduce and characterize the HBP, more precisely the {\sc Time}-space process of the HBP, as it appears in \cite{Per02}, Proposition II.2.5.
%%%%%%%%%%%%%%%%%%%%%%%%%%%%%%%%%%%%%%%%%%%%%%%%%%%%%%%%%%%%%%%%%%%%%%%%%%%%%%%%%%%%%%%%%%%%%%%%%%% The associated path process of a {\duttfamily BSMP}

Let
\begin{equation}
 \overline{X} =  \li( \overline{\Omega} , \overline{\mathcal{F}}, (\overline{\mathcal{F}}_{t})_{t \ge 0},(\overline{X}_{t})_{t \ge 0}, \{\overline{\mathbb{P}}_{\overline{\eta}} : \overline{\eta} \in \overline{\mathcal{E}}\}\re) 
\end{equation}
be the Borel strong Markov process (depending on $J \subset I$) that describes the BP and is characterized as given in Subsection \ref{subsecAnalyticalCharacterization} (hence recall $\overline{L}$ and $\overline{\mathcal{K}}$). In order to define the path process associated to $\overline{X}$ let
\begin{equation}
\overline{\mathcal{E}}^{\diamond} = \li\{\overline{\eta}^{\diamond}  = \li((\overline{\eta}^{\diamond})^{\mbox{\tiny \sc Time}} ,(\overline{\eta}^{\diamond})^{\mbox{\tiny $\mathcal{D}$}}\re) \in \mathbb{R} \times \mathcal{D}(\mathbb{R}, \overline{\mathcal{E}}): (\overline{\eta}^{\diamond})^{\mbox{\tiny $\mathcal{D}$}}_{s} = \overline{\eta}^{*} \; \mbox{ for all } s \ge (\overline{\eta}^{\diamond})^{\mbox{\tiny \sc Time}}  \re\}\;,
\end{equation}
where
\begin{equation}
\overline{\eta}^{*} :=  (\overline{\eta}^{\diamond})^{\mbox{\tiny $\mathcal{D}$}}_{(\overline{\eta}^{\diamond})^{\mbox{\tiny \sc Time}}} \in \overline{\mathcal{E}}
\end{equation}
is the projection on the state of the HBP at {\sc Time} $(\overline{\eta}^{\diamond})^{\mbox{\tiny \sc Time}}$.
%%%%%%%%%%%%%%%%%%%%%%%%%%%%%%%%%%%%%%%%%%%%%%%%%%%%%%%%%%%%%%%%%%%%%%%%%%%%%%%%%%%%%%%%%%%%%%%%%%% HBP, Definition

Moreover, let $\overline{\Omega}^{\diamond} = \mathcal{D}([0,\infty), \overline{\mathcal{E}}^{\diamond})$, $(\overline{X}^{\diamond}_{t})_{t \ge 0}$ be the canonical $\overline{\mathcal{E}}^{\diamond}$-valued process on $(\overline{\Omega}^{\diamond}, \overline{\mathcal{F}}^{\diamond})$ with canonical right continuous filtration $(\overline{\mathcal{F}}^{\diamond}_{t})_{t \ge 0}$ and $\{\overline{\mathbb{P}}^{\diamond}_{\overline{\eta}^{\diamond}} : \overline{\eta}^{\diamond} \in \overline{\mathcal{E}}^{\diamond}\}$ the collection of probability measures on $\overline{\Omega}^{\diamond}$ defined by
\begin{equation}
\overline{\mathbb{P}}^{\diamond}_{\overline{\eta}^{\diamond}} : \overline{\mathcal{F}}^{\diamond} \to [0,1],\; C  \mapsto  \overline{\mathbb{P}}_{\overline{\eta}}\li( \li((\overline{\eta}^{\diamond})^{\mbox{\tiny \sc Time}} + t, \li[(\overline{\eta}^{\diamond})^{\mbox{\tiny $\mathcal{D}$}}\, \li|\, (\overline{\eta}^{\diamond})^{\mbox{\tiny \sc Time}} \, \li| \, \li(\overline{X}_{\cdot\wedge t}\re) \re. \re. \re]\re)_{t \ge 0} \in  C\re)
\end{equation}
where for $\overline{\omega} \in \overline{\Omega}$ and $t \ge 0$,
\begin{equation}
\li[(\overline{\eta}^{\diamond})^{\mbox{\tiny $\mathcal{D}$}}\, \li|\, (\overline{\eta}^{\diamond})^{\mbox{\tiny \sc Time}} \, \li| \, (\overline{\omega}_{\cdot\wedge t}) \re. \re. \re]\in  \mathcal{D}(\mathbb{R}, \overline{\mathcal{E}})
\end{equation}
is the path we obtain by continuing the path $(\overline{\eta}^{\diamond})^{\mbox{\tiny $\mathcal{D}$}}$ with the path $(\overline{\omega}_{\cdot\wedge t}) = (\overline{\omega}_{s\wedge t})_{s \ge 0}$, the path $\overline{\omega}$ stopped at time $t$, from {\sc Time} $(\overline{\eta}^{\diamond})^{\mbox{\tiny \sc Time}}$ on. Formally,  the map
\begin{equation}
(\overline{\eta}^{\diamond}, \overline{\omega}) \mapsto \li[(\overline{\eta}^{\diamond})^{\mbox{\tiny $\mathcal{D}$}}\, \li|\, (\overline{\eta}^{\diamond})^{\mbox{\tiny \sc Time}} \, \li| \, (\overline{\omega}_{\cdot\wedge t}) \re. \re. \re] \; \; \mbox{ from }\; \overline{\mathcal{E}}^{\diamond} \times \overline{\Omega} \;\mbox{ to }\; \mathcal{D}(\mathbb{R}, \overline{\mathcal{E}})
\end{equation}
is defined by
\begin{equation}
\li[(\overline{\eta}^{\diamond})^{\mbox{\tiny $\mathcal{D}$}}\, \li|\, (\overline{\eta}^{\diamond})^{\mbox{\tiny \sc Time}} \, \li| \, (\overline{\omega}_{\cdot\wedge t}) \re. \re. \re]_{s} = \li\{ \begin{array}{ccc}(\overline{\eta}^{\diamond})^{\mbox{\tiny $\mathcal{D}$}}_{s} &,&  s < (\overline{\eta}^{\diamond})^{\mbox{\tiny \sc Time}} \\ \overline{\omega}_{(s-(\overline{\eta}^{\diamond})^{\mbox{\tiny \sc Time}})\wedge t} &,&  s \ge (\overline{\eta}^{\diamond})^{\mbox{\tiny \sc Time}} \end{array} \re. \;.
\end{equation}
%%%%%%%%%%%%%%%%%%%%%%%%%%%%%%%%%%%%%%%%%%%%%%%%%%%%%%%%%%%%%%%%%%%%%%%%%%%%%%%%%%%%%%%%%%%%%%%%%%% HBP, Definition and relation to BP

Since $\overline{L}$, the generator corresponding to $\overline{X}$, is a bounded linear operator, we obtain that
\begin{equation}
 \overline{X}^{\diamond} =  \li(\overline{\Omega}^{\diamond} , \overline{\mathcal{F}}^{\diamond}, (\overline{\mathcal{F}}^{\diamond}_{t})_{t \ge 0},(\overline{X}^{\diamond}_{t})_{t \ge 0}, \{\overline{\mathbb{P}}^{\diamond}_{\overline{\eta}^{\diamond}} : \overline{\eta}^{\diamond} \in \overline{\mathcal{E}}^{\diamond}\} \re) 
\end{equation}
is a Borel strong Markov process such that for $T \ge 0$ and $\overline{\eta} \in \overline{\mathcal{E}}$,
\begin{equation}
\overline{\mathbb{P}}_{\overline{\eta}} \circ \li( ( \overline{X}_{t})_{t \in \li[0,T\re]}\re)^{-1} = \overline{\mathbb{P}}^{\diamond}_{\overline{\eta}^{\diamond}} \circ \li( ((\overline{X}^{\diamond}_{T})^{\mbox{\tiny $\mathcal{D}$}}_{t})_{t \in \li[(\overline{X}^{\diamond}_{0})^{\mbox{\tiny \sc Time}},(\overline{X}^{\diamond}_{T})^{\mbox{\tiny \sc Time}}\re]}\re)^{-1}
\end{equation}
if $\overline{\eta}^{*} = \overline{\eta}$. 
%%%%%%%%%%%%%%%%%%%%%%%%%%%%%%%%%%%%%%%%%%%%%%%%%%%%%%%%%%%%%%%%%%%%%%%%%%%%%%%%%%%%%%%%%%%%%%%%%%% HBP, analytic characterization

Finally we give an analytic characterization we use in the proof of the  Feynman-Kac duality between the HMM and the HBP. Namely, we specify a set $\overline{\mathcal{A}}^{\diamond} \subset M_{b}(\overline{\mathcal{E}}^{\diamond})$ (measurable and bounded functions on $\overline{\mathcal{E}}^{\diamond}$) and a map 
\begin{equation}
\overline{L}^{\diamond}:\overline{\mathcal{A}}^{\diamond} \to M_{b}(\overline{\mathcal{E}}^{\diamond})
\end{equation}
such that for each $\overline{\eta}^{\diamond} \in \overline{\mathcal{E}}^{\diamond}$, $\overline{\mathbb{P}}^{\diamond}_{\overline{\eta}^{\diamond}}$ is a solution of the $\overline{\Omega}^{\diamond}$-martingale problem for $(\overline{L}^{\diamond}, \delta_{\overline{\eta}^{\diamond}})$ w.r.t$.$ $\overline{\mathcal{A}}^{\diamond}$.

The set $\overline{\mathcal{A}}^{\diamond}$ is specified by
\begin{eqnarray}
f \in \overline{\mathcal{A}}^{\diamond} &\iff& g(\overline{\eta}^{\diamond}) := \lim_{\epsilon \downarrow 0} \frac{f((\overline{\eta}^{\diamond})^{\mbox{\tiny \sc Time}} + \epsilon, (\overline{\eta}^{\diamond})^{\mbox{\tiny $\mathcal{D}$}}) - f(\overline{\eta}^{\diamond})}{\epsilon} \;\mbox{ exists for all } \overline{\eta}^{\diamond} \in \overline{\mathcal{E}}^{\diamond},\;\; \\
&&f((\overline{\eta}^{\diamond})^{\mbox{\tiny \sc Time}} + t, (\overline{\eta}^{\diamond})^{\mbox{\tiny $\mathcal{D}$}}) - f(\overline{\eta}^{\diamond}) = \int\limits_{0}^{t} g((\overline{\eta}^{\diamond})^{\mbox{\tiny \sc Time}} + s, (\overline{\eta}^{\diamond})^{\mbox{\tiny $\mathcal{D}$}}) ds\; \forall t \ge 0.
\end{eqnarray} 
The map $\overline{L}^{\diamond} = \overline{L}^{\diamond,\mbox{\tiny \sc Time}} + \overline{L}^{\diamond, \mbox{\tiny $\mathcal{D}$}}:\overline{\mathcal{A}}^{\diamond} \to M_{b}(\overline{\mathcal{E}}^{\diamond})$ is defined by
\begin{equation}
\overline{L}^{\diamond,\mbox{\tiny \sc Time}}f(\overline{\eta}^{\diamond}) = \lim_{\epsilon \downarrow 0} \frac{f((\overline{\eta}^{\diamond})^{\mbox{\tiny \sc Time}} + \epsilon, (\overline{\eta}^{\diamond})^{\mbox{\tiny $\mathcal{D}$}}) - f(\overline{\eta}^{\diamond})}{\epsilon}
\end{equation}
and
\begin{equation}
\overline{L}^{\diamond,\mbox{\tiny $\mathcal{D}$}}f(\overline{\eta}^{\diamond}) = \sum_{\overline{\zeta} \in \overline{\mathcal{E}}} \overline{\mathcal{K}}(\overline{\eta}^{*}, \overline{\zeta})\li[f(\overline{\eta}^{\diamond:\overline{\zeta}}) - f(\overline{\eta}^{\diamond})\re]\;,
\end{equation}
where $\overline{\eta}^{\diamond:\overline{\zeta}} \in \overline{\mathcal{E}}^{\diamond}$ is the element given by
\begin{equation}
(\overline{\eta}^{\diamond:\overline{\zeta}})^{\mbox{\tiny \sc Time}} := (\overline{\eta}^{\diamond})^{\mbox{\tiny \sc Time}} \mbox{ and }(\overline{\eta}^{\diamond:\overline{\zeta}})^{\mbox{\tiny \sc $\mathcal{D}$}}_{s} := \li\{ \begin{array}{ccl} (\overline{\eta}^{\diamond})^{\mbox{\tiny \sc $\mathcal{D}$}}_{s} &,& s <  (\overline{\eta}^{\diamond})^{\mbox{\tiny \sc Time}}\\ \overline{\zeta} &,&  s \ge  (\overline{\eta}^{\diamond})^{\mbox{\tiny \sc Time}}\end{array} \re. \;.
\end{equation}
%%%%%%%%%%%%%%%%%%%%%%%%%%%%%%%%%%%%%%%%%%%%%%%%%%%%%%%%%%%%%%%%%%%%%%%%%%%%%%%%%%%%%%%%%%%%%%%%%%%      The HBP, duality

\subsection{Feynman-Kac duality between the HMM and the HBP}\label{subsecFeynmanKacHMMandHBP}
\noindent In this subsection we state the Feynman-Kac duality, or rather a whole collection of Feynman-Kac dualities, between the HMM and the HBP, where it is helpful to recall the notation introduced in Subsection \ref{subsecAnalyticalCharacterization}. For this purpose we introduce a collection of duality functions denoted by $\mathcal{H}$ and the Feynman-Kac function denoted by $V^{\diamond}$. 

An element in $\mathcal{H}$ is a bounded and continuous function
\begin{equation}\label{equationH}
H: \mathcal{E} \times \overline{\mathcal{E}}^{\diamond} \to \mathbb{R}
\end{equation}
(depending on $J \subset I$) of the form
\begin{equation}
H(\eta, \overline{\eta}^{\diamond}) = g(\eta^{\mbox{\tiny \sc Time}} + (\overline{\eta}^{\diamond})^{\mbox{\tiny \sc Time}}) H^{*}(\eta, \overline{\eta}^{*}) H^{\mbox{\tiny $\mathcal{D}$}}(\eta, \overline{\eta}^{\diamond})\;,
\end{equation}
where $g \in C^{1}_{b}(\mathbb{R})$, $H^{*}$ as given in (\ref{equationHstar}) and
\begin{equation}
H^{\mbox{\tiny $\mathcal{D}$}}(\eta, \overline{\eta}^{\diamond}) =  \prod_{\gamma \in \Gamma(\overline{\eta}^{*})}\prod_{i \in \gamma}\prod_{n = 1}^{m_{i}}\int_{r_{i}^{n}}^{t_{i}^{n}} \li(F_{i}^{n}(\eta^{\mbox{\tiny $\mathcal{D}$}}_{l,s})\mathbbm{1}\{s \le \eta^{\mbox{\tiny \sc Time}}\} + F_{i}^{n}( ((\overline{\eta}^{\diamond})^{\mbox{\tiny $\mathcal{D}$}}_{\overline{s}})^{\mbox{\tiny $J$}}_{i})\mathbbm{1}\{s > \eta^{\mbox{\tiny \sc Time}}\}\re)ds 
\end{equation}
with 
\begin{equation}
l = (\overline{\eta}^{*})^{\mbox{\tiny $J$}}_{\gamma,I}\; \mbox{ and }\; \overline{s} = \eta^{\mbox{\tiny \sc Time}} + (\overline{\eta}^{\diamond})^{\mbox{\tiny \sc Time}} - s\;,
\end{equation}
where $r_{j}^{1} < t_{j}^{1} < \dots < r_{j}^{m_{j}} < t_{j}^{m_{j}}$ and $F_{j}^{n} \in C_{b}(K \times I)$ for all $j \in J$. 

In the function $H^{\mbox{\tiny $\mathcal{D}$}}$ we evaluate the paths (each of these in terms of integrals over different time intervals) that arise by connecting certain extended ancestral lines with suitable paths of the HBP in a special way depending on the {\sc Time}. Namely, for each $\gamma \in \Gamma(\overline{\eta}^{*})$ (which represents an ancestor occupying the life-site $l = (\overline{\eta}^{*})^{\mbox{\tiny $J$}}_{\gamma,I}$) and for each $i \in \gamma$ (which represents a descendent of $l$)  we  connect  $\eta^{\mbox{\tiny $\mathcal{D}$}}_{l,\cdot}$ (from $-\infty$ up to $\eta^{\mbox{\tiny \sc Time}}$) with $((\overline{\eta}^{\diamond})^{\mbox{\tiny $\mathcal{D}$}}_{\cdot})^{\mbox{\tiny $J$}}_{i}$ (from $(\overline{\eta}^{\diamond})^{\mbox{\tiny \sc Time}}$ back to $-\infty$). 

\begin{Remark}\point If the element $\overline{\eta}^{\diamond} \in \overline{\mathcal{E}}^{\diamond}$ satisfies $(\overline{\eta}^{\diamond})^{\mbox{\tiny $\mathcal{D}$}}_{s} = \overline{\xi^{*}}$ (for all $s \in \mathbb{R}$) with $\xi^{*} \in K^{J}$, then
\begin{equation}
H(\eta, \overline{\eta}^{\diamond}) = g(\eta^{\mbox{\tiny \sc Time}}) \li(\prod_{j \in J}\prod_{n = 1}^{m_{j}}\int_{r_{j}^{n}}^{t_{j}^{n}}F_{j}^{n}(\eta^{\mbox{\tiny $\mathcal{D}$}}_{j,s})ds\re) \mathbbm{1}\{ (\eta^{*}_{j})_{j \in J} = \xi^{*}\}\;.
\end{equation}
This means, recall (\ref{equationFunctionsInA}), that for this choice of the parameters the function $H$ is an element of $\mathcal{A}$, but note that a general $H$ is only an element of the {\it bp}-closure of $\mathcal{A}$.
\end{Remark}

The Feynman-Kac function is given by
\begin{equation}\label{equationVdiamond}
V^{\diamond}: \overline{\mathcal{E}}^{\diamond} \to \mathbb{R}, \; \overline{\eta}^{\diamond} \mapsto V(\overline{\eta}^{*})\;,
\end{equation}
where $V$ is given as in (\ref{equationV}).

\begin{Theorem}[Feynman-Kac duality for the HMM] \label{TheoremFeynmanKacDualityHMMandHBP}\pointpar
\noindent Let $J \subset I$, $\mu$ be a general distribution on $\mathcal{E}$, $(X_{t})_{t \ge 0}$ the canonical $\mathcal{E}$-valued process on $(\Omega, \mathcal{F})$ and $L$ and $\mathcal{A}$ as in Subsection \ref{subsecAnalyticalCharacterization}. In addition, let $\overline{X}^{\diamond}$ be the Borel strong Markov process that describes the HBP, $V^{\diamond}$ as in (\ref{equationVdiamond}) and $H$ a function in $\mathcal{H}$.

For any solution $\mathbb{Q}$ of the $\Omega$-martingale problem for $(L,\mu)$ w.r.t$.$ $\mathcal{A}$ we have that
\begin{equation}
\int\limits_{\Omega}H(X_{t}(\omega), \overline{\eta}^{\diamond})\,\mathbb{Q}(d\omega) = \overline{\mathbb{E}}^{\diamond}_{\overline{\eta}^{\diamond}}\li[ \li(\int\limits_{\mathcal{E}}H(\eta, \overline{X}^{\diamond}_{t})\mu(d\eta)\re)e^{\int_{0}^{t} V^{\diamond}(\overline{X}^{\diamond}_{s})ds}\re]
\end{equation}
for all $\overline{\eta}^{\diamond} \in \overline{\mathcal{E}}^{\diamond}$ and all $t \ge 0$.
\end{Theorem}
%%%%%%%%%%%%%%%%%%%%%%%%%%%%%%%%%%%%%%%%%%%%%%%%%%%%%%%%%%%%%%%%%%%%%%%%%%%%%%%%%%%%%%%%%%%%%%%%%%%      Proofs, 

\setcounter{equation}{0}
\section{Proofs}\label{secProofs}
\noindent This section contains all proofs. We start with the proof of the Feynman-Kac duality between the HMM and the HBP  (Theorem \ref{TheoremFeynmanKacDualityHMMandHBP}). This duality ensures the uniqueness of the martingale problem which is the key property to obtain the statements of Theorem \ref{TheoremAnalyticalHMM}. In addition, we obtain the Feynman-Kac duality between the type information of the HMM and the BP (Proposition \ref{PropositionDualityTypInfoHMM}) and the stochastic representation for the extended ancestral lines (Theorem \ref{TheoremStochasticRepresentation}). 

Then we show the analytical characterization of the transformed BP (Theorem \ref{TheoremTransformedBP}), where we shall use that the BP is a Borel strong Markov process with finite state space and bounded generator. After that we apply the transformed BP together with the stochastic representation for the extended ancestral lines in order to obtain the strong stochastic representation (Theorem \ref{TheoremStrongStochasticRepresentation}) which in turn can be used together with the analytical characterization of the transformed BP to prove the longtime behaviour of the HMM (Theorem \ref{TheoremLongtime}).

Finally we give the proofs of Proposition \ref{PropositionCAT} (representation for the stationary type distribution of the CAT), Proposition \ref{PropositionGenalogicalDistance} (representation for the conditioned genealogical distance in equilibrium) and  Proposition \ref{Proposition} (the genealogical distance of two individuals in equilibrium).
%%%%%%%%%%%%%%%%%%%%%%%%%%%%%%%%%%%%%%%%%%%%%%%%%%%%%%%%%%%%%%%%%%%%%%%%%%%%%%%%%%%%%%%%%%%%%%%%%%%      Proofs, Proof of Theorem \ref{TheoremFeynmanKacDualityHMMandHBP}

\subsection{Proof of Theorem \ref{TheoremFeynmanKacDualityHMMandHBP}}
\noindent Let $\mu$ be a general distribution on $\mathcal{E}$, $J \subset I$, $H \in \mathcal{H}$ and $\overline{\xi}^{\diamond} \in \overline{\mathcal{E}}^{\diamond}$. The proof is carried out in two steps. In the first step we verify the assumptions of Theorem 4.4.11 in \cite{EK86} in order to obtain that
\begin{eqnarray}
&& \int\limits_{\Omega}H(X_{t}, \overline{\xi}^{\diamond})d\mathbb{Q} -  \overline{\mathbb{E}}^{\diamond}_{\overline{\xi}^{\diamond}}\li[ \li(\int\limits_{\mathcal{E}}H(\eta, \overline{X}^{\diamond}_{t})\mu(d\eta)\re)e^{\int_{0}^{t} V^{\diamond}(\overline{X}^{\diamond}_{s})ds}\re]\\
&=& \int\limits_{0}^{t}\li(\;\int\limits_{\Omega \times \overline{\Omega}^{\diamond}} \mathbb{H}(X_{s},\overline{X}^{\diamond}_{t-s})e^{\int_{0}^{t-s}V^{\diamond}(\overline{X}^{\diamond}_{r})dr}d(\mathbb{Q} \otimes \overline{\mathbb{P}}^{\diamond}_{\overline{\xi}^{\diamond}})\re)ds\;,
\end{eqnarray}
where
\begin{equation}
\mathbb{H}(\eta,\overline{\eta}^{\diamond}) = [LH(\, \cdot \,,\overline{\eta}^{\diamond})](\eta) - [\overline{L}^{\diamond}H(\eta,\, \cdot\,)](\overline{\eta}^{\diamond}) - V^{\diamond}(\overline{\eta}^{\diamond})H(\eta,\overline{\eta}^{\diamond})
\end{equation}
for all $\eta \in \mathcal{E}$ and all $\overline{\eta}^{\diamond} \in \overline{\mathcal{E}}^{\diamond}$. In the second step we show that
\begin{equation}
\mathbb{H}(\eta,\overline{\eta}^{\diamond}) = 0
\end{equation}
for all $\eta \in \mathcal{E}$ and all $\overline{\eta}^{\diamond} \in \overline{\mathcal{E}}^{\diamond}$.
%%%%%%%%%%%%%%%%%%%%%%%%%%%%%%%%%%%%%%%%%%%%%%%%%%%%%%%%%%%%%%%%%%%%%%%%%%%%%%%%%%%%%%%%%%%%%%%%%%%      Proof of Theorem \ref{TheoremFeynmanKacDualityHMMandHBP}

Throughout these two steps we assume (for simplicity) that $g \equiv 1$ and $m_{j} = 1$ for all $j \in J$, that is,  the function $H(\eta, \overline{\eta}^{\diamond})$ has the form
\begin{equation}
H^{*}(\eta, \overline{\eta}^{*})  \prod_{\gamma \in \Gamma(\overline{\eta}^{*})} \prod_{i \in \gamma}\int_{r_{i}}^{t_{i}}(F_{i}(\eta^{\mbox{\tiny $\mathcal{D}$}}_{l,s})\mathbbm{1}\{s \le \eta^{\mbox{\tiny \sc Time}}\} + F_{i}( ((\overline{\eta}^{\diamond})^{\mbox{\tiny $\mathcal{D}$}}_{\overline{s}})^{\mbox{\tiny $J$}}_{i})\mathbbm{1}\{s > \eta^{\mbox{\tiny \sc Time}}\})ds \;,
\end{equation}
where $l = (\overline{\eta}^{*})^{\mbox{\tiny $J$}}_{\gamma,I}$, $ \overline{s} = \eta^{\mbox{\tiny \sc Time}} + (\overline{\eta}^{\diamond})^{\mbox{\tiny \sc Time}} - s$ and $H^{*}$ is the duality function  defined in (\ref{equationHstar}). \bigskip
%%%%%%%%%%%%%%%%%%%%%%%%%%%%%%%%%%%%%%%%%%%%%%%%%%%%%%%%%%%%%%%%%%%%%%%%%%%%%%%%%%%%%%%%%%%%%%%%%%%      Proof of Theorem \ref{TheoremFeynmanKacDualityHMMandHBP}, Step 1

{\it The first step:} We have to show that:
\begin{enumerate}
\item For each $T >0$ the random variable  
\begin{equation}
\sup_{s,r,t \le T}(|[LH(\cdot,\overline{X}^{\diamond}_{t})](X_{s})| + |[\overline{L}^{\diamond}H(X_{s},\cdot)](\overline{X}^{\diamond}_{t})| +|H(X_{s},\overline{X}^{\diamond}_{t})|)(|V(\overline{X}^{\diamond}_{r})| +1)
\end{equation}
is integrable with respect to $\mathbb{Q} \otimes \overline{\mathbb{P}}^{\diamond}_{\overline{\xi}^{\diamond}}$. But this holds since
\begin{equation}
\sup\li\{(|[LH(\cdot,\overline{\eta}^{\diamond})](\eta)| + |[\overline{L}^{\diamond}H(\eta,\cdot)](\overline{\eta}^{\diamond})| +\|H\|_{\infty}) (\|V\|_{\infty}
+1): \eta \in \mathcal{E}, \overline{\eta}^{\diamond} \in \overline{\mathcal{E}}^{\diamond} \re\}
\end{equation}
is bounded by a real number depending on $\|H\|_{\infty}$, $N$, $B$ and $S$.
\item The probability measure $\overline{\mathbb{P}}^{\diamond}_{\overline{\xi}^{\diamond}}$ is a solution of the $\overline{\Omega}^{\diamond}$-martingale problem for $(\overline{L}^{\diamond}, \delta_{\overline{\xi}^{\diamond}})$ w.r.t$.$ $\{\overline{\eta}^{\diamond} \mapsto H(\eta,\overline{\eta}^{\diamond}): \eta \in \mathcal{E}\}$. 
\item The probability measure $\mathbb{Q}$ is a solution of the $\Omega$-martingale problem for $(L, \mu)$ w.r.t$.$ $\{\eta \mapsto H(\eta,\overline{\eta}^{\diamond}): \overline{\eta}^{\diamond} \in \overline{\mathcal{E}}^{\diamond}\}$.
\end{enumerate}
%---------------------------------------------------------------------------
In order to see 2 and 3 we rewrite $H(\eta, \overline{\eta}^{\diamond})$ as
\begin{equation}
H^{*}(\eta, \overline{\eta}^{*}) \prod_{\gamma \in \Gamma(\overline{\eta}^{*})}\prod_{i \in \gamma}\li(\int_{r_{i}}^{t_{i}} [F_{i}(\eta^{\mbox{\tiny $\mathcal{D}$}}_{l,s}) + F_{i}( ((\overline{\eta}^{\diamond})^{\mbox{\tiny $\mathcal{D}$}}_{\eta^{\mbox{\tiny \sc Time}} + (\overline{\eta}^{\diamond})^{\mbox{\tiny \sc Time}} - s})^{\mbox{\tiny $J$}}_{i})]ds  - F_{i}(\eta^{*}_{l})(t_{i} -r_{i})\re)
\end{equation}
which is possible due to the form of $H^{*}$ and the fact that  $\eta^{*}_{l} = (\overline{\eta}^{*})^{\mbox{\tiny $J$}}_{\gamma,K}$ for each $\gamma \in \Gamma(\overline{\eta}^{*})$. With this representation one  first of all gets that $\{ \overline{\eta}^{\diamond} \mapsto H(\eta,\overline{\eta}^{\diamond}): \eta \in \mathcal{E}\} \subset \overline{\mathcal{A}}^{\diamond}$ which yields 2. Furthermore, one can recognize that in general  $\{\eta \mapsto H(\eta,\overline{\eta}^{\diamond}): \overline{\eta}^{\diamond} \in \overline{\mathcal{E}}^{\diamond}\}$ is indeed not a subset of $\mathcal{A}$, but a subset of the {\it bp}-closure of $\mathcal{A}$ which gives 3. 
%%%%%%%%%%%%%%%%%%%%%%%%%%%%%%%%%%%%%%%%%%%%%%%%%%%%%%%%%%%%%%%%%%%%%%%%%%%%%%%%%%%%%%%%%%%%%%%%%%%      Proof of Theorem \ref{TheoremFeynmanKacDualityHMMandHBP}, Step 1

Finally, we obtain that
\begin{equation}
[L^{\mbox{\tiny \sc Time}}H(\, \cdot \,, \overline{\eta}^{\diamond})](\eta) = [\overline{L}^{\diamond, \mbox{\tiny \sc Time}}H(\eta, \, \cdot \,)](\overline{\eta}^{\diamond})
\end{equation}
for all $\eta \in \mathcal{E}$ and all $\overline{\eta}^{\diamond} \in \overline{\mathcal{E}}^{\diamond}$. Hence in the second step we have to show that
\begin{equation}\label{equationSecondStepThm6}
[L^{\mbox{\tiny $\mathcal{D}$}}H(\, \cdot \,,\overline{\eta}^{\diamond})](\eta) = \sum_{\overline{\zeta} \in \overline{\mathcal{E}}} \overline{\mathcal{K}}(\overline{\eta}^{*}, \overline{\zeta})\li[H(\eta,\overline{\eta}^{\diamond:\overline{\zeta}}) - H(\eta,\overline{\eta}^{\diamond})\re] + V(\overline{\eta}^{*})
H(\eta,\overline{\eta}^{\diamond})
\end{equation}
for all $\eta \in \mathcal{E}$ and all $\overline{\eta}^{\diamond} \in \overline{\mathcal{E}}^{\diamond}$, where $V$ is the Feynman-Kac term defined in (\ref{equationV}).\bigskip
%%%%%%%%%%%%%%%%%%%%%%%%%%%%%%%%%%%%%%%%%%%%%%%%%%%%%%%%%%%%%%%%%%%%%%%%%%%%%%%%%%%%%%%%%%%%%%%%%%%      Proof of Theorem \ref{TheoremFeynmanKacDualityHMMandHBP}, Step 2

{\it The second step:} For showing (\ref{equationSecondStepThm6}) we write $\overline{\eta}$ instead of $\overline{\eta}^{*}$ (projection on the BP). Furthermore, recall the notation used in the Definitions \ref{DefinitionHMM} and \ref{DefinitionBP} and observe that $H(\eta, \overline{\eta}^{\diamond})$  has the form
\begin{equation}
\li(\prod_{\gamma \in \Gamma(\overline{\eta})}  \mathbbm{1}\{\eta^{*}_{l} =\overline{\eta}^{\mbox{\tiny $J$}}_{\gamma,K}\} H^{\mbox{\tiny $\mathcal{D}$}}_{\eta, \overline{\eta}^{\diamond}}(l,\gamma)\re)\li(\prod_{i \in \tilde{\Gamma}(\overline{\eta})} \mathbbm{1}\{\eta^{*}_{i} \in \overline{\eta}^{\mbox{\tiny $I$}}_{i}\}\re)
\end{equation}
with 
\begin{equation}
H^{\mbox{\tiny $\mathcal{D}$}}_{\eta, \overline{\eta}^{\diamond}}(l,\gamma) := \prod_{i \in \gamma}\int\limits_{r_{i}}^{t_{i}}(F_{i}(\eta^{\mbox{\tiny $\mathcal{D}$}}_{l,s})\mathbbm{1}\{s \le \eta^{\mbox{\tiny \sc Time}}\} + F_{i}(((\overline{\eta}^{\diamond})^{\mbox{\tiny $\mathcal{D}$}}_{\overline{s}})^{\mbox{\tiny $J$}}_{i})\mathbbm{1}\{s > \eta^{\mbox{\tiny \sc Time}}\})ds
\end{equation}
and $l = (\overline{\eta})^{\mbox{\tiny $J$}}_{\gamma,I}$.
%%%%%%%%%%%%%%%%%%%%%%%%%%%%%%%%%%%%%%%%%%%%%%%%%%%%%%%%%%%%%%%%%%%%%%%%%%%%%%%%%%%%%%%%%%%%%%%%%%%      Proof of Theorem \ref{TheoremFeynmanKacDualityHMMandHBP}, Step 2

The first task is to decompose the left hand side of (\ref{equationSecondStepThm6}) according to the different kinds of transitions (\ref{itemdualmutationJ}, \ref{itemdualmutationI}, \ref{itemdualresamplingJJ}, \ref{itemdualresamplingJI} \ref{itemdualresamplingIJ} and \ref{itemdualresamplingII}) of the BP. With $l = (\overline{\eta})^{\mbox{\tiny $J$}}_{\gamma,I}$ and $l' = (\overline{\eta})^{\mbox{\tiny $J$}}_{\gamma',I}$ we have (observe that some terms cancel each other out) 

\begin{eqnarray} 
&& [L^{\mbox{\tiny $\mathcal{D}$}}H(\,\cdot \,,\overline{\eta}^{\diamond})](\eta)\\
\label{equationLHSmutationJ}&=\stackrel{\ref{itemdualmutationJ}}{}& B\sum_{\gamma \in \Gamma(\overline{\eta})}\sum_{u \in K} b(\eta^{*}_{l},u )\li[H(\eta^{l;u}, \overline{\eta}^{\diamond})-H(\eta, \overline{\eta}^{\diamond})\re] \\
\label{equationMutationI}&\stackrel{\ref{itemdualmutationI}}{}& + B\sum_{i \in \tilde{\Gamma}(\overline{\eta})}\sum_{u \in K} b(\eta^{*}_{i},u )\li[H(\eta^{i;u}, \overline{\eta}^{\diamond})-H(\eta, \overline{\eta}^{\diamond})\re]\\
\label{equationResamplingJJ} &\stackrel{\ref{itemdualresamplingJJ}}{}& + \sum_{\gamma \not = \gamma' \in \Gamma(\overline{\eta})}(\tfrac{1}{2} + \tfrac{S}{2N}[\chi(\eta^{*}_{l'}) - \chi(\eta^{*}_{l})])H(\eta^{l' \to l}, \overline{\eta}^{\diamond}) -H(\eta, \overline{\eta}^{\diamond})\sum_{\gamma \not = \gamma' \in \Gamma(\overline{\eta})}\tfrac{1}{2}\\
\label{equationResamplingJI} &\stackrel{\ref{itemdualresamplingJI}}{}& + \sum_{i \in \tilde{\Gamma}(\overline{\eta}), \gamma \in \Gamma(\overline{\eta})}(\tfrac{1}{2} + \tfrac{S}{2N}[\chi(\eta^{*}_{i}) - \chi(\eta^{*}_{l})])H(\eta^{i \to  l}, \overline{\eta}^{\diamond})\\
\label{equationResamplingIJ} &\stackrel{\ref{itemdualresamplingIJ}}{}& + \sum_{\gamma \in \Gamma(\overline{\eta}), i \in \tilde{\Gamma}(\overline{\eta})} (\tfrac{1}{2} + \tfrac{S}{2N}[\chi(\eta^{*}_{l}) - \chi(\eta^{*}_{i})])H(\eta^{l\to i}, \overline{\eta}^{\diamond})\\
\label{equationResamplingJIandIJminus} &\stackrel{\ref{itemdualresamplingJI}, \ref{itemdualresamplingIJ}}{}& -H(\eta, \overline{\eta}^{\diamond})\sum_{ i \in \tilde{\Gamma}(\overline{\eta}), \gamma \in \Gamma(\overline{\eta})}\tfrac{1}{2} \;\; - H(\eta, \overline{\eta}^{\diamond})\sum_{\gamma \in \Gamma(\overline{\eta}), i \in \tilde{\Gamma}(\overline{\eta})}\tfrac{1}{2}\\\
\label{equationResamplingII} &\stackrel{\ref{itemdualresamplingII}}{}& + \sum_{ i\not = j \in \tilde{\Gamma}(\overline{\eta})}\mathbbm{1}\{\overline{\eta}^{\mbox{\tiny $I$}}_{i} \cap \overline{\eta}^{\mbox{\tiny $I$}}_{j}   \not\in \{\emptyset, K\} \}(\tfrac{1}{2} + \tfrac{S}{2N}[\chi(\eta^{*}_{i}) - \chi(\eta^{*}_{j})])H(\eta^{i\to j}, \overline{\eta}^{\diamond})\\
\label{equationResamplingIIminus} &\stackrel{\ref{itemdualresamplingII}}{}& -H(\eta, \overline{\eta}^{\diamond})\sum_{ i\not = j \in \tilde{\Gamma}(\overline{\eta})}\tfrac{1}{2}\mathbbm{1}\{\overline{\eta}^{\mbox{\tiny $I$}}_{i} \cap \overline{\eta}^{\mbox{\tiny $I$}}_{j}   \not= K \}\;.
\end{eqnarray}

In order to get to the right hand side of (\ref{equationSecondStepThm6}) we now write out this right hand side and use again that some terms cancel each other out. Then we explain step by step how the different parts relate to each other, where again $l = (\overline{\eta})^{\mbox{\tiny $J$}}_{\gamma,I}$ and $l' = (\overline{\eta})^{\mbox{\tiny $J$}}_{\gamma',I}$.

We have that
\begin{eqnarray}
&&  \sum_{\overline{\zeta} \in \overline{\mathcal{E}}} \overline{\mathcal{K}}(\overline{\eta}^{*}, \overline{\zeta})[H(\eta,\overline{\eta}^{\diamond:\overline{\zeta}}) - H(\eta,\overline{\eta}^{\diamond})] + V(\overline{\eta}^{*})
H(\eta,\overline{\eta}^{\diamond})\\
\label{equationRHSmutationJ} &=\stackrel{\ref{itemdualmutationJ}}{}& B\sum_{\gamma \in \Gamma(\overline{\eta})}\sum_{u \in K} b(u,\overline{\eta}^{\mbox{\tiny $J$}}_{\gamma,K})\li[H(\eta,\overline{\eta}^{\diamond:\overline{\eta}^{\gamma;u}})-H(\eta,\overline{\eta}^{\diamond})\re]\\
\label{equationVmutationJ}&\stackrel{V}{}&  + H(\eta,\overline{\eta}^{\diamond})B\sum_{\gamma \in \Gamma(\overline{\eta})}\li(\sum_{u \in K}b(u,\overline{\eta}^{\mbox{\tiny $J$}}_{\gamma,K})-1\re)\\
\label{equationMutationIcup}&\stackrel{\ref{itemdualmutationIcup}}{}& + B\sum_{i \in \tilde{\Gamma}(\overline{\eta})}\sum_{v \in \overline{\eta}^{\mbox{\tiny $I$}}_{i}}\sum_{u \not \in \overline{\eta}^{\mbox{\tiny $I$}}_{i}} b(u,v)\li[H(\eta,\overline{\eta}^{\diamond:\overline{\eta}^{i;\cup \{u\}}})-H(\eta,\overline{\eta}^{\diamond})\re]\\
\label{equationMutationIsetminus}&\stackrel{\ref{itemdualmutationIsetminus}}{}& + B\sum_{i \in \tilde{\Gamma}(\overline{\eta})}\mathbbm{1}\{|\overline{\eta}^{\mbox{\tiny $I$}}_{i}| > 1\}\sum_{v \in \overline{\eta}^{\mbox{\tiny $I$}}_{i}}\sum_{u \not \in \overline{\eta}^{\mbox{\tiny $I$}}_{i}} b(u,v)\li[H(\eta,\overline{\eta}^{\diamond:\overline{\eta}^{i;\setminus \{u\}}})-H(\eta,\overline{\eta}^{\diamond})\re]\\
\label{equationVmutationI}&\stackrel{V}{}& - H(\eta,\overline{\eta}^{\diamond})B\sum_{i \in \tilde{\Gamma}(\overline{\eta})}\mathbbm{1}\{|\overline{\eta}^{\mbox{\tiny $I$}}_{i}| = 1\}\sum_{u \in \overline{\eta}^{\mbox{\tiny $I$}}_{i}}\sum_{v \in K \setminus \{u\}}b(u,v)\\[0.2cm]
&& + \li(\mbox{terms corresponding to \ref{itemdualresamplingJJK} - \ref{itemdualresamplingIIw}}\re),
\end{eqnarray}
where
%%%%%%%%%%%%%%%%%%%%%%%%%%%%%%%%%%%%%%%%%%%%%%%%%%%%%%%%%%%%%%%%%%%%%%%%%%%%%%%%%%%%%%%%%%%%%%%%%%%      Proof of Theorem \ref{TheoremFeynmanKacDualityHMMandHBP}, Step 2

\begin{eqnarray}
&&  \li(\mbox{terms corresponding to \ref{itemdualresamplingJJK} - \ref{itemdualresamplingIIw}}\re)\\[0.2cm]
%-----------------------------------Resamling----------2a----------------
\label{equationResamplingJJK}&=\stackrel{\ref{itemdualresamplingJJK}}{}&  \sum_{\gamma \not = \gamma'\in \Gamma(\overline{\eta})}\mathbbm{1}\{\overline{\eta}^{\mbox{\tiny $J$}}_{\gamma,K} =\overline{\eta}^{\mbox{\tiny $J$}}_{\gamma',K}\}\li(\tfrac{1}{2} + \tfrac{S}{2N}[\chi(\overline{\eta}^{\mbox{\tiny $J$}}_{\gamma,K}) -1]\re) H(\eta, \overline{\eta}^{\diamond:\overline{\eta}^{\gamma \to \gamma'}})\\
\label{equationResamplingJJw}&\stackrel{\ref{itemdualresamplingJJw}}{}& + \sum_{\gamma \not = \gamma'\in \Gamma(\overline{\eta})}\mathbbm{1}\{\overline{\eta}^{\mbox{\tiny $J$}}_{\gamma,K} =\overline{\eta}^{\mbox{\tiny $J$}}_{\gamma',K}\}\sum_{w = 0}^{d-2}\tfrac{S}{2N}[\chi(w+1) -\chi(w)]H(\eta, \overline{\eta}^{\diamond:\overline{\eta}^{\gamma \stackrel{w}{\to} \gamma'}})\\ 
\label{equationVresamplingJJ} &\stackrel{\ref{itemdualresamplingJJ}, V}{}& - H(\eta,\overline{\eta}^{\diamond}) \sum_{\gamma \not = \gamma'\in \Gamma(\overline{\eta})} \tfrac{1}{2}\\
 %----------------------------------Resamling----------2b----------------
\label{equationResamplingJIK}&\stackrel{\ref{itemdualresamplingJIK}}{}& + \sum_{\gamma \in \Gamma(\overline{\eta}), i \in \tilde{\Gamma}(\overline{\eta})} \mathbbm{1}\{\overline{\eta}^{\mbox{\tiny $J$}}_{\gamma,K} \in \overline{\eta}^{\mbox{\tiny $I$}}_{i}\}\li(\tfrac{1}{2} + \tfrac{S}{2N}[\chi(\overline{\eta}^{\mbox{\tiny $J$}}_{\gamma,K}) -1]\re)H(\eta, \overline{\eta}^{\diamond:\overline{\eta}^{\gamma \to i}})\\
\label{equationResamplingJIw}&\stackrel{\ref{itemdualresamplingJIw}}{}& + \sum_{\gamma \in \Gamma(\overline{\eta}), i \in \tilde{\Gamma}(\overline{\eta})}\mathbbm{1}\{\overline{\eta}^{\mbox{\tiny $J$}}_{\gamma,K} \in \overline{\eta}^{\mbox{\tiny $I$}}_{i}\}\sum_{w = 0}^{d-2}\tfrac{S}{2N}[\chi(w+1) -\chi(w)]H(\eta, \overline{\eta}^{\diamond:\overline{\eta}^{\gamma \stackrel{w}{\to} i}})\\ 
 %----------------------------------Resamling----------2c----------------
\label{equationResamplingIJK}&\stackrel{\ref{itemdualresamplingIJK}}{}& + \sum_{i \in \tilde{\Gamma}(\overline{\eta}), \gamma \in \Gamma(\overline{\eta})} \mathbbm{1}\{\overline{\eta}^{\mbox{\tiny $J$}}_{\gamma,K} \in \overline{\eta}^{\mbox{\tiny $I$}}_{i}\}\li(\tfrac{1}{2} + \tfrac{S}{2N}[\chi(\overline{\eta}^{\mbox{\tiny $J$}}_{\gamma,K}) -1]\re)H(\eta, \overline{\eta}^{\diamond:\overline{\eta}^{i \to \gamma }})\\
\label{equationResamplingIJw} &\stackrel{\ref{itemdualresamplingIJw}}{}& + \sum_{i \in \tilde{\Gamma}(\overline{\eta}), \gamma \in \Gamma(\overline{\eta})}\mathbbm{1}\{\overline{\eta}^{\mbox{\tiny $J$}}_{\gamma,K} \in \overline{\eta}^{\mbox{\tiny $I$}}_{i}\}\sum_{w = 0}^{d-2}\tfrac{S}{2N}[\chi(w+1) -\chi(w)]H(\eta, \overline{\eta}^{\diamond:\overline{\eta}^{i \stackrel{w}{\to} \gamma}})\\ 
\label{equationVresamplingJIandIJ}&\stackrel{\ref{itemdualresamplingJI},\ref{itemdualresamplingIJ}, V}{}& - H(\eta,\overline{\eta}^{\diamond}) \sum_{\gamma \in \Gamma(\overline{\eta}), i \in \tilde{\Gamma}(\overline{\eta})}\tfrac{1}{2} -H(\eta,\overline{\eta}^{\diamond}) \sum_{i \in \tilde{\Gamma}(\overline{\eta}), \gamma \in \Gamma(\overline{\eta})}\tfrac{1}{2}\\
 %----------------------------------Resamling----------2d----------------
\label{equationResamplingIIK}&\stackrel{\ref{itemdualresamplingIIK}}{}& + \sum_{i \not= j\in \tilde{\Gamma}(\overline{\eta})}\mathbbm{1}\{\overline{\eta}^{\mbox{\tiny $I$}}_{i} \cap \overline{\eta}^{\mbox{\tiny $I$}}_{j}   \not\in \{\emptyset, K\}\} (\tfrac{1}{2}  + \tfrac{S}{2N}[\chi(\min \overline{\eta}^{\mbox{\tiny $I$}}_{i} \cap \overline{\eta}^{\mbox{\tiny $I$}}_{j}) - 1])H(\eta, \overline{\eta}^{\diamond:\overline{\eta}^{i \cap j}})\\
\label{equationResamplingIIw}&\stackrel{\ref{itemdualresamplingIIw}}{}& + \sum_{i\not=j\in \tilde{\Gamma}(\overline{\eta})}\mathbbm{1}\{\overline{\eta}^{\mbox{\tiny $I$}}_{i} \cap \overline{\eta}^{\mbox{\tiny $I$}}_{j}   \not\in \{\emptyset, K\} \}\sum_{w = 0}^{d-2}\tfrac{S}{2N}[\chi(w+1) -\chi(w)]H(\eta, \overline{\eta}^{\diamond:\overline{\eta}^{i \stackrel{w}{\cap} j}})\\ 
\label{equationResamplingIIv}&\stackrel{\ref{itemdualresamplingIIv}}{}& + \sum_{i\not=j \in \tilde{\Gamma}(\overline{\eta})}\mathbbm{1}\{\overline{\eta}^{\mbox{\tiny $I$}}_{i} \cap \overline{\eta}^{\mbox{\tiny $I$}}_{j}   \not\in \{\emptyset, K\} \}\sum_{\stackrel{v \in \overline{\eta}^{\mbox{\tiny $I$}}_{i} \cap \overline{\eta}^{\mbox{\tiny $I$}}_{j}:}{v \not= \min \overline{\eta}^{\mbox{\tiny $I$}}_{i} \cap \overline{\eta}^{\mbox{\tiny $I$}}_{j}}}\tfrac{S}{2N}[\chi(v) -\chi(v^{<})]H(\eta, \overline{\eta}^{\diamond:\overline{\eta}^{i \stackrel{v}{\cap} j}})\\ 
\label{equationVresamplingII}&\stackrel{\ref{itemdualresamplingII}, V}{}& - H(\eta,\overline{\eta}^{\diamond}) \sum_{i\not=j \in \tilde{\Gamma}(\overline{\eta})} \tfrac{1}{2}\mathbbm{1}\{\overline{\eta}^{\mbox{\tiny $I$}}_{i} \cap \overline{\eta}^{\mbox{\tiny $I$}}_{j}   \not= K \}\; .
\end{eqnarray}
%%%%%%%%%%%%%%%%%%%%%%%%%%%%%%%%%%%%%%%%%%%%%%%%%%%%%%%%%%%%%%%%%%%%%%%%%%%%%%%%%%%%%%%%%%%%%%%%%%%      Proof of Theorem \ref{TheoremFeynmanKacDualityHMMandHBP}, Step 2

The equation $(\ref{equationLHSmutationJ}) = (\ref{equationRHSmutationJ}) + (\ref{equationVmutationJ})$ holds since for each $\gamma \in \Gamma(\overline{\eta})$,
\begin{eqnarray}
&& \sum_{u \in K} b(\eta^{*}_{l},u )[\mathbbm{1}\{u =\overline{\eta}^{\mbox{\tiny $J$}}_{\gamma,K}\}-\mathbbm{1}\{\eta^{*}_{l} =\overline{\eta}^{\mbox{\tiny $J$}}_{\gamma,K}\}]  + \mathbbm{1}\{\eta^{*}_{l} =\overline{\eta}^{\mbox{\tiny $J$}}_{\gamma,K}\}\\
&=&  \sum_{u \in K} b(u,\overline{\eta}^{\mbox{\tiny $J$}}_{\gamma,K})[\mathbbm{1}\{\eta^{*}_{l} = u \}-\mathbbm{1}\{\eta^{*}_{l} =\overline{\eta}^{\mbox{\tiny $J$}}_{\gamma,K}\}] + \mathbbm{1}\{\eta^{*}_{l} =\overline{\eta}^{\mbox{\tiny $J$}}_{\gamma,K}\}\sum_{u \in K}b(u,\overline{\eta}^{\mbox{\tiny $J$}}_{\gamma,K})\; .
\end{eqnarray}
%-------------------------------------------------------------------------
The equation $(\ref{equationMutationI}) = (\ref{equationMutationIcup}) + (\ref{equationMutationIsetminus}) + (\ref{equationVmutationI})$ holds since for each $i \in \tilde{\Gamma}(\overline{\eta})$,
\begin{eqnarray}
&& \sum_{v \in K} b(\eta^{*}_{i},v)[\mathbbm{1}\{v \in \overline{\eta}^{\mbox{\tiny $I$}}_{i}\}- \mathbbm{1}\{\eta^{*}_{i} \in \overline{\eta}^{\mbox{\tiny $I$}}_{i}\}]\\
&=& \mathbbm{1}\{\eta^{*}_{i} \not\in \overline{\eta}^{\mbox{\tiny $I$}}_{i}\}\sum_{v \in \overline{\eta}^{\mbox{\tiny $I$}}_{i}} b(\eta^{*}_{i},v) - \mathbbm{1}\{\eta^{*}_{i} \in \overline{\eta}^{\mbox{\tiny $I$}}_{i}\}\sum_{v \not\in \overline{\eta}^{\mbox{\tiny $I$}}_{i}} b(\eta^{*}_{i},v)\\
&=& \sum_{v \in \overline{\eta}^{\mbox{\tiny $I$}}_{i}}\sum_{u \not \in \overline{\eta}^{\mbox{\tiny $I$}}_{i}} b(u,v)[\mathbbm{1}\{\eta^{*}_{i} \in \overline{\eta}^{\mbox{\tiny $I$}}_{i}\cup \{u\}\}- \mathbbm{1}\{\eta^{*}_{i} \in \overline{\eta}^{\mbox{\tiny $I$}}_{i}\}]\\
&&+\mathbbm{1}\{|\overline{\eta}^{\mbox{\tiny $I$}}_{i}| > 1\}\sum_{v \not\in \overline{\eta}^{\mbox{\tiny $I$}}_{i}}\sum_{u \in \overline{\eta}^{\mbox{\tiny $I$}}_{i}} b(u,v)[\mathbbm{1}\{\eta^{*}_{i} \in \overline{\eta}^{\mbox{\tiny $I$}}_{i}\setminus \{u\}\}- \mathbbm{1}\{\eta^{*}_{i} \in \overline{\eta}^{\mbox{\tiny $I$}}_{i}\}]\\
&&-\sum_{i \in \tilde{\Gamma}(\overline{\eta})}\mathbbm{1}\{|\overline{\eta}^{\mbox{\tiny $I$}}_{i}| = 1\}\sum_{u \in \overline{\eta}^{\mbox{\tiny $I$}}_{i}}\sum_{v \in K \setminus \{u\}}b(u,v)\;.
\end{eqnarray}
%-------------------------------------------------------------------------
The equation $(\ref{equationResamplingJJ}) = (\ref{equationResamplingJJK}) + (\ref{equationResamplingJJw}) + (\ref{equationVresamplingJJ})$ holds since for each $\gamma, \gamma' \in \Gamma(\overline{\eta})$ with $\gamma \not = \gamma'$,
\begin{eqnarray}
&& \li(\tfrac{1}{2} + \tfrac{S[\chi(\eta^{*}_{l'}) - \chi(\eta^{*}_{l})]}{2N}\re)\mathbbm{1}\{\eta^{*}_{l'} =\overline{\eta}^{\mbox{\tiny $J$}}_{\gamma,K}\}H^{\mbox{\tiny $\mathcal{D}$}}_{\eta, \overline{\eta}^{\diamond}}(l',\gamma) \mathbbm{1}\{\eta^{*}_{l'} =\overline{\eta}^{\mbox{\tiny $J$}}_{\gamma',K}\}H^{\mbox{\tiny $\mathcal{D}$}}_{\eta, \overline{\eta}^{\diamond}}(l',\gamma') \\
&=& \mathbbm{1}\{\overline{\eta}^{\mbox{\tiny $J$}}_{\gamma,K} =\overline{\eta}^{\mbox{\tiny $J$}}_{\gamma',K}\}\li(\tfrac{1}{2} + \tfrac{S[\chi(\overline{\eta}^{\mbox{\tiny $J$}}_{\gamma,K}) -1]}{2N}\re) \mathbbm{1}\{\eta^{*}_{l'} =\overline{\eta}^{\mbox{\tiny $J$}}_{\gamma,K}\}H^{\mbox{\tiny $\mathcal{D}$}}_{\eta, \overline{\eta}^{\diamond}}(l',\gamma\cup\gamma') \mathbbm{1}\{\eta^{*}_{l} \in K\}\\
 &&+ \mathbbm{1}\{\overline{\eta}^{\mbox{\tiny $J$}}_{\gamma,K} =\overline{\eta}^{\mbox{\tiny $J$}}_{\gamma',K}\}\sum_{w = 0}^{d-2}\tfrac{S[\chi(w+1) -\chi(w)]}{2N}\mathbbm{1}\{\eta^{*}_{l'} =\overline{\eta}^{\mbox{\tiny $J$}}_{\gamma,K}\}H^{\mbox{\tiny $\mathcal{D}$}}_{\eta, \overline{\eta}^{\diamond}}(l',\gamma\cup\gamma') \mathbbm{1}\{\eta^{*}_{l} \le w\}.\;\;
\end{eqnarray}
%------------------------------------------------------------------------
The equation $(\ref{equationResamplingJI}) + (\ref{equationResamplingIJ}) + (\ref{equationResamplingJIandIJminus}) = (\ref{equationResamplingJIK}) + (\ref{equationResamplingJIw}) + (\ref{equationResamplingIJK}) + (\ref{equationResamplingIJw}) + (\ref{equationVresamplingJIandIJ})$ holds since for each $\gamma\in \Gamma(\overline{\eta})$ and each $i \in \tilde{\Gamma}(\overline{\eta})$,
\begin{eqnarray}
&& \li(\tfrac{1}{2} + \tfrac{S}{2N}[\chi(\eta^{*}_{i}) - \chi(\eta^{*}_{l})]\re)\mathbbm{1}\{\eta^{*}_{i} =\overline{\eta}^{\mbox{\tiny $J$}}_{\gamma,K}\}H^{\mbox{\tiny $\mathcal{D}$}}_{\eta, \overline{\eta}^{\diamond}}(i,\gamma) \mathbbm{1}\{\eta^{*}_{i} \in \overline{\eta}^{\mbox{\tiny $I$}}_{i}\} \\
&=& \mathbbm{1}\{\overline{\eta}^{\mbox{\tiny $J$}}_{\gamma,K} \in \overline{\eta}^{\mbox{\tiny $I$}}_{i}\}\li(\tfrac{1}{2} + \tfrac{S}{2N}[\chi(\overline{\eta}^{\mbox{\tiny $J$}}_{\gamma,K}) -1]\re) \mathbbm{1}\{\eta^{*}_{i} =\overline{\eta}^{\mbox{\tiny $J$}}_{\gamma,K}\}H^{\mbox{\tiny $\mathcal{D}$}}_{\eta, \overline{\eta}^{\diamond}}(i,\gamma)\mathbbm{1}\{\eta^{*}_{l} \in K\}\\
 &&+ \mathbbm{1}\{\overline{\eta}^{\mbox{\tiny $J$}}_{\gamma,K} \in \overline{\eta}^{\mbox{\tiny $I$}}_{i}\}\sum_{w = 0}^{d-2}\tfrac{S}{2N}[\chi(w+1) -\chi(w)]\mathbbm{1}\{\eta^{*}_{i} =\overline{\eta}^{\mbox{\tiny $J$}}_{\gamma,K}\}H^{\mbox{\tiny $\mathcal{D}$}}_{\eta, \overline{\eta}^{\diamond}}(i,\gamma)\mathbbm{1}\{\eta^{*}_{l} \le w\}
\end{eqnarray}
and 
\begin{eqnarray}
&& \li(\tfrac{1}{2} + \tfrac{S}{2N}[\chi(\eta^{*}_{l}) - \chi(\eta^{*}_{i})]\re)\mathbbm{1}\{\eta^{*}_{l} =\overline{\eta}^{\mbox{\tiny $J$}}_{\gamma,K}\}\mathbbm{1}\{\eta^{*}_{l} \in \overline{\eta}^{\mbox{\tiny $I$}}_{i}\} \\
&=& \mathbbm{1}\{\overline{\eta}^{\mbox{\tiny $J$}}_{\gamma,K} \in \overline{\eta}^{\mbox{\tiny $I$}}_{i}\}\li(\tfrac{1}{2} + \tfrac{S}{2N}[\chi(\overline{\eta}^{\mbox{\tiny $J$}}_{\gamma,K}) -1]\re)\mathbbm{1}\{\eta^{*}_{l} =\overline{\eta}^{\mbox{\tiny $J$}}_{\gamma,K}\}\mathbbm{1}\{\eta^{*}_{i} \in K\}\\
 &&+ \mathbbm{1}\{\overline{\eta}^{\mbox{\tiny $J$}}_{\gamma,K} \in \overline{\eta}^{\mbox{\tiny $I$}}_{i}\}\sum_{w = 0}^{d-2}\tfrac{S}{2N}[\chi(w+1) -\chi(w)]\mathbbm{1}\{\eta^{*}_{l} =\overline{\eta}^{\mbox{\tiny $J$}}_{\gamma,K}\}\mathbbm{1}\{\eta^{*}_{i} \le w\}\;.
\end{eqnarray}
%-------------------------------------------------------------------------
The equation $(\ref{equationResamplingII}) + (\ref{equationResamplingIIminus}) = (\ref{equationResamplingIIK}) + (\ref{equationResamplingIIw}) + (\ref{equationResamplingIIv}) + (\ref{equationVresamplingII})$ holds  since for each $i,j \in \tilde{\Gamma}(\overline{\eta})$ with $\overline{\eta}^{\mbox{\tiny $I$}}_{i} \cap \overline{\eta}^{\mbox{\tiny $I$}}_{j}   \not\in \{\emptyset, K\}$,
\begin{eqnarray}
&& \li(\tfrac{1}{2} + \tfrac{S}{2N}[\chi(\eta^{*}_{i}) - \chi(\eta^{*}_{j})]\re)\mathbbm{1}\{\eta^{*}_{i} \in \overline{\eta}^{\mbox{\tiny $I$}}_{i}\}\mathbbm{1}\{\eta^{*}_{i} \in \overline{\eta}^{\mbox{\tiny $I$}}_{j}\} \\
&=& \li(\tfrac{1}{2} + \tfrac{S}{2N}[\chi(\min \overline{\eta}^{\mbox{\tiny $I$}}_{i} \cap \overline{\eta}^{\mbox{\tiny $I$}}_{j}) -1]\re) \mathbbm{1}\{\eta^{*}_{i} \in \overline{\eta}^{\mbox{\tiny $I$}}_{i} \cap \overline{\eta}^{\mbox{\tiny $I$}}_{j}\}\mathbbm{1}\{\eta^{*}_{j} \in K\}\\
 &&+ \sum_{w = 0}^{d-2}\tfrac{S}{2N}[\chi(w+1) -\chi(w)]\mathbbm{1}\{\eta^{*}_{i} \in\overline{\eta}^{\mbox{\tiny $I$}}_{i} \cap \overline{\eta}^{\mbox{\tiny $I$}}_{j}\}\mathbbm{1}\{\eta^{*}_{j} \le w\}\\
 &&+ \sum_{v \in \overline{\eta}^{\mbox{\tiny $I$}}_{i} \cap \overline{\eta}^{\mbox{\tiny $I$}}_{j}:v \not= \min \overline{\eta}^{\mbox{\tiny $I$}}_{i} \cap \overline{\eta}^{\mbox{\tiny $I$}}_{j}}\tfrac{S}{2N}[\chi(v) -\chi(v^{<})]\mathbbm{1}\{\eta^{*}_{i} \in \{v, \dots\}\}\mathbbm{1}\{\eta^{*}_{j} \in K\}\,.
\end{eqnarray}
%%%%%%%%%%%%%%%%%%%%%%%%%%%%%%%%%%%%%%%%%%%%%%%%%%%%%%%%%%%%%%%%%%%%%%%%%%%%%%%%%%%%%%%%%%%%%%%%%%%      Proof of Theorem \ref{TheoremAnalyticalHMM}

\subsection{Proof of Theorem \ref{TheoremAnalyticalHMM}}
\noindent We show part \ref{TheoremAnalyticalHMMa}) and obtain part \ref{TheoremAnalyticalHMMb}) and \ref{TheoremAnalyticalHMMc}) by the Theorems 4.4.2 and 4.4.6 in \cite{EK86}. The proof of part \ref{TheoremAnalyticalHMMa}) includes existence and uniqueness of the martingale problem, where the proof for existence also implies part \ref{TheoremAnalyticalHMMd}).\bigskip

{\it Existence:} For each distribution $\mu$ on $\mathcal{E}$ let $\mathbb{Q}^{\mbox{\tiny \sc Time}}_{\mu}(\cdot) := \int \mathbb{Q}^{\mbox{\tiny \sc Time}}_{\eta}(\cdot) \mu(d\eta)$  with
\begin{equation}
\mathbb{Q}^{\mbox{\tiny \sc Time}}_{\eta}\li(X_{t} = (\eta^{\mbox{\tiny \sc Time}} + t,\eta^{\mbox{\tiny $\mathcal{D}$}}) \; \mbox{ for all } t \ge 0\re) = 1 
\end{equation}
for all $\eta \in \mathcal{E}$. Since for every distribution $\mu$ on $\mathcal{E}$, $\mathbb{Q}^{\mbox{\tiny \sc Time}}_{\mu}$ is a solution of the $\Omega$-martingale problem for $(L^{\mbox{\tiny \sc Time}},\mu)$ w.r.t$.$ $\mathcal{A}$, the  statement follows from Proposition 4.10.2 in \cite{EK86}. Observe that the construction for the solution of the martingale problem in the proof of Proposition 4.10.2 fits with our description in Definition \ref{DefinitionHMM}, that is, if $\mu$ is defined as in (\ref{equationMu}), then the law of the piecewise deterministic Markov jump process from Definition \ref{DefinitionHMM} is a solution of the $\Omega$-martingale problem for $(L,\mu)$ w.r.t$.$ $\mathcal{A}$. \bigskip

{\it Uniqueness:} Since for $J = I$ the family $\{ \eta \mapsto H(\eta,\overline{\eta}^{\diamond}): \overline{\eta}^{\diamond} \in \overline{\mathcal{E}}^{\diamond}, H \in \mathcal{H}\}$ is measure determining, uniqueness follows from Theorem \ref{TheoremFeynmanKacDualityHMMandHBP} and Theorems 4.4.2 in \cite{EK86}.
%%%%%%%%%%%%%%%%%%%%%%%%%%%%%%%%%%%%%%%%%%%%%%%%%%%%%%%%%%%%%%%%%%%%%%%%%%%%%%%%%%%%%%%%%%%%%%%%%%%      Proof of Theorem \ref{TheoremStochasticRepresentation}

\subsection{Proofs of Proposition \ref{PropositionDualityTypInfoHMM} and Theorem \ref{TheoremStochasticRepresentation}}
\noindent In both cases we can apply Theorem \ref{TheoremFeynmanKacDualityHMMandHBP}. To see Proposition \ref{PropositionDualityTypInfoHMM} we note that $H = H^{*}$ if $g \equiv 1$ and $F_{j}^{n} \equiv 1$ for all $j \in J$.

To see Theorem \ref{TheoremStochasticRepresentation} let $\overline{\eta} \in \overline{\mathcal{E}}$ and  $\mu$ be as in (\ref{equationMu}) with $c = -T$ and general $\mu^{*}$. Now, if the element $\overline{\eta}^{\diamond} \in \overline{\mathcal{E}}^{\diamond}$ satisfies $(\overline{\eta}^{\diamond})^{\mbox{\tiny \sc Time}} = 0$ and $\overline{\eta}^{*} = \overline{\xi^{*}}$ with $\xi^{*} \in K^{J}$, $g \equiv 1$, $-T \le r_{j}^{1} < t_{j}^{1} < \dots < r_{j}^{m_{j}} < t_{j}^{m_{j}} \le 0$ and $F_{j}^{n} \in C_{b}(K \times I)$, then
\begin{equation}
\mathbb{E}_{\mu}\li[H(X_{T}, \overline{\eta}^{\diamond})\re] =(\ref{equationFunctionalHMM})
\end{equation} 
which is the left hand side of the equation in Theorem \ref{TheoremStochasticRepresentation} and 
\begin{eqnarray}
&& \overline{\mathbb{E}}^{\diamond}_{\overline{\eta}^{\diamond}}\li[ \mathbb{E}_{\mu}[H(X_{0}, \overline{X}^{\diamond}_{T})]e^{\int_{0}^{T} V^{\diamond}\li(\overline{X}^{\diamond}_{s}\re)ds}\re]\\
&=& \overline{\mathbb{E}}^{\diamond}_{\overline{\eta}^{\diamond}} \li[\li(\prod_{j \in J}\prod_{n = 1}^{m_{j}} \int_{r_{j}^{n}}^{t_{j}^{n}}F_{j}^{n}((( \overline{X}^{\diamond}_{T})^{\mbox{\tiny $\mathcal{D}$}}_{-s})^{\mbox{\tiny $J$}}_{j})ds\re) \mathbb{E}_{\mu}[H^{*}(X_{0}, \overline{X}^{*}_{T})]e^{\int_{0}^{T}V(\overline{X}^{*}_{s})ds} \re]\\
&=& \overline{\mathbb{E}}_{\overline{\xi^{*}}}\li[\li(\prod_{j \in J}\prod_{n = 1}^{m_{j}}\int_{r_{j}^{n}}^{t_{j}^{n}} F_{j}^{n}((\overline{X}_{-s})^{\mbox{\tiny $J$}}_{j})ds\re) \mathbb{E}_{\mu}[H^{*}(X_{0}, \overline{X}_{T})]e^{\int_{0}^{T}V(\overline{X}_{s})ds} \re]
\end{eqnarray}
which is the right hand side of the equation in Theorem \ref{TheoremStochasticRepresentation}.
%%%%%%%%%%%%%%%%%%%%%%%%%%%%%%%%%%%%%%%%%%%%%%%%%%%%%%%%%%%%%%%%%%%%%%%%%%%%%%%%%%%%%%%%%%%%%%%%%%%      Proof of Theorem \ref{TheoremTransformedBP}

\subsection{Proof of Theorem \ref{TheoremTransformedBP}}
\noindent Although it is also a consequence of the time-inhomogeneous case (part \ref{TheoremTransformedBPitemA})) we briefly explain why the statements hold in the time-homogeneous case (part \ref{TheoremTransformedBPitemB})). If $h:\overline{\mathcal{E}} \to  (0,\infty)$ and $h^{T}(t,\cdot) =h(\cdot)$ for all $0 \le  t < T < \infty$, then
\begin{equation}
\li(\mathbb{E}_{\mu}\li[H^{*}(X_{0}, \overline{X}_{t})\re]e^{\int_{0}^{t}V(\overline{X}_{s})ds}\re)_{t \ge 0}
\end{equation}
is a martingale under each $\overline{\mathbb{P}}_{\overline{\eta}}$. This means (see \cite{FS04}) the time-homogeneous transformed BP is a compensated $h$-transform corresponding to  a strongly continuous semigroup with bounded generator
\begin{equation}
\overline{L}^{h}f(\overline{\eta}) =  \tfrac{1}{h(\overline{\eta})}\li(\overline{L}[fh](\overline{\eta}) - f(\overline{\eta})\overline{L}h(\overline{\eta}) \re)
\end{equation}
which gives the right hand side of (\ref{equationLh}). 
%%%%%%%%%%%%%%%%%%%%%%%%%%%%%%%%%%%%%%%%%%%%%%%%%%%%%%%%%%%%%%%%%%%%%%%%%%%%%%%%%%%%%%%%%%%%%%%%%%%      Proof of Theorem \ref{TheoremTransformedBP}, time-inhomogeneous

In the time-inhomogeneous case we use the abbreviation $g(\overline{\eta}) := \mathbb{E}_{\mu}[H^{*}(X_{0}, \overline{\eta})]$ and consider instead of a semigroup the family 
\begin{equation}
\{\overline{S}^{h^{T}}_{t,s}: 0 \le t \le s < T\} \; \mbox{ defined by }\; \overline{S}^{h^{T}}_{t,s}:C_{b}(\overline{\mathcal{E}}) \to C_{b}(\overline{\mathcal{E}}), \; f(\overline{\eta}) \mapsto  \overline{\mathbb{E}}^{h^{T}}_{t,\overline{\eta}}\li[f(\overline{X}_{s})\re]\;.
\end{equation}
This family is a {\it strongly continuous propagator}, see \cite{Kol10},  because it satisfies the Chapman-Kolmogorov equation, i.e$.$
\begin{eqnarray}
\overline{S}^{h^{T}}_{t,s}f(\overline{\eta}) &=& \tfrac{1}{h^{T}(t,\overline{\eta})}\overline{\mathbb{E}}_{\overline{\eta}}\li[f(\overline{X}_{s-t})g(\overline{X}_{T-t})e^{\int_{0}^{T-t}V(\overline{X}_{\rho})d\rho}\re]\\
&=&\tfrac{1}{h^{T}(t,\overline{\eta})}\overline{\mathbb{E}}_{\overline{\eta}}\li[ \overline{\mathbb{E}}_{\overline{X}_{r-t}}\li[f(\overline{X}_{s-r}) g(\overline{X}_{T-r})e^{\int_{0}^{T-r}V(\overline{X}_{\rho})d\rho}\re] e^{\int_{0}^{r-t}V(\overline{X}_{\rho})d\rho} \re]\\
&=&\tfrac{1}{h^{T}(t,\overline{\eta})}\overline{\mathbb{E}}_{\overline{\eta}}\li[ \overline{S}^{h^{T}}_{r,s}f(\overline{X}_{r-t}) h^{T}(r,\overline{X}_{r-t}) e^{\int_{0}^{r-t}V(\overline{X}_{\rho})d\rho} \re]\\
&=&\tfrac{1}{h^{T}(t,\overline{\eta})}\overline{\mathbb{E}}_{\overline{\eta}}\li[ \overline{S}^{h^{T}}_{r,s}f(\overline{X}_{r-t}) g(\overline{X}_{T-t})e^{\int_{0}^{T-t}V(\overline{X}_{\rho})d\rho}\re]\\
&=&  \overline{S}^{h^{T}}_{t,r}\overline{S}^{h^{T}}_{r,s}f(\overline{\eta})\;,
\end{eqnarray}
and is strongly continuous in $(t,s)$, where the latter property is due to the fact that $\overline{X}$ is a Borel strong Markov process with finite state space and bounded generator $\overline{L}$.  Hence the time-inhomogeneous transformed BP is a time-inhomogeneous Borel strong Markov process. (Note that the basic properties of a time-inhomogeneous Borel strong Markov process are explicitly stated in Theorem 9 of \cite{Seidel15}, but are not used in the present paper)
%%%%%%%%%%%%%%%%%%%%%%%%%%%%%%%%%%%%%%%%%%%%%%%%%%%%%%%%%%%%%%%%%%%%%%%%%%%%%%%%%%%%%%%%%%%%%%%%%%%      Proof of Theorem \ref{TheoremTransformedBP}, time-inhomogeneous

What we need (for understanding the evolution of the time-inhomogeneous transformed BP and proving  Theorem \ref{TheoremLongtime}) is the form of the family of generators in (\ref{equationFamilyLhT}) corresponding to the strongly continuous propagator. In order to obtain (\ref{equationLhT}) we show that
\begin{equation}\label{equationLhTcalculation}
\overline{L}^{h^{T}}_{t}f(\overline{\eta}) := \lim_{s \downarrow t}\tfrac{\overline{S}^{h^{T}}_{t,s}f(\overline{\eta})- f(\overline{\eta})}{s-t} =  \tfrac{1}{h^{T}(t,\overline{\eta})}\li(\overline{L}[f(\cdot)h^{T}(t, \cdot)](\overline{\eta}) - f(\overline{\eta})\overline{L}[h^{T}(t, \cdot)](\overline{\eta}) \re)
\end{equation}
which finishes the proof.

Recall that $g(\overline{\eta}) = \mathbb{E}_{\mu}[H^{*}(X_{0}, \overline{\eta})]$ and  let $\{\overline{S}^{V}_{r}: r \ge 0\}$  be the strongly continuous semigroup on $C_{b}(\overline{\mathcal{E}})$  given by 
\begin{equation}
\overline{S}^{V}_{r}f(\overline{\eta}) :=\exp(r\overline{L}^{V}f(\overline{\eta})) :=\exp(r\overline{L}f(\overline{\eta}) + rV(\overline{\eta})f(\overline{\eta}))\;. 
\end{equation}
Note that
\begin{equation}
h^{T}(t, \overline{\eta}) = \overline{S}^{V}_{T-t}g(\overline{\eta})\;\;\mbox{ and }\; \overline{S}^{h^{T}}_{t,s}f(\overline{\eta}) = \frac{\overline{S}^{V}_{s-t}[f\overline{S}^{V}_{T-s}g](\overline{\eta})}{\overline{S}^{V}_{T-t}g(\overline{\eta})} \;. 
\end{equation}
Now (\ref{equationLhTcalculation}) follows due to
\begin{eqnarray}
&&\overline{S}^{V}_{T-t}g(\overline{\eta})\overline{L}^{h^{T}}_{t}f(\overline{\eta})\\ 
&=& \lim_{s \downarrow t}\tfrac{1}{s-t}\li(\overline{S}^{V}_{s-t}[f\overline{S}^{V}_{T-s}g](\overline{\eta})- f(\overline{\eta})\overline{S}^{V}_{T-t}g(\overline{\eta})\re)\\
&=& \lim_{s \downarrow t}\tfrac{f(\overline{\eta})}{s-t}\li(\overline{S}^{V}_{T-s}g(\overline{\eta})- \overline{S}^{V}_{T-t}g(\overline{\eta})\re) + \overline{L}^{V}[f\overline{S}^{V}_{T-t}g](\overline{\eta})\\
&=& f(\overline{\eta})[-\overline{L}[\overline{S}^{V}_{T-t}g](\overline{\eta}) - V(\overline{\eta})\overline{S}^{V}_{T-t}(\overline{\eta})] +  \overline{L}[f\overline{S}^{V}_{T-t}g](\overline{\eta}) + V(\overline{\eta})f(\overline{\eta})\overline{S}^{V}_{T-t}(\overline{\eta})\\
&=& \overline{L}[f\overline{S}^{V}_{T-t}g](\overline{\eta}) -f(\overline{\eta})\overline{L}[\overline{S}^{V}_{T-t}g](\overline{\eta})\;.
\end{eqnarray}
%%%%%%%%%%%%%%%%%%%%%%%%%%%%%%%%%%%%%%%%%%%%%%%%%%%%%%%%%%%%%%%%%%%%%%%%%%%%%%%%%%%%%%%%%%%%%%%%%%%      Proof of Theorem \ref{TheoremStrongStochasticRepresentation}

\subsection{Proof of Theorem \ref{TheoremStrongStochasticRepresentation}}
\noindent  Since
\begin{eqnarray}
&&\frac{\mathbb{E}_{\mu}\li[\li(\prod_{j \in J}\prod_{n = 1}^{m_{j}}\int_{r_{j}^{n}}^{t_{j}^{n}}F_{j}^{n}((X_{T})^{\mbox{\tiny $\mathcal{D}$}}_{j,s})ds\re) \mathbbm{1}\{((X_{T})^{*}_{j})_{j \in J}  = \xi^{*}\}\re]}{\mathbb{P}_{\mu}(((X_{T})^{*}_{j})_{j \in J} = \xi^{*})} \\
&=& \overline{\mathbb{E}}^{h^{T}}_{\overline{\xi^{*}}}\li[\prod_{j \in J}\prod_{n = 1}^{m_{j}}\int^{t_{j}^{n}}_{r_{j}^{n}} F_{j}^{n}((\overline{X}_{-s})^{\mbox{\tiny $J$}}_{j})ds\re]
\end{eqnarray}
for all $-T \le r_{j}^{1} < t_{j}^{1} < \dots < r_{j}^{m_{j}} < t_{j}^{m_{j}} \le 0$ and all $F_{j}^{n} \in C_{b}(K \times I)$ due to Theorem \ref{TheoremStochasticRepresentation}, the statements in Theorem \ref{TheoremStrongStochasticRepresentation} directly follow from the definition of $\overline{\mathbb{P}}^{h^{T}}_{\overline{\xi^{*}}}$, the special map $\mathbb{F}_{T}$ defined in (\ref{equationFT}) and the fact that the collection  of all these functionals is measure determining.
%%%%%%%%%%%%%%%%%%%%%%%%%%%%%%%%%%%%%%%%%%%%%%%%%%%%%%%%%%%%%%%%%%%%%%%%%%%%%%%%%%%%%%%%%%%%%%%%%%%      Proof of Theorem \ref{TheoremLongtime}

\subsection{Proof of Theorem \ref{TheoremLongtime}}
\noindent Let $\xi^{*} \in K^{J}$ and $\overline{\xi^{*}}$ be the initial state defined in (\ref{equationoverlinexi}). The statements in \ref{TheoremLongtimeItemA}) and \ref{TheoremLongtimeItemB}) hold if we have shown that
\begin{equation}\label{equationPhTtoPh}
\lim_{T \to \infty}\overline{\mathbb{P}}_{0,\overline{\xi^{*}}}^{h^{T}}\li( \mathbb{F}(\li((\overline{X}_{t\wedge T})^{\mbox{\tiny $J$}}\re)_{t \ge 0}) \in \, \cdot \,\re) = \overline{\mathbb{P}}_{\overline{\xi^{*}}}^{h}\li( \mathbb{F}(\li((\overline{X}_{t})^{\mbox{\tiny $J$}}\re)_{t \ge 0}) \in \, \cdot \,\re) \;,
\end{equation}
where we shall verify that in  the limit  $T \to \infty$ the generators $\overline{L}^{h^{T}}_{t}$ (corresponding to the strongly continuous propagator defined in the proof of Theorem \ref{TheoremTransformedBP}) converge to the bounded generator $\overline{L}^{h}$ uniformly for bounded $t$.

Since $h>0$ and since for each $t \ge 0$,
\begin{equation}
\lim_{T \to \infty}\sup_{\overline{\eta} \in \overline{\mathcal{E}}}\li|h^{T}(t,\overline{\eta}) - h(\overline{\eta})\re| = \lim_{T \to \infty}\sup_{\overline{\eta} \in \overline{\mathcal{E}}}\li|\mathbb{E}_{\mu}\li[H^{*}(X_{T-t},\overline{\eta})\re] - \mathbb{E}\li[H^{*}(X_{\infty},\overline{\eta})\re]\re| = 0\;, 
\end{equation}
we obtain that 
\begin{equation}
\lim_{T \to \infty}\sup_{\overline{\eta} \in \overline{\mathcal{E}}}\li|\overline{L}^{h^{T}}_{t}(\overline{\eta}) - \overline{L}^{h}(\overline{\eta})\re| =  0 \;\; \mbox{ uniformly for bounded $t$ }.
\end{equation}
Now (\ref{equationPhTtoPh}) follows due to the fact that $\mathbb{F}$ is continuous.
%%%%%%%%%%%%%%%%%%%%%%%%%%%%%%%%%%%%%%%%%%%%%%%%%%%%%%%%%%%%%%%%%%%%%%%%%%%%%%%%%%%%%%%%%%%%%%%%%%%      Proof of Proposition \ref{PropositionCAT} and its applications

\subsection{Proofs of Propositions \ref{PropositionCAT} and \ref{PropositionGenalogicalDistance}}
\noindent The transition rates of both the functional (\ref{equationReducedBPforCAT}) and the functional (\ref{equationReducedBPforDistances}) can be obtained from the transition rates of the  time-homogeneous transformed BP. Remember, see also Subsubsection \ref{subsubTransformeddualprocess}, that the time-homogeneous transformed BP has the same transitions (\ref{itemdualmutationJ} - \ref{itemdualresamplingIIv}) as the BP, but the rates are changed by the time-homogeneous potential $h$ defined in Theorem \ref{TheoremLongtime} and therefore by probabilities of the form (\ref{equationDefinitionPN}). 

We start with the functional (\ref{equationReducedBPforCAT}).  In this case there is only one partition element in the first component of the BP, where $u$ represents its mark in $K$ and $n$ the number of active life-sites with set $\{0\}$.
\begin{enumerate}
\item The transition $(u,n) \to (u,n+1)$ occurs if an active life-site with set $K$ becomes an active life-site with set $\{0\}$. On the one hand this can happen if either the interaction \ref{itemdualresamplingJIw} or the interaction \ref{itemdualresamplingIJw} occurs between the partition element and an active life-site with set $K$. On the other hand this can happen if the interaction \ref{itemdualresamplingIIw} occurs between an active life-site with set $\{0\}$ and an active life-site with set $K$.
\item The transition $(u,n) \to (u,n-1)$ occurs if an active life-site with set $\{0\}$ becomes an active life-site with set $K$. First, this can happen if the transition \ref{itemdualmutationIcup} occurs for an active life-site with set $\{0\}$. Second, this can happen if either the interaction \ref{itemdualresamplingJIK} or the interaction \ref{itemdualresamplingIJK} occurs between the partition element and an active life-site with set $\{0\}$ (observe that this is only possible if the mark of the partition element is $0$). Third,  this can happen if the interaction \ref{itemdualresamplingIIK} occurs between two active life-sites with set $\{0\}$.
\item The transition $(u,n) \to (1-u,n)$ occurs if the first component of the BP is changed by the transition in \ref{itemdualmutationJ}. 
\end{enumerate}
Hence the transition rates of the functional (\ref{equationReducedBPforCAT}) are composed as listed in Figure \ref{FigureReducedBPforDistances}, where $N-1-n$ represents the number of active life-sites with set $K$.

\renewcommand{\arraystretch}{1.5}
\begin{figure}[h]\caption{Transitions rates of the functional (\ref{equationReducedBPforCAT})}
\label{FigureReducedBPforCAT}
\begin{center}
\begin{tabular}{|l|l|}\hline
transition &  rate  \\ \hline 
%-------------------------------n to n+1---------------------------------
$(u,n) \to (u,n+1)$  & $[2(N-1-n)n\frac{S}{2N} + 2(N-1-n)\frac{S}{2N}] \frac{P_{N}(1^{u},0^{n+2-u})}{P_{N}(1^{u},0^{n+1-u})}$ \\ \hline
%-------------------------------n to n-1---------------------------------
$(u,n) \to (u,n-1)$ & $[Bb_{0}n + 2(1-u)n(\frac{1}{2} - \frac{S}{2N}) +n(n-1)(\frac{1}{2} - \frac{S}{2N})]\frac{P_{N}(1^{u},0^{n-u})}{P_{N}(1^{u},0^{n+1-u})}$ \\ \hline
%-------------------------------u to 1-u---------------------------------
$(u,n) \to (1-u,n)$ & $B[ub_{1}+ (1-u)b_{0}]\frac{P_{N}(1^{1-u},0^{n+u})}{P_{N}(1^{u},0^{n+1-u})}$ \\ \hline
\end{tabular} 
\end{center}
\end{figure}
%%%%%%%%%%%%%%%%%%%%%%%%%%%%%%%%%%%%%%%%%%%%%%%%%%%%%%%%%%%%%%%%%%%%%%%%%%%%%%%%%%%%%%%%%%%%%%%%%%%      Proof of Corollary \ref{Corollarydistances} and its applications

In analogy, the transition rates of the functional (\ref{equationReducedBPforDistances}) are composed of the transition rates of the  time-homogeneous transformed BP as listed in Figure \ref{FigureReducedBPforDistances}, where a transition to the absorbing state $\bigtriangleup$ can only occur if  the two partition elements in the first component of the BP coalesce, either by \ref{itemdualresamplingJJK} or by \ref{itemdualresamplingJJw}.

\renewcommand{\arraystretch}{1.5}
\begin{figure}[h]\caption{Transitions rates of the functional (\ref{equationReducedBPforDistances})}
\label{FigureReducedBPforDistances}
\begin{center}
\begin{tabular}{|l|l|}\hline
transition & the rate consists of\\ \hline
%-------------------------------n to n+1---------------------------------
$(\Circle,n) \to (\Circle,n+1)$ & $[2(N-2-n)n\frac{S}{2N}+2\cdot 2(N-2-n)\frac{S}{2N}] \frac{P_{N}(0^{n+3})}{P_{N}(0^{n+2})}$ \\ \hline
$(\CIRCLE,n) \to (\CIRCLE,n+1)$ & $[2(N-2-n)n\frac{S}{2N} + 2\cdot 2(N-2-n)\frac{S}{2N}]\frac{P_{N}(1^{2},0^{n+1})}{P_{N}(1^{2},0^{n})}$ \\ \hline
$(\RIGHTcircle,n) \to (\RIGHTcircle,n+1)$ & $[2(N-2-n)n\frac{S}{2N} + 2\cdot 2(N-2-n)\frac{S}{2N}]\frac{P_{N}(1,0^{n+2})}{P_{N}(1,0^{n+1})}$  \\ \hline
%-------------------------------n to n-1---------------------------------
$(\Circle,n) \to (\Circle,n-1)$ & $[Bb_{0}n + 2\cdot 2n(\frac{1}{2} - \frac{S}{2N}) +n(n-1)(\frac{1}{2} - \frac{S}{2N})]\frac{P_{N}(0^{n+1})}{P_{N}(0^{n+2})} $  \\ \hline
$(\CIRCLE,n) \to (\CIRCLE,n-1)$ & $[Bb_{0}n + 0\cdot 2n(\frac{1}{2} - \frac{S}{2N})+ n(n-1)(\frac{1}{2} - \frac{S}{2N})]\frac{P_{N}(1^{2},0^{n-1})}{P_{N}(1^{2},0^{n})} $  \\ \hline
$(\RIGHTcircle,n) \to (\RIGHTcircle,n-1)$ & $[Bb_{0}n + 1\cdot 2n(\frac{1}{2} - \frac{S}{2N}) +n(n-1)(\frac{1}{2} - \frac{S}{2N})]\frac{P_{N}(1,0^{n})}{P_{N}(1,0^{n+1})} $ \\ \hline
%--------------------- \Circle,\CIRCLE to \RIGHTcircle ------------------
$(\Circle,n),(\CIRCLE,n) \to (\RIGHTcircle,n)$ & $[Bb_{0}+Bb_{0}]\frac{P_{N}(1,0^{n+1})}{P_{N}(0^{n+2})},\;\; [Bb_{1}+Bb_{1}]\frac{P_{N}(1,0^{n+1})}{P_{N}(1^{2},0^{n})}$\\ \hline
%-------------------- \RIGHTcircle to \Circle, \CIRCLE ------------------
$(\RIGHTcircle,n) \to (\Circle,n),(\CIRCLE,n)$ & $Bb_{1}\frac{P_{N}(0^{n+2})}{P_{N}(1,0^{n+1})}, \;\; Bb_{0}\frac{P_{N}(1^{2},0^{n})}{P_{N}(1,0^{n+1})}$ \\ \hline
%------------------ \Circle, \CIRCLE to \bigtriangleup ------------------
$(\Circle,n), (\CIRCLE,n) \to \bigtriangleup$ & $2(\frac{1}{2} - \frac{S}{2N})\frac{P_{N}(0^{n+1})}{P_{N}(0^{n+2})} + 2\frac{S}{2N},\;\; 2\frac{1}{2}\frac{P_{N}(1,0^{n})}{P_{N}(1^{2},0^{n})} + 2\frac{S}{2N}\frac{P_{N}(1,0^{n+1})}{P_{N}(1^{2},0^{n})}$ \\ \hline
\end{tabular} 
\end{center}
\end{figure}
%%%%%%%%%%%%%%%%%%%%%%%%%%%%%%%%%%%%%%%%%%%%%%%%%%%%%%%%%%%%%%%%%%%%%%%%%%%%%%%%%%%%%%%%%%%%%%%%%%%      Proof of  Proposition \ref{Proposition}, S= 0

\subsection{Proof of Proposition \ref{Proposition}}
\noindent In order to analyze $f_{t}(y,n)$ for all $(y,n) \in \{\Circle, \CIRCLE, \RIGHTcircle\} \times \mathbb{N}_{0}$ as well as $pf_{t}$, we use that 
\begin{equation} 
\frac{\partial f_{t}(y,n)}{\partial t} = L^{\overline{Y}^{2}}f_{t}(y,n) \;\; \mbox{ for all }\; t \ge 0\;,
\end{equation}
where $L^{\overline{Y}^{2}}$ is the generator corresponding to the jump process $\overline{Y}^{2}$. Thus, it is not hard to see that
\begin{equation}\label{equationProofPropHelp1}
\li.\frac{\partial f_{t}(\Circle,n)}{\partial t}\re|_{t = 0}  =  -\frac{E[0^{n+1}]}{E[0^{n+2}]}, \; \li.\frac{\partial f_{t}(\CIRCLE,n)}{\partial t}\re|_{t = 0}  =  -\frac{E[1,0^{n}]}{E[1^{2},0^{n}]} \mbox{ and } \li.\frac{\partial f_{t}(\RIGHTcircle,n)}{\partial t}\re|_{t = 0}  =  0
\end{equation}
which shows (set $n=0$) the statements concerning $f_{t}(\Circle,0)$, $f_{t}(\CIRCLE,0)$ and $f_{t}(\RIGHTcircle,0)$ in the case with selection.

In addition, for $n \in \mathbb{N}_{0}$ we set
\begin{equation}
pf_{t}(n) := E[0^{n+2}]f_{t}(\Circle,n) + E[1^{2},0^{n}]f_{t}(\CIRCLE,n) + 2E[1,0^{n+1}]f_{t}(\RIGHTcircle,n) \;\; \mbox{ for all } t \ge 0\;,
\end{equation}
that is, $pf_{t} = pf_{t}(0)$ and 
\begin{equation}\label{equationProofPropHelp2}
 pf_{0}(n)  = E[0^{n}] \;\;\mbox{ as well as } \;\; \li. \frac{\partial pf_{t}(n)}{\partial t}\re|_{t = 0}  = -E[0^{n}]
\end{equation}
for each $n \in \mathbb{N}_{0}$.
%%%%%%%%%%%%%%%%%%%%%%%%%%%%%%%%%%%%%%%%%%%%%%%%%%%%%%%%%%%%%%%%%%%%%%%%%%%%%%%%%%%%%%%%%%%%%%%%%%%      Proof of  Proposition \ref{Proposition}, Lemma mixed moments

We start with a lemma in which we derive a system of ODE's for the functions $f_{t}(y,n)$, $(y,n) \in \{\Circle, \CIRCLE, \RIGHTcircle\} \times \mathbb{N}_{0}$, and the functions  $pf_{t}(n)$, $n \in \mathbb{N}_{0}$. This lemma is used to show the remaining statements of Proposition \ref{Proposition} and additionally provides a starting point for a proof that genealogical distances are stochastically smaller under selection.
\begin{Lemma}\label{LemmaMixedMoments}\point  For each $n \in \mathbb{N}_{0}$ and each $t \ge 0$ we have that
\begin{eqnarray}
E[0^{n+2}]\tfrac{\partial f_{t}(\Circle,n)}{\partial t} &=&  (nBb_{0}+\tfrac{(n+2)(n+1)}{2} - 1)E[0^{n+1}]f_{t}(\Circle,n-1) \\
&&+ 2Bb_{0}E[1,0^{n+1}]f_{t}(\RIGHTcircle,n)+ (2+n)SE[0^{n+3}]f_{t}(\Circle,n+1) \\
&&- (\tfrac{(2+n)(n+1+2B+2S)}{2} -2B +2Bb_{1})E[0^{n+2}]f_{t}(\Circle,n),
\end{eqnarray}
\begin{eqnarray}
E[1^{2},0^{n}]\tfrac{\partial f_{t}(\CIRCLE,n)}{\partial t}&=&(nBb_{0}+\tfrac{n(n-1)}{2})E[1^{2},0^{n-1}]f_{t}(\CIRCLE,n-1) \\
&& + 2Bb_{1}E[1,0^{n+1}]f_{t}(\RIGHTcircle,n) + (n+2)SE[1^{2},0^{n+1}]f_{t}(\CIRCLE,n+1)\;\;\\
&&- (\tfrac{(2+n)(n+1+2B+2S)}{2}- 2S-2B + 2Bb_{0})E[1^{2},0^{n}]f_{t}(\CIRCLE,n)
\end{eqnarray}
and 
\begin{eqnarray}
E[1,0^{n+1}]\tfrac{\partial f_{t}(\RIGHTcircle,n)}{\partial t} &=& (nBb_{0}+\tfrac{n(n+1)}{2})E[1, 0^{n}]f_{t}(\RIGHTcircle,n-1) + Bb_{1}E[0^{n+2}]f_{t}(\Circle,n) \;\\
&& Bb_{0}E[1^{2},0^{n}]f_{t}(\CIRCLE,n) + (n+2)SE[1,0^{n+2}]f_{t}(\RIGHTcircle,n+1)\\
&&- (\tfrac{(2+n)(n+1+2B+2S)}{2} -S -B)E[1,0^{n+1}]f_{t}(\RIGHTcircle,n)\;.
\end{eqnarray}
This means
\begin{eqnarray}
\tfrac{\partial pf_{t}(n)}{\partial t} &=& \tfrac{n}{2}(n-1+2Bb_{0})pf_{t}(n-1) -[\tfrac{(n+2)}{2}(n+1+2B+2S)-2B]pf_{t}(n)\\ 
&&+(n+2)Spf_{t}(n+1) + \mathcal{R}_{t}(n)
\end{eqnarray}
for all $t \ge 0$,  where
\begin{eqnarray}
\mathcal{R}_{t}(n) &=&  2nE[0^{n+1}]f_{t}(\Circle,n-1) + n2E[1,0^{n}]f_{t}(\RIGHTcircle,n-1)\\
&& + 2SE[1^{2},0^{n}]f_{t}(\CIRCLE,n) + S2E[1,0^{n+1}]f_{t}(\RIGHTcircle,n)\;.
\end{eqnarray}
\end{Lemma}
%%%%%%%%%%%%%%%%%%%%%%%%%%%%%%%%%%%%%%%%%%%%%%%%%%%%%%%%%%%%%%%%%%%%%%%%%%%%%%%%%%%%%%%%%%%%%%%%%%%      Proof of Proposition, Proof Lemma mixed moments

\noindent {\bf Proof of Lemma \ref{LemmaMixedMoments}:} Consider the following three relations which hold, see for example \cite{F02}, for the mixed moments  $E[0^{n+2}]$, $E[1,0^{n+1}]$ and $E[1^{2},0^{n}]$ defined in (\ref{equationDefinitionMixedMoments}). For each $n \in \mathbb{N}_{0}$ one has that 
\begin{equation}
(n+1+2B+2S)E[0^{n+2}] = (n+1+2Bb_{0})E[0^{n+1}] + 2SE[0^{n+3}]\;,
\end{equation}
\begin{eqnarray}
&&([n+2][n+1+2B+2S]-2S)E[1,0^{n+1}]\\
&=& (n+1)(n+2Bb_{0})E[1,0^{n}] + 2Bb_{1}E[0^{n+1}] + (n+2)2SE[1,0^{n+2}]
\end{eqnarray}
and
\begin{eqnarray}
&&([n+2][n+1+2B+2S]-4S)E[1^{2},0^{n}] \\
&=&n(n+2Bb_{0}-1)E[1^{2},0^{n-1}] + 2(1+2Bb_{1})E[1,0^{n}] + (n+2)2SE[1^{2},0^{n+1}].\;
\end{eqnarray}
Using these relations and the generator of $\overline{Y}^{2}$ one obtains the stated equations. \hfill $\blacksquare$ \bigskip

Now we come to the proof of the statements in Proposition \ref{Proposition}, where we distinguish between $S = 0$  (without selection) and $S >0$ (with selection).
\bigskip
%%%%%%%%%%%%%%%%%%%%%%%%%%%%%%%%%%%%%%%%%%%%%%%%%%%%%%%%%%%%%%%%%%%%%%%%%%%%%%%%%%%%%%%%%%%%%%%%%%%      Proof of Proposition, S = 0

{\it The case $S = 0$:} In this case $\overline{Y}^{2}$ is a jump process on $\{(\Circle,0), (\CIRCLE,0), (\RIGHTcircle,0), \bigtriangleup \}$ and
\begin{itemize}
\item[] $E[0] = b_{0}$, $E[1] = b_{1}$, $E[0^{2}] = \frac{(2Bb_{0} +1)b_{0}}{2B +1}$, $E[1^{2}] = \frac{(2Bb_{1} +1)b_{1}}{2B+1}$ and $2E[1,0] = 2  \frac{2Bb_{0}b_{1}}{2B+1}$.
\end{itemize}
Thus, using these explicit values together with Lemmma \ref{LemmaMixedMoments} for $n = 0$, we obtain the following system of ODE'S:
\begin{enumerate}
\item $\tfrac{\partial f_{t}(\Circle,0)}{\partial t} =\frac{4B^{2}b_{1}b_{0}}{1+2Bb_{0}}f_{t}(\RIGHTcircle,0)  - (1+2Bb_{1})f_{t}(\Circle,0) $
\item $\tfrac{\partial f_{t}(\CIRCLE,0)}{\partial t} = \frac{4B^{2}b_{1}b_{0}}{1+2Bb_{1}}f_{t}(\RIGHTcircle,0)  -(1+2Bb_{0})f_{t}(\CIRCLE,0)$
\item $\tfrac{\partial f_{t}(\RIGHTcircle,0)}{\partial t} = (1+2Bb_{0})\frac{1}{2}f_{t}(\Circle,0) + (1+2Bb_{1})\frac{1}{2}f_{t}(\CIRCLE,0)   - (1+B)f_{t}(\RIGHTcircle,0)$
\end{enumerate}
Solving this system we get that
\begin{enumerate}
\item $f_{t}(\Circle,0) =\frac{(1+2B)b_{0}e^{-t}}{(1+2Bb_{0})}  + \frac{b_{1}e^{-(1+2B)t}}{(1+2Bb_{0})} = e^{-t}(1 + \frac{b_{1}(e^{-2Bt} -1)}{(1+2Bb_{0})})$
\item $f_{t}(\CIRCLE,0) = \frac{(1+2B)b_{1}e^{-t}}{(1+2Bb_{1})}  + \frac{b_{0}e^{-(1+2B)t}}{(1+2Bb_{1})} = e^{-t}(1 + \frac{b_{0}(e^{-2Bt} -1)}{(1+2Bb_{1})})$
\item $f_{t}(\RIGHTcircle,0) = \frac{(1+2B)}{2B}e^{-t} - \frac{1}{2B}e^{-(1+2B)t} = e^{-t}(1 - \frac{(e^{-2Bt}-1)}{2B}) $
\end{enumerate}
which is the desired result.
%%%%%%%%%%%%%%%%%%%%%%%%%%%%%%%%%%%%%%%%%%%%%%%%%%%%%%%%%%%%%%%%%%%%%%%%%%%%%%%%%%%%%%%%%%%%%%%%%%%      Proof of  Proposition \ref{Proposition}, S > 0

\bigskip
{\it The case $S >0$:} The plan is to expand $pf_{t}(0)$ into its Taylor series at $0$ in order to compare it with $e^{-t}$. We do this expansion up to  degree $3$, where we use Lemma \ref{LemmaMixedMoments} as well as (\ref{equationProofPropHelp1}) and (\ref{equationProofPropHelp2}).

Degree $1$: Since
\begin{equation}
\tfrac{\partial pf_{t}(0)}{\partial t} = -(1+2S)pf_{t}(0) +2S[pf_{t}(1) + E[1^{2}]f_{t}(\CIRCLE,0) +E[1,0]f_{t}(\RIGHTcircle,0)]\;,
\end{equation}
we have that 
\begin{equation}
\li. \tfrac{\partial pf_{t}(0)}{\partial t}\re|_{t = 0} = -(1+2S) +2S(E[0] + E[1^{2}] +E[1,0]) = -1\;.
\end{equation}
%%%%%%%%%%%%%%%%%%%%%%%%%%%%%%%%%%%%%%%%%%%%%%%%%%%%%%%%%%%%%%%%%%%%%%%%%%%%%%%%%%%%%%%%%%%%%%%%%%%      Proof of  Proposition \ref{Proposition}, S > 0, Taylor, degree 2

Degree $2$: We have that  
\begin{equation}
\li.\tfrac{\partial^{2}pf_{t}(0)}{\partial t^{2}}\re|_{t = 0} =  -(1+2S)(-1) +2S(-E[0] -E[1] +0) = 1  \;.
\end{equation}
%%%%%%%%%%%%%%%%%%%%%%%%%%%%%%%%%%%%%%%%%%%%%%%%%%%%%%%%%%%%%%%%%%%%%%%%%%%%%%%%%%%%%%%%%%%%%%%%%%%      Proof of  Proposition \ref{Proposition}, S > 0, Taylor, degree 2

Degree $3:$  First we have that
\begin{equation}
\li.\tfrac{\partial^{3}pf_{t}(0)}{\partial t^{3}}\re|_{t = 0} =  -(1+2S) +2S(\tfrac{\partial^{2}pf_{t}(1)}{\partial t^{2}}|_{t = 0} + E[1^{2}] \tfrac{\partial^{2}f_{t}(\CIRCLE,0)}{\partial t^{2}}|_{t = 0} + E[1,0] \tfrac{\partial^{2}f_{t}(\RIGHTcircle,0)}{\partial t^{2}}|_{t = 0}).
\end{equation}
Since
\begin{eqnarray}
\tfrac{\partial pf_{t}(1)}{\partial t} &=& Bb_{0}pf_{t}(0) - (3+B+3S) pf_{t}(1) +3S pf_{t}(2) + 2E[0^{2}]f_{t}(\Circle,0) \\
&&  + 2E[1,0]f_{t}(\RIGHTcircle,0) + 2SE[1^{2},0]f_{t}(\CIRCLE,1) + S2E[1,0^{2}]f_{t}(\RIGHTcircle,1)
\end{eqnarray}
we have that 
\begin{equation}
\li. \tfrac{\partial^{2} pf_{t}(1)}{\partial t^{2}}\re|_{t = 0} = (3+B+3S)E[0] -Bb_{0} -3SE[0^{2}] - 2E[0] - 2SE[1,0]\;.
\end{equation}
Since
\begin{eqnarray}
E[1^{2}]\tfrac{\partial f_{t}(\CIRCLE,0)}{\partial t}  &=& 2Bb_{1}E[1,0]f_{t}(\RIGHTcircle,0) +2SE[1^{2},0]f_{t}(\CIRCLE,1)\\
&& - (1+ 2Bb_{0})E[1^{2}]f_{t}(\CIRCLE,0)
\end{eqnarray}
and since
\begin{eqnarray}
E[1,0] \tfrac{\partial f_{t}(\RIGHTcircle,0)}{\partial t} &=& Bb_{1}E[0^{2}]f_{t}(\Circle,0) + Bb_{0}E[1^{2}] f_{t}(\CIRCLE,0) \\
&& + 2SE[1,0^{2}]f_{t}(\RIGHTcircle,1) - (1+B+S)E[1,0]f_{t}(\RIGHTcircle,0)\;,
\end{eqnarray}
we have that
\begin{equation}
\li. E[1^{2}]\tfrac{\partial^{2} f_{t}(\CIRCLE,0)}{\partial t^{2}}\re|_{t = 0} + \li. E[1,0] \tfrac{\partial^{2} f_{t}(\RIGHTcircle,0)}{\partial t^{2}}\re|_{t = 0}  =  (1+ Bb_{0})E[1] - 2SE[1,0] -Bb_{1}E[0]\,.
\end{equation}
Hence
\begin{eqnarray}
\li.\tfrac{\partial^{3}pf_{t}(0)}{\partial t^{3}}\re|_{t = 0} &=& -1-2S +2S([1+B+3S]
E[0]-Bb_{0} -3SE[0^{2}] -2SE[1,0]) \\
&& + 2S([1+Bb_{0}]E[1] - Bb_{1}E[0]-2SE[1,0])\\
&=& -1-2S +2S(1+ 3SE[0]-3SE[0^{2}] -4SE[1,0])\\
&=& -1 +2S^{2}(3E[0]-3(E[0^{2}] +E[1,0]) -E[1,0])\\
&=& -1 -2S^{2}E[1,0] < -1
\end{eqnarray}
and the proof is finished.
%%%%%%%%%%%%%%%%%%%%%%%%%%%%%%%%%%%%%%%%%%%%%%%%%%%%%%%%%%%%%%%%%%%%%%%%%%%%%%%%%%%%%%%%%%%%%%%%%%% literature

\bibliography{literatur}
\addcontentsline{toc}{section}{References}
\end{document}